\documentclass[hidelinks,onefignum,onetabnum]{siamart251216}

\usepackage{lipsum}
\usepackage{amsfonts}
\usepackage{graphicx}
\usepackage{epstopdf}
\usepackage{algorithmic}
\ifpdf
  \DeclareGraphicsExtensions{.eps,.pdf,.png,.jpg}
\else
  \DeclareGraphicsExtensions{.eps}
\fi

\newcommand{\bu}{{\bf u}}
\newcommand{\bw}{{\bf w}}

\newsiamremark{remark}{Remark}
\newsiamremark{hypothesis}{Hypothesis}
\crefname{hypothesis}{Hypothesis}{Hypotheses}
\newsiamthm{claim}{Claim}
\newsiamremark{fact}{Fact}
\crefname{fact}{Fact}{Facts}

\headers{High-order in-cell discontinuous reconstruction path-conservative methods for conservative systems in primitive variables}{E. Pimentel-García, M. Ricchiuto, and C. Parés}

\title{High-order in-cell discontinuous reconstruction path-conservative methods for conservative systems in primitive variables\thanks{Submitted to the editors 14/03/2026.
\funding{This work was funded by PID2022-137637NB-C21, PCI2024-155061-2 and PPRO-FQM216-2023}}}

\author{Ernesto Pimentel-García\thanks{Departamento Matemática Aplicada, Universidad de Málaga, C/ Dr. Ortiz Ramos s/n. Campus de Teatinos. 29071, Málaga, Spain 
  (\email{erpigar@uma.es}, \url{www.erpigar.es}).}
\and Mario Ricchiuto\thanks{Inria, Univ. Bordeaux, CNRS, Bordeaux INP, IMB, UMR 5251, 200 Avenue de la Vieille Tour, 33405 Talence cedex, France
  (\email{mario.ricchiuto@inria.fr}).}
\and Carlos Parés\thanks{Departamento  de  An\'alisis  Matem\'atico,  Estad\'istica  e  Investigaci\'on  Operativa,  y  Matem\'atica  aplicada,  Universidad  de  M\'alaga,Bulevar Louis Pasteur, 31, 29010, M\'alaga, Spain. 
  (\email{pares@uma.es}).}}

\usepackage{amsopn}

\ifpdf
\hypersetup{
  pdftitle={High-order in-cell discontinuous reconstruction path-conservative methods for conservative systems in primitive variables},
  pdfauthor={E. Pimentel-García, M. Ricchiuto, and C. Parés}
}
\fi

\begin{document}

\maketitle

\begin{abstract}
We present a path-conservative methodology for the numerical approximation of conservative systems formulated exclusively in primitive variables, the most physically relevant set of variables. The approach is designed to preserve the correct weak solution while avoiding the spurious oscillations near contact discontinuities that may arise in conservative formulations. Nonconservative products are treated through the selection of a suitable family of paths combined with an in-cell discontinuous reconstruction strategy, which is locally activated near shocks. To choose the path we introduce a  compatibility condition between the Rankine-Hugoniot conditions 
of the non-conservative form, with the exact ones, assumed to be known. 
In smooth regions and across contact discontinuities, a standard path-conservative scheme is employed. The method is further extended to achieve high-order accuracy. Numerical results for the Euler multimaterial and shear interface model demonstrates correct shock-capturing property, suppression of spurious oscillations due to contact discontinuities, and convergence without using hybrid formulations.
\end{abstract}



\section{Introduction: conservative and non-conservative forms}

Let us consider a PDE system of conservation laws of the form:

\begin{equation}\label{conservative_system}
    \bu_{t} + F(\bu)_{x} = 0,
\end{equation}
where $\bu(x,t)$ is the conservative vector of unknowns that takes values in the convex set $\Omega$ of $\mathbb{R}^{N}$ and $F(\bu)$ is the physical flux that goes from $\Omega$ to $\mathbb{R}^{N}$.
In this work, we focus on a general family of systems in which the components of the vector $\bu$ are mass densities, momenta, and energies, plus additional variables characterizing the local
composition of the flow. This type of structure is encountered very often both in multi-phase or multi-material flows, as well as in  many shallow water  models.

The main motivation of this work is to design numerical methods
able to solve certain fluid models using the so-called primitive variables $\bw$ instead of the conserved ones $\bu$ without altering shock waves. In some cases, the use of primitive variables is mandatory due to the spurious oscillations that appear close to contact discontinuities when solving the problem in conserved variables (see, for example,  \cite{Quirk1996129}).

If the variables $\bu$ and $\bw$ are related by a smooth one-to-one transformation  $T(\bw) = \bu$, the system in the new variables is written as follows:
\begin{equation}\label{nonconservative_system}
    \bw_{t} + \mathcal{A}(\bw)\bw_{x} = 0
\end{equation}
where
$$\mathcal{A}(\bw) = \frac{\partial \bw}{\partial \bu}\frac{\partial F}{\partial \bu}\frac{\partial \bu}{\partial \bw},$$
having denoted with  ${\partial F}/{\partial \bu}$, ${\partial \bw}/{\partial \bu}$, ${\partial \bu}/{\partial \bw}$ the Jacobians of $F$, $T$,   and $T^{-1}$ respectively. Therefore the system of conservation laws  \eqref{conservative_system} is equivalent to  the nonconservative system \eqref{nonconservative_system} for smooth solutions.

The objective of this work is to solve \eqref{conservative_system} by applying an adequate numerical method to \eqref{nonconservative_system}  so that no oscillations occur at the contact discontinuities and the shock waves are captured correctly. To do this, path-conservative numerical methods combined with the so-called in-cell discontinuous reconstruction technique (see \cite{chalons2019path, pimentel2022cell, pimentel2024high}) will be used. 

Previous works are based on a hybrid approach (see, for instance, \cite{karni1996hybrid,  abgrall2001computations, GABURRO2024117311}), in which primitive variables are used to handle smooth regions (rarefactions) and contact discontinuities, while conservative variables are employed when shocks appear in the solution. The main advantage of the present work is that the use of conservative variables is no longer required.

A first difficulty to overcome stems from the definition of weak solution of \eqref{nonconservative_system}. Various mathematical frameworks provide ways to define the nonconservative product \(\mathcal{A}(\bw)\,\partial_x \bw\) for discontinuous functions $\bw(x)$. Here, the approach by Dal Maso, LeFloch, and Murat \cite{Maso1995} is applied, introducing a weak solution concept that satisfies \eqref{nonconservative_system} in the sense of Borel measures. This approach relies on a family of Lipschitz continuous paths, \(\Phi: [0,1] \times \Omega \times \Omega \to \Omega\), which must satisfy certain regularity and compatibility conditions:
\[
\Phi(0; \bw_l, \bw_r) = \bw_l, \quad \Phi(1; \bw_l, \bw_r) = \bw_r, \quad \forall \bw_l, \bw_r \in \Omega,
\]
\[
\Phi(s; \bw, \bw) = \bw, \quad \forall \bw \in \Omega.
\]

Let  $\bw$ be a function with bounded total variation. According to the definition of the nonconservative product \(\mathcal{A}(\bw)\,\partial_x \bw\) as a Borel measure, the total measure of an interval $(a,b)$, which will de denoted by
\[
\dashint_a^b \mathcal{A}(\bw(x))\,\partial_x \bw(x) \, dx,
\]
is given by     
\begin{equation*}
\begin{aligned}
\dashint_a^b \mathcal{A}(\bw(x))\,\partial_x \bw(x) \, dx
&= \int_a^b \mathcal{A}(\bw(x))\,\partial_x \bw(x) \, dx  \\
&\quad + \sum_m \int_0^1 \mathcal{A}(\Phi(s; \bw_m^-, \bw_m^+))
      \frac{\partial \Phi}{\partial s}(s; \bw_m^-, \bw_m^+) \, ds.
\end{aligned}
\end{equation*}
where $\partial_x \mathbf{w}$ in the integral on the right-hand side represents the derivative at almost every point, and \(\bw_m^-\) and \(\bw_m^+\) represent the values to the left and right of the $m$th jump discontinuity.  Observe that the family of paths is used to determine the weights of the Dirac masses that accounts for the derivatives of  $\mathbf{w}$ at the jumps.

If nonconservative products are defined in this way, the admissible weak solutions of System \eqref{nonconservative_system} have to satisfy the  generalized Rankine-Hugoniot condition:
\begin{equation}\label{gen_R-H}
\int_0^1 \mathcal{A}(\Phi(s; \bw^-, \bw^+)) \frac{\partial \Phi}{\partial s}(s; \bw^-, \bw^+) \, ds = \sigma (\bw^+ - \bw^-),
\end{equation}
where \(\sigma\) is the speed of propagation of the discontinuity, and \(\bw^-\) and \(\bw^+\) are the values on either side. The notion of weak solutions depends thus on the choice of family of paths.
Therefore, this family has to be adequately chosen so that the jump condition
\eqref{gen_R-H} is equivalent to the  Rankine-Hugoniot condition for \eqref{conservative_system}:
\begin{equation}\label{con_R-H}
F(\bu^{+})-F(\bu^{-}) = \sigma (\bu^+ - \bu^-).
\end{equation}
Although the general question of whether there exists a family of paths that ensures the equivalence of \eqref{gen_R-H} and \eqref{con_R-H} is complex and involves the study of viscous profiles (see Appendix \ref{A1}),  in some cases, as we shall see in Section \ref{sec:Euler}, such a family of paths can be calculated algebraically.

After prescribing the family of paths, defining a concept of entropy solution becomes necessary: the notion of entropy induced by the one of the conservative system (that can be given by an entropy pair or by Lax's condition, for instance) will be assumed here, i.e. a weak solution $\bw(x,t)$ of \eqref{nonconservative_system} will be considered an entropy solution if $\bu(x,t) = T(\bw(x,t))$ is an entropy solution for \eqref{conservative_system}.

Once the family of paths and the notion of entropy have been prescribed so that the weak solutions of \eqref{conservative_system} and \eqref{nonconservative_system} are the same, defining numerical methods for the latter that  converge to the correct weak solutions is not a simple task.  The framework of path-conservative methods introduced in \cite{pares2006numerical} makes easy the design of numerical methods that are consistent with the definition of weak-solutions  but,  if the numerical dissipation is not controlled, the limits of the numerical solutions can be weak solutions satisfying a different jump condition: see \cite{Castro2008}, \cite{Abgrall2010}. The in-cell reconstruction technique will be used here to control the numerical dissipation so that isolated shock waves are correctly captured and \eqref{con_R-H} is satisfied: see \cite{chalons2019path, pimentel2022cell, pimentel2024high}.

The paper is organized as follows. In Section \ref{sec:2} we present the general methodology developed in this work. We first introduce a path-conservative semi-discrete formulation of the method, detailing the incorporation of the in-cell reconstruction operator, to finally obtain the fully discrete method. In Section \ref{sec:Euler} we focus on the Euler equations with an interface and shear, providing all the necessary ingredients to apply the proposed method. Section \ref{sec:4} demonstrates the numerical performance of the method, comparing it with existing approaches to highlight the improvements achieved. Finally, Section \ref{sec:5} includes some final conclusions and remarks.

\section{High-order in-cell discontinuous reconstruction path-conservative methods for conservative systems in primitive variables}\label{sec:2}

Following \cite{pares2006numerical}, high-order path-conservative numerical schemes in semi-discrete form:
\begin{equation}\label{eq:hopc}
    {\bf w}_{i}'(t) = -\frac{1}{\Delta x}\left(\mathcal{D}_{i+\frac{1}{2}}^{-}(t) + \mathcal{D}_{i-\frac{1}{2}}^{+}(t) + \dashint_{x_{i-\frac{1}{2}}}^{x_{i+\frac{1}{2}}} \mathcal{A}(P_{i}^{n}(x,t))\frac{\partial}{\partial x} P_{i}^{n}(x,t)dx\right),\quad t\geq t_n
\end{equation}
will be used to solve \eqref{nonconservative_system},
where:
\begin{itemize}
    \item $\Delta x$ and $\Delta t$ are the space and the time steps respectively. They are assumed to be constant in the description of the numerical method for simplicity.
    \item $I_{i} = \left[x_{i-\frac{1}{2}}, x_{i-\frac{1}{2}}\right]$ are the computational cells with lenght $\Delta x$.
    \item $t_n = n\Delta t$, $n=0,1,2,...$.
	\item $\bw_{i}(t) \cong \displaystyle \frac{1}{\Delta x}\int_{x_{i+\frac{1}{2}}}^{x_{i-\frac{1}{2}}} \bw(x,t) \,dx$ is the approximation to the cell average of the solution at the $i$-th cell in time $t$;
    \item $P_{i}^{n}(x,t)$ is a space-time reconstruction operator, which may be discontinuous, to be defined in Subsection \ref{sec:reconstruction_operator}. As it will be seen, it is designed so that the method is high-order accurate in smooth regions while suppressing viscosity near shocks, thereby ensuring the correct jump conditions are recovered;
	\item $\mathcal{D}^\pm_{i+1/2}(t) =\mathcal{D}^\pm(P^n_{i}(x_{i+1/2}, t),P^n_{i+1}(x_{i+1/2}, t))$, where $\mathcal{D}^-$ and $\mathcal{D}^+$ are two Lipschitz continuous functions from $\Omega \times \Omega$ to $\Omega$ that satisfy
\begin{equation}\label{pc1}
\mathcal{D}^{\pm}(\bw,\bw)=0, \quad \forall \bw\in \Omega,
\end{equation}
and
\begin{equation}\label{pc2}
\mathcal{D}^-(\bw_l,\bw_r) + \mathcal{D}^+(\bw_l, \bw_r) = \int _0^1 \mathcal{A}\big(\Phi (s;\bw_l,\bw_r)\big)\frac{\partial
\Phi }{\partial s} (s;\bw_l,\bw_r)\,ds,
\end{equation}
for every  set $\bw_{l}, \bw_r \in \Omega$. See \cite{castro2017well} for details.
 \end{itemize}

Once the numerical approximations $\bw_i^n$ at time $t^n = n \Delta t$ have been computed, the ODE system \eqref{eq:hopc} has to be solved with an adequate time stepping to update the numerical solution.

\subsection{Definition of reconstruction operators $P_{i}^{n}(x,t)$}\label{sec:reconstruction_operator}
The first step is to identify the cells $I_i$ where a shock wave is present or may be generated. 
To do this, the Riemann problem with initial data 
\begin{equation}\label{eq:Riemann}
\bw(x,0) =
\begin{cases}
\bw_{i-i}^n, & x < 0, \\
\bw_{i+1}^n, & x > 0,
\end{cases}
\end{equation}
is used as a detector: if its solution involves at least one  shock wave, the cell is marked. If the only wave detected is a shock,  it will be said that the shock wave is isolated.
Let us denote by $\mathcal{M}_n$ the set of indices corresponding to the marked cells, namely
\begin{equation}
\mathcal{M}_n = \left\{\,
\begin{aligned}
  & i \ \ \text{s.t. the solution of the Riemann problem \eqref{nonconservative_system},} \\
  & \text{\eqref{eq:Riemann} contains one or more shock waves}
\end{aligned}
\,\right\}.
\end{equation}

The second step is to select a space-time reconstruction operator $Q_i$ that, given the cell averages $\bw^n_i$ at the cells  of a smooth function $\bw(x,t)$ at time $t_n$, provides smooth approximations of $\bw$ at the cells 
$$
\bw(x,t) \approx Q_i(x,t) = Q_i(x,t; \{ \bw^n_j \}_{j \in \mathcal{S}_{i} }), \quad x \in [x_{i-1/2}, x_{i+1/2}], \ t \in [t_n, t_{n+1}],
$$
where $\mathcal{S}_{i}$ is the stencil, i.e. the set of indices of the cells whose averages are used to compute $Q_i$. In the case of first order methods, the piecewise constant reconstruction operator will be considered. For higher order methods,  Taylor expansions and the Cauchy-Kowaleskaya procedure will be followed to compute these reconstruction operators as in the ADER approach \cite{titarev2005ader, dumbser2009ader}.

On the basis of the marked cells and the standard reconstruction operator, we define $P_{i}^{n}(x,t)$ as follows:
\begin{itemize}

\item If $i \in \mathcal{M}_n$ then the reconstruction operator is defined by
        \begin{equation}\label{eq:in_cell_reconstruction}
        P_{i}^{n}(x,t) = 
\begin{cases}
\bw^n_{i,l} & \text{ if $x \leq x_{i - 1/2} + d_{i}^{n} \Delta x + \sigma_i^n (t -t^n) $,}\\
\bw^n_{i,r} & \text{otherwise},
\end{cases}.
\end{equation}
where 
\begin{itemize}

\item  $\bw^n_{i,l}$, $\bw^n_{i,r}$ are two states and $\sigma_i^n$ a scalar to be chosen so that if $\bw^n_{i-1}$ and
$\bw^n_{i+1}$ can be  linked by an admissible shock with speed $\sigma$, then
\begin{equation}\label{speedstates}
\bw^n_{i,l} = \bw^n_{i-1}, \quad \bw^n_{i,r} = \bw^n_{i+1}, \quad \sigma_i^n = \sigma.
\end{equation}

\item $d_i^n$ is chosen so that
\begin{equation}\label{conservation}
d_i^n w_{i,l, k}^n + (1 -  d_i^n) w_{i,r, k}^n = w_{i, k}^n, 
\end{equation}
where $k \in \{ 1, \dots, N \}$ is a predetermined index (in the case of fluids, the one corresponding to the mass conservation equation).

Observe that this in-cell discontinuous reconstruction can only be done
if $0 \leq d_i^n \leq 1$, i.e. if
\begin{equation}\label{eq:d_property}
    0 \leq \frac{w_{i,r,k}^n - w_{i,k}^n}{w_{i,r,k}^n - w_{i,l,k}^n} \leq 1,
\end{equation}
otherwise the index $i$ is removed from the set $\mathcal{M}_n$.

\end{itemize}
\begin{remark}
The idea behind this reconstruction operator is the following (see Figure \ref{fig:In_cell_discontinuous_reconstruction_operator}):
\begin{itemize}
\item Two intermediate states of the solution of the Riemann problem \eqref{eq:Riemann} (or approximations of them), $\bw^n_{i,l}$, $\bw^n_{i,r}$, linked by an admissible shock with speed $\sigma_i^n$ are selected.
\item At time $t^n$, the constant state $\bw_i^n$ is replaced in the cell $I_i$ by a discontinuity placed at $x_{i-1/2} + d_i^n \Delta x$ linking $\bw^n_{i,l}$ and $\bw^n_{i,r}$,
where $d_i^n$ is chosen so that the mass is conserved.
\item The reconstruction operator $P_i^n$ represents the displacement at velocity $\sigma_i^n$ of this discontinuity in the cell $I_i$.
\end{itemize}
The states $\bw^n_{i,l}$, $\bw^n_{i,r}$ are such that, if $\bw_{i-1}^n$ and $\bw_{i+1}^n$ can be linked by an admissible shock traveling at speed $\sigma$, the reconstruction coincides with the exact solution of the Riemann problem. 
\end{remark}

\item If $i \notin \mathcal{M}_n$ and $\mathcal{S}_{i}\cap  \mathcal{M}_n \neq \emptyset$, then 
    \begin{equation}\label{eq:Constant_reconstruction_operator}
        P_{i}^{n}(x,t)= \bw_{i}^{n}.
    \end{equation}
\begin{remark}
The idea behind this reconstruction operator is the following (see Figure \ref{fig:Constant_reconstruction_operator}): suppose that $\bw_i^n$ and $\bw_{i+2}^n$ can be linked by an admissible shock. Then, $I_{i+1}$ will be a marked cell and the reconstructions at $I_i$, $I_{i+2}$ will be the constant states $\bw_i^n$ and $\bw_{i+2}^n$ respectively: no wave interaction occurs at the intercells then and the reconstructions at the three cells coincide with the solution of the Riemann problem. 
\end{remark}
\item Otherwise (i.e.,  if $\mathcal{S}_{i}\cap  \mathcal{M}_{n} = \emptyset$)
$$
P_i^n(x,t) = Q_i^n(x,t; \{ \bw_j^n \}_{j \in \mathcal{S}_i} ),
$$
i.e. the standard reconstruction operator is considered in the cells where no shocks are detected.

\begin{figure}[htpb]
    \centering
    \subfloat[In-cell discontinuous reconstruction operator at time $t = t^{n}$.]{
        \includegraphics[width=0.45\linewidth]{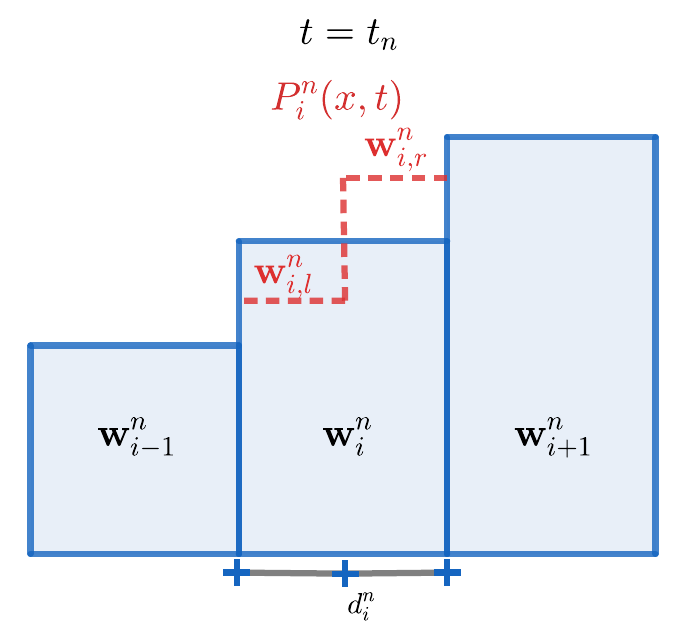}
    }
    \hfill
    \subfloat[In-cell discontinuous reconstruction operator at intermediate time $t^{n} < t < t^{n+1}$.]{
        \includegraphics[width=0.45\linewidth]{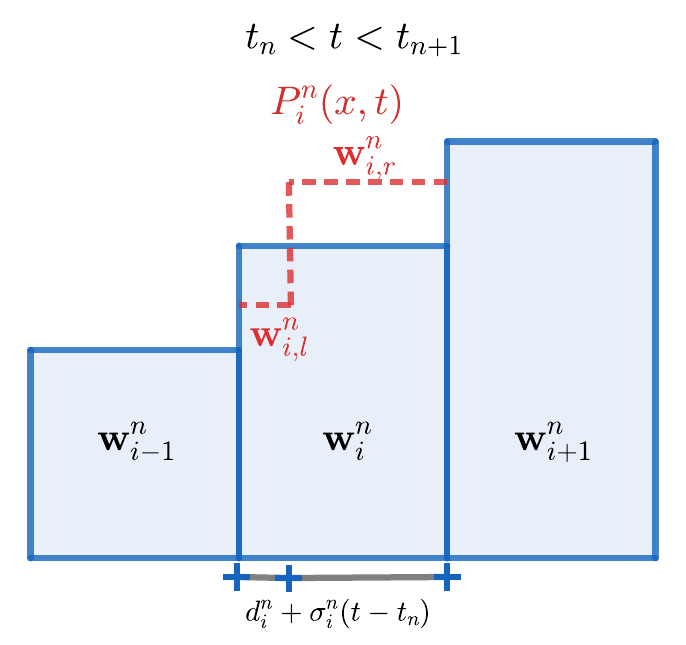}
    }
    \caption{In-cell discontinuous reconstruction operator.}
    \label{fig:In_cell_discontinuous_reconstruction_operator}
\end{figure}

\begin{figure}[htpb]
\centering
\includegraphics[width=0.5\linewidth]{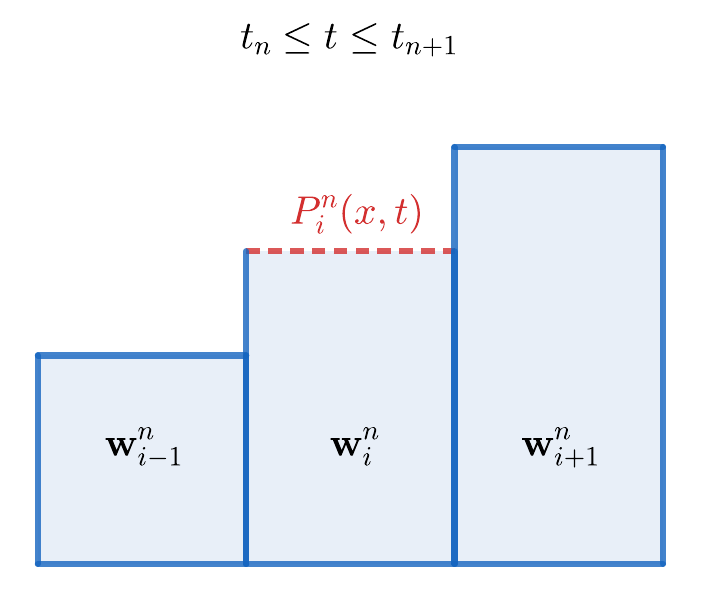}
\caption{Constant reconstruction operator.}
\label{fig:Constant_reconstruction_operator}
\end{figure}
\end{itemize}



\subsection{Fully discrete form}

$\bw_i^{n+1}$ is computed by integrating \eqref{eq:hopc} in the interval $t \in [t_{n}, t_{n+1}]$: 
\begin{equation}
    {\bf w}_{i}^{n+1} = \bw_{i}^{n}-\int_{t_{n}}^{t_{n+1}}\frac{1}{\Delta x}\left(\mathcal{D}_{i+\frac{1}{2}}^{-}(t) + \mathcal{D}_{i-\frac{1}{2}}^{+}(t) + \dashint_{x_{i-\frac{1}{2}}}^{x_{i+\frac{1}{2}}} \mathcal{A}(P_{i}^{n}(x,t))\frac{\partial}{\partial x} P_{i}^{n}(x,t)dx\right)dt,
\end{equation}
where $t_{n+1}=t_n + \Delta t_n$ will be defined in Subsection \ref{sub:time_step}.  Two quadrature formulas
$$
\int_{x_{i-1/2}}^{x_{i+1/2}} f(x) \, dx \approx \Delta x \sum_{m = 1}^{M_{x}} \beta_{m}^{x} f(x_i^{m}), \quad
\int_{t^{n}}^{t^{n+1}} g(t) \, dt \approx \Delta t_n \sum_{l = 1}^{M_{t}} \beta_{l}^{t} g(t^{n,l})
$$
are used to approximate the integrals in space and time, respectively. Taking into account the expression of the reconstruction operators we obtain:

\begin{equation}\label{eq:discrete1}
{\footnotesize{
    \bw_{i}^{n+1} = \begin{cases}
        \displaystyle \bw_{i}^{n} - \frac{\Delta t_n}{\Delta x}\left(\mathcal{D}_{i+\frac{1}{2}}^{-} + \mathcal{D}_{i-\frac{1}{2}}^{+} \right) & \\
        \qquad \displaystyle - \frac{\Delta t_n}{\Delta x}  \int_0^1 \mathcal{A}(\Phi(s; \bw_{i,l}^n, \bw_{i,r}^n)) \frac{\partial \Phi}{\partial s}(s; \bw_{i,l}^n, \bw_{i,r}^n) \, ds  & \text{if $i \in \mathcal{M}_n$,}\\
        \displaystyle \bw_{i}^{n} - \frac{\Delta t_n}{\Delta x}\left(\mathcal{D}_{i+\frac{1}{2}}^{-} + \mathcal{D}_{i-\frac{1}{2}}^{+}\right) & \text{if $i \notin \mathcal{M}_n$ and  $\mathcal{S}_{i}\cap \mathcal{M}_n \neq \emptyset$,}\\
        \displaystyle \bw_{i}^{n} - \frac{\Delta t_n}{\Delta x}\sum_{l=1}^{M_{t}} \beta_{l}^{t}\left(\mathcal{D}_{i+\frac{1}{2}}^{-}(t^{n,l}) + \mathcal{D}_{i-\frac{1}{2}}^{+}(t^{n,l})\right) & \\
        \qquad \displaystyle- \Delta t_n \sum_{l=1}^{M_{t}} \sum_{m=1}^{M_{x}} \beta_{l}^{t} \beta_{m}^{x}\mathcal{A}(P_{i}^{n}(x_{i}^{m},t^{n,l}))
    \partial_x P_{j}^{n}(x_{i}^{m},t^{n,l}) & \text{if $\mathcal{S}_{i} \cap \mathcal{M}_{n} = \emptyset$.} 
    \end{cases}
    }}
\end{equation}

\begin{remark}
    Observe that this method is basically a standard high-order path-conservative method in the smooth region and the standard in-cell discontinuous reconstruction developed in \cite{chalons2019path, pimentel2022cell, pimentel2024high} in the presence of shocks for pure nonconservative systems.
\end{remark}

\subsection{How to mark a cell? Choices of $\sigma_{i}^{n}, \bw_{i,l}^{n}$ and $\bw_{i,r}^{n}$}

Let us assume that a path-conservative Roe linearization for the nonconservative system $\mathcal{A}_\Phi(\bw_l, \bw_r)$ is available (see \cite{toumi1992weak}),  i.e. a matrix-valued function
$\mathcal{A}_\Phi\colon\Omega \times\Omega \mapsto\mathcal{M}_{N\times
N}(\mathbb{R})$ satisfying the following properties:
\begin{itemize}
\item[(R1)]  for each $\bw_l,\bw_r\in\Omega$,
$\mathcal{A}_\Phi(\bw_l,\bw_r)$ has $N$ distinct real eigenvalues
$\lambda_1(\bw_l, \bw_r)$, \dots, $\lambda_N(\bw_l, \bw_r)$;
\item[(R2)] $\mathcal{A}_\Phi(\bw,\bw)=\mathcal{A}(\bw)$, for every
$\bw\in\Omega$;
\item[(R3)]  for any $\bw_l,\bw_r\in\Omega$,
\begin{equation}\label{Roe-gen}
\mathcal{A}_\Phi(\bw_l,\bw_r)\,
(\bw_r-\bw_l)=\int_0^1\mathcal{A}(\Phi(s;\bw_l,\bw_r))\frac{\partial\Phi}{\partial
s}(s;\bw_l,\bw_r)\,ds.
\end{equation}
\end{itemize}
Then, for every $i$ the coordinates $\{ \alpha_k \}$ of  $\bw^{n}_{i+1} - \bw^{n}_{i-1}$ in the basis of eigenvectors of $\mathcal{A}_\Phi(\bw^{n}_{i-1}, \bw^{n}_{i+1})$, are computed, i.e.
$$
\bw^n_{i+1} - \bw^n_{i-1} = \sum_{k=1}^N \alpha_k R^{\Phi}_k(\bw^n_{i-1}, \bw^n_{i+1}).$$
and the intermediate states
$$
\bw^{(k)}_i = \bw^n_{i-1} + \sum_{s=1}^k \alpha_s R_s^{\Phi} (\bw^n_{i-1}, \bw^n_{i+1}),
\quad k = 0, \dots, N,
$$
are defined.
The cell is then marked if there exists an index $k^*$ of those corresponding to genuinely nonlinear fields such that
\begin{equation}\label{eq:choice_k_shock}
    \frac{|\alpha_{k^{*}}|}{\sum_{k=1}^{N} |\alpha_k|} >\overline{\alpha} \quad \text{and}\quad \lambda_{k^*}(\bw_{i-1}^{n}) >\lambda_{k^*}(\bw_{i+1}^{n}),
\end{equation}
where $\overline{\alpha} \in (0,1]$ is a parameter to be selected. Observe that, if \eqref{eq:choice_k_shock} is satisfied with $\overline{\alpha} \approx 1$, then the solution of the linear Riemann problem essentially consists of only one shock wave that, due to the Roe property (R3), is also an admissible discontinuity of the system. In practice we use $\overline{\alpha} = 0.9$. Moreover, the second condition in \eqref{eq:choice_k_shock} allows one to discard non-entropy shocks. In this case, the in-cell discontinuous reconstruction is used \eqref{eq:in_cell_reconstruction} with the following choices:
\begin{equation}\label{eq:reconstructed_states_and_speed}
   \bw^n_{i,l} = \bw^{(k^*-1)}_i, \quad \bw^n_{i,r} = \bw^{(k^*)}_i, \quad   \sigma_i^n = \lambda_{k^*}(\bw^n_{i-1}, \bw^n_{i+1}).
\end{equation}
Summing up, the set $\mathcal{M}_n$ is given by
$$
\mathcal{M}_{n} = \{ i \ \text{ s.t. there exists $k^{*} \in \{1, \dots, N\}$ satisfying \eqref{eq:choice_k_shock} and \eqref{eq:d_property}.} \}
$$

\subsection{Time step}\label{sub:time_step}
Following \cite{pimentel2022cell, pimentel2024high} the time step $\Delta t_n$ is chosen as follows:
\begin{equation}\label{restricted_time_step}
    \Delta t_n = \min(\Delta t^c_n, \Delta t^r_n).
\end{equation}
Here
\begin{equation}
    \Delta t_n^c = CFL\cdot \min_i\left( \frac{ \Delta x}{ \max_{l}|\lambda_{i,l}|} \right)
\end{equation}
where $CFL \in (0,1)$ is the stability parameter, $\lambda_{i,l}, \dots, \lambda_{i,N}$ represent  the eigenvalues of
$\mathcal{A}(\bw_i^n)$, and
\begin{equation}
    \Delta t_n^r = \min_{i \in \mathcal{M}_n} \begin{cases}
    \displaystyle \frac{1 - d_i^n}{|\sigma_i^n|}\Delta x, & \quad if \quad \sigma_i^n>0, \\ \\
    \displaystyle \frac{d_i^n}{|\sigma_i^n|}\Delta x, & \quad if \quad \sigma_i^n<0,
    \end{cases} 
\end{equation}
where $d_{i}^{n}$ and $\sigma_{i}^{n}$
are given in \eqref{conservation} and \eqref{eq:reconstructed_states_and_speed}, respectively, only for marked cells. Observe that, besides the stability requirement, this choice of time step ensures that no discontinuous reconstruction leaves a marked cell. Note that this time restriction could be avoided as in \cite{chalons2019path} but, for simplicity, we keep it.

\subsection{Path-conservative fluctuations}
The fluctuations corresponding to Roe method will be considered here
\begin{equation}\label{NCFR} 
\mathcal{D}^\pm (\bw_l, \bw_r) = \frac{1}{2}
\mathcal{A}_\Phi(\bw_l, \bw_r) \, (\bw_r - \bw_l) \pm \frac{1}{2}
|\mathcal{A}_\Phi(\bw_l, \bw_r)| \, (\bw_r - \bw_l),
\end{equation}
where  $\mathcal{A}_\Phi(\bw_l, \bw_r)$ is the  Roe linearization used to mark cells. 
As usual $|\mathcal{A}^\pm_\Phi(\bw_l, \bw_r)|$ represent  the matrices whose
 eigenvalues are the absolute value of $\lambda_1(\bw_l, \bw_r)$, \dots, $\lambda_N(\bw_l, \bw_r)$
with same eigenvectors. 

\subsection{First-order method}

In this case of first-order methods, the space-time reconstruction operator is the piecewise constant one, i.e.
$$Q_{i}^{n}(x,t) = \bw_{i}^{n},$$
so that the first-order in-cell discontinuous reconstruction path-conservative method \eqref{eq:discrete1} method is reduced to:
\begin{equation}\label{eq:discretefo}
    \bw_{i}^{n+1} = \begin{cases}
        \displaystyle \bw_{i}^{n} - \frac{\Delta t_n}{\Delta x}\left(\mathcal{D}_{i+\frac{1}{2}}^{-} + \mathcal{D}_{i-\frac{1}{2}}^{+}\right) & \\
        \qquad \displaystyle - \frac{\Delta t_n}{\Delta x}  \int_0^1 \mathcal{A}(\Phi(s; \bw_{i,l}^n, \bw_{i,r}^n)) \frac{\partial \Phi}{\partial s}(s; \bw_{i,l}^n, \bw_{i,r}^n) \, ds & \text{if $i \in \mathcal{M}_n$,}\\
        \displaystyle \bw_{i}^{n} - \frac{\Delta t_n}{\Delta x}\left(\mathcal{D}_{i+\frac{1}{2}}^{-} + \mathcal{D}_{i-\frac{1}{2}}^{+}\right) & \text{if $i \notin \mathcal{M}_n$.}\\
    \end{cases}
\end{equation}

\subsection{Second-order method}

The MUSCL-Hancock space-time reconstruction operator (see\cite{van1974towards, van1984relation}) is considered in this case:
\begin{equation}
Q_{i}^{n}(x,t)= \bw_{i}^{n} + \widetilde{\partial_x\bw}_{i}^{n}(x-x_{i})- \mathcal{A}(\bw_{i}^{n})\widetilde{\partial_x\bw}_{i}^{n}(t-t^{n}),
    \end{equation}
and the stencil is $\mathcal{S}_i = \{i-1,i,i+1\}$.
Here, $\widetilde{\partial_x\bw}_{j}^{n}$ is the \textit{minmod} approximation of the first order spatial derivative of 
$\bw$ at $x_i$ at time $t_n$, whose $k$th component is given by
$$
\left(\widetilde{\partial_x\bw}_{i}^{n}\right)_k = {\rm{minmod}}\left( \frac{w_{i+1,k}^{n}-w_{i,k}^{n}}{\Delta x}, \frac{w_{i,k}^{n}-w_{i-1,k}^{n}}{\Delta x}\right),$$
where $w_{i,k}^n$ represents the $k$th component of $\bw_i^n$ and 
$$\rm{minmod}(a,b,c) = \begin{cases}
     \min\{a,b,c\} & \text{if $a,b,c>0$,}  \\
     \max\{a,b,c\} & \text{if $a,b,c<0$,}  \\
     0 & \text{otherwise.}
     \end{cases}$$
The mid-point rule is chosen as quadrature formula both in time and space, so that the second-order in-cell discontinuous reconstruction path-conservative method \eqref{eq:discrete1} is given by:
\begin{equation}\label{eq:discrete2}
{\small{
    \bw_{i}^{n+1} = \begin{cases}
        \displaystyle \bw_{i}^{n} - \frac{\Delta t_n}{\Delta x}\left(\mathcal{D}_{i+\frac{1}{2}}^{-} + \mathcal{D}_{i-\frac{1}{2}}^{+}\right) & \\
        \qquad \displaystyle - \frac{\Delta t_n}{\Delta x} \int_0^1 \mathcal{A}(\Phi(s; \bw_{i,l}^n, \bw_{i,r}^n)) \frac{\partial \Phi}{\partial s}(s; \bw_{i,l}^n, \bw_{i,r}^n) \, ds & \text{if $i \in \mathcal{M}_n$,}\\
        \displaystyle \bw_{i}^{n} - \frac{\Delta t_n}{\Delta x}\left(\mathcal{D}_{i+\frac{1}{2}}^{-} + \mathcal{D}_{i-\frac{1}{2}}^{+}\right) & \text{if $i \notin \mathcal{M}_n$, $\mathcal{S}_{i} \cap \mathcal{M}_{n} \neq \emptyset$,}\\
        \displaystyle \bw_{i}^{n} - \frac{\Delta t_n}{\Delta x}\left(\mathcal{D}_{i+\frac{1}{2}}^{-}(t_{n+\frac{1}{2}}) + \mathcal{D}_{i-\frac{1}{2}}^{+}(t_{n+\frac{1}{2}}) \right) & \\
        \qquad \displaystyle - \Delta t_n \mathcal{A}(\bw_{i}^{n+\frac{1}{2}})\widetilde{\partial_x\bw}_{i}^{n} & \text{if $\mathcal{S}_{i} \cap \mathcal{M}_{n} = \emptyset$,} 
    \end{cases}
    }}
\end{equation}
where $\displaystyle t_{n+\frac{1}{2}} = t_n + \frac{\Delta t_n}{2}$ so
    $$\mathcal{D}_{i+\frac{1}{2}}^{\pm}(t_{n+\frac{1}{2}}) = \mathcal{D}^{\pm}\left(P_{i}^{n}\left(x_{i+\frac{1}{2}},t_n + \frac{\Delta t_n}{2}\right),P_{i+1}^{n}\left(x_{i+\frac{1}{2}},t_n + \frac{\Delta t_n}{2}\right)\right) $$
    and
    $$\bw_{i}^{n+\frac{1}{2}} =  P_{i}^{n}\left(x_{i},t_n + \frac{\Delta t_n}{2}\right).$$

\subsection{Shock-capturing and Rankine-Hugoniot compatibility}

We recall the theorem in \cite{chalons2019path, pimentel2022cell, pimentel2024high} that states that the scheme \eqref{eq:discrete1} captures exactly isolated shocks in the following way:

\begin{theorem}\label{th}
Assume that $\bu_l$ and $\bu_r$ can be joined by an entropy shock of speed $\sigma$ satisfying \eqref{con_R-H}. Then, the numerical method  \eqref{eq:discrete1} provides an exact numerical solution
of the Riemann problem with initial conditions
$$
\bw(x,0) = \begin{cases}
\bw_l= T^{-1}(\bu_l) & \text{if $x < 0$,} \\
\bw_r=T^{-1}(\bu_r) & \text{otherwise,}
\end{cases} 
$$
in the sense that
\begin{equation}\label{presave}
\bw_i^n = \frac{1}{\Delta x}\int_{x_{i-1/2}}^{x_{i+1/2}} \bw(x,t_n) \, dx, \quad \forall i, n
\end{equation}
where $\bw(x,t)$ is the exact solution.

\end{theorem}

\begin{proof}
    Since \eqref{con_R-H} and \eqref{gen_R-H} are equivalents for the given family of paths, the proof reduces to the one in \cite{pimentel2022cell}.
\end{proof}

\subsection{Extensions}

In this work only first- and second-order methods are considered, but this strategy can be applied to higher order methods as it has been done in \cite{pimentel2024high}. Moreover in the mentioned paper the methods are combined with the MOOD approach. The method can also be extended to problem with source terms in a well-balanced manner as it was done in \cite{pimentel2022cell}.  

\section{Euler equations with an interface and shear}\label{sec:Euler}

As an example of a non-linear system 
we consider the Euler equations of gas dynamics for two immiscible materials, plus a shear wave
modeling the presence of transverse flow.
The resulting  system of PDEs  is a good representative example for many other applications 
in compressible flows as well as shallow water flows.
It is well known  (see e.g. 
\cite{Quirk1996129,ABGRALL1996150} and may others) that the use of conservative variables produce spurious oscillations close to contact discontinuities related to material or velocity interfaces.

\subsection{Conservative  form}

We consider the Euler equations  enhanced with a marker function to detect moving interfaces and two velocities:

\begin{equation}\label{euler0}
\partial_t\left[
\begin{array}{c}
\rho\\ \rho u\\ \rho v \\ \rho E\\ \rho c
\end{array}
\right]+\partial_x\left[
\begin{array}{c}
\rho u\\ \rho u^2 + p\\\rho u v \\ \rho u H\\ \rho u c
\end{array}
\right]=0
\end{equation}
with  
\begin{equation}\label{euler1}
 E = e + k  \;,\quad  H= h + k   \;,\quad h = e +p/\rho, \quad k = \frac{1}{2}(u^2+v^2).
\end{equation}

The fluid is assumed to be composed of two immiscible ideal gases with equations of states
$$
p = (\gamma -1) \rho e, \quad \gamma = c\gamma_p  + (1-c)\gamma_s,
$$
Typical values are  for example $\gamma_p =1.4$ for air, 
and $\gamma_s=2.8$ for water.
The marker function $c$ takes  values in $\{0,1\}$, $c=0$ denoting one material and $c=1$ the other one.
The equation of state of the fluid is then given by 
$$
p = c (\gamma_p-1)\rho e  + (1-c)( \gamma_s-1)\rho e\;.
$$

Two cases will be considered here
\begin{itemize}

\item two materials, one velocity ($v = 0$),

\item one material ($c = 1$), two velocities.

\end{itemize}




Rankine-Hugoniot conditions for the variables $\rho$, $\rho u$, and $E$ writes as follows:
\begin{equation}\label{RK_conEuler}
    \begin{cases}
        \sigma (\rho_r-\rho_l) = \rho_r u_r - \rho_l u_l,\\
        \sigma (\rho_r u_r - \rho_l u_l) = \rho_r u_r^2 + p_r - \rho_l u_l^2 - p_l,\\
        \sigma (\rho_r E_r - \rho_l E_l) = \rho_r u_r H_r - \rho_l u_l H_l.
    \end{cases}
\end{equation}
From the first one, we derive:
$$ m:= \rho_r(u_r - \sigma) = \rho_l (u_l - \sigma).$$
Therefore, from the Rankine-Hugoniot condition for $c$:
$$
\sigma (\rho_r c_r - \rho_l c_l) = \rho_r c_r u_r - \rho_l c_l u_l$$
we derive
$$
m(c^+ - c^-) = 0.$$
Since for shock waves $m \not= 0$ we obtain that $c$ is constant through shocks, so that the RH condition for $\rho c$ is the same that the one for $\rho$. The same happens for the variable $v$. Therefore, \eqref{RK_conEuler} are the RH conditions to be considered for shocks. The variables $\rho c$ and $\rho v$ can only be discontinuous at contact discontinuities that travel at speed $u := u^+  = u^-$.

\subsection{Nonconservative form with $(\rho, u , v, p, c)$}
We can consider the nonconservative form of the system reading
\begin{equation}\label{euler_nonconservative}
\partial_t \bw + \mathcal{A}(\bw) \partial_x \bw=0\;,\quad \bw = \left[\begin{array}{c}
\rho\\  u\\ v \\ p \\   c
\end{array}
\right] \;,\quad
\mathcal{A}(\bw) =\left(\begin{array}{ccccc}
u & \rho & 0 & 0 & 0\\
0 & u & 0 & 1/\rho& 0\\
0 & 0 & u & 0 & 0 \\
0 & \gamma p & 0 & u& 0\\
0 & 0&0& 0 & u
\end{array}
\right),
\end{equation}
The matrix $\mathcal{A}$ has  5 eigenvalues given by 
$$ \lambda_1 = u - a, \ \ \lambda_{2,3,4} = u, \ \ \lambda_{5} = u +a$$
with 
$$
a^2 = \dfrac{\gamma p}{\rho}.
$$
The corresponding eigenvectors of $A$ are given by
$$R_{1} = \left(\begin{array}{c}
    \frac{\rho}{\gamma p}  \\
     -\frac{1}{\sqrt{\gamma\rho p}} \\
     0
     \\
     1 \\
     0
\end{array}\right), \ R_{2} = \left(\begin{array}{c}
     1  \\
     0 \\
     0\\
     0 \\
     0
\end{array}\right), \ R_{3} = \left(\begin{array}{c}
     0  \\
     0 \\
     1\\
     0 \\
     0
\end{array}\right), \ R_{4} = \left(\begin{array}{c}
     0  \\
     0 \\
     0\\
     0 \\
     1
\end{array}\right), \ R_{5} = \left(\begin{array}{c}
      \frac{\rho}{\gamma p}  \\ 
      \frac{1}{\sqrt{\gamma\rho p}} \\
     1 \\
     0
\end{array}\right).$$

\subsection{Correct choice of the family of paths}

Given a family of paths 
$$\Phi(s;\bw_l, \bw_r) =  \left(\begin{array}{c}
     \Phi_\rho(s;\bw_l, \bw_r)  \\
     \Phi_u(s;\bw_l, \bw_r) \\
     \Phi_v(s;\bw_l, \bw_r)\\
     \Phi_p(s;\bw_l, \bw_r) \\
     \Phi_c(s;\bw_l, \bw_r)
\end{array}\right),$$
the generalized Rankine-Hugoniot conditions \eqref{gen_R-H} for the variables $\rho$, $u$, $p$ are given by:
\begin{equation}\label{RK_nonconEuler}
    \begin{cases}
        \sigma (\rho_r - \rho_l) = \rho_r u_r - \rho_l u_l, \\
        \displaystyle \sigma (u_r - u_l) = \frac{u_r^2}{2} - \frac{u_l^2}{2} + \int_0^1 \frac{\Phi_p'(s)}{\Phi_\rho(s)} ds,\\
        \sigma (p_r - p_l) = \int_0^1 (\gamma \Phi_p(s) \Phi'_u(s) + \Phi_u(s) \Phi'_p(s)) ds.\\
    \end{cases}
\end{equation}
Again, $c$ (and thus $\gamma$) and $v$ are constant though shocks, so that their corresponding RH conditions are equal to the one of $\rho$.

\begin{theorem}\label{th:path_Euler}
    Given $\bw_l$ and $\bw_r$, if the family of paths
    \begin{equation}\label{eq:correct_path}
        \Phi(s;\bw_l, \bw_r) =  \left(\begin{array}{c}
     \displaystyle\frac{1}{\frac{1}{\rho_l} + s \left(\frac{1}{\rho_r} - \frac{1}{\rho_l}\right)}  \\
     u_l + s(u_r-u_l) \\
     v_l + s(v_r-v_l) \\
     p_l + s(p_r-p_l) \\
     c_l + s(c_r-c_l)
\end{array}\right)
    \end{equation}
is chosen, then the Rankine-Hugoniot conditions \eqref{RK_conEuler} and \eqref{RK_nonconEuler} are equivalent.
\end{theorem}

\begin{proof}
    The proof is given in Appendix \ref{appendix_proof}.
\end{proof}

\subsubsection{Nonconservative Roe matrix}
Given the family of paths \eqref{eq:correct_path}, a Roe matrix $\mathcal{A}_{\Phi}(\bw_{l}, \bw_{r})$ must verify:
\begin{equation}\label{Roe_matrix_condition}
    \mathcal{A}_{\Phi}(\bw_{l}, \bw_{r})(\bw_{r}-\bw_{l}) = \int_{0}^{1} \mathcal{A}(\Phi(s; \bw_{l}, \bw_{r}))\partial_{s} \Phi(s; \bw_{l}, \bw_{r}). 
\end{equation} 
It can be easily checked that a Roe matrix verifying \eqref{Roe_matrix_condition} is given by
    \begin{equation}\label{eq:Roe_matrix_Euler}
        \mathcal{A}_{\Phi}(\bw_l, \bw_r) = \left(\begin{array}{ccccc}
\bar{u} & \bar{\rho} & 0 & 0 & 0\\
0 & \bar{u} & 0 &  \overline{(\frac{1}{\rho})}& 0\\
0 & 0 & \bar{u} & 0 & 0\\
0 & \bar{\gamma} \bar{p} & 0 & \bar{u}& 0\\
0 & 0&0& 0&  \bar{u}
\end{array}
\right)  
    \end{equation}
where
\[
\begin{aligned}
&\bar{\rho} = \frac{\rho_r + \rho_l}{2},\quad
\bar{u} = \frac{u_r + u_l}{2},\quad
\bar{p} = \frac{p_r + p_l}{2},\\[4pt]
&
\bar{\gamma} = \bar{c}\gamma_p + (1-\bar{c})\gamma_s,
\quad \bar{c} = \frac{c_r + c_l}{2},\quad
\overline{\left(\frac{1}{\rho}\right)} = \frac{1}{2}\left(\frac{1}{\rho_l} + \frac{1}{\rho_r}\right).
\end{aligned}
\]
The eigenvalues of  $\mathcal{A}_{\Phi}$ \eqref{eq:Roe_matrix_Euler} are given by:
        $$ \lambda_1 = \bar{u} - \bar{a}, \ \ \lambda_{2,3} = \bar{u}, \ \ \lambda_{4} = \bar{u} +\bar{a},$$
        where $\bar{a} = \sqrt{\bar{\gamma} \,\bar{p} \overline{\left(\frac{1}{\rho}\right)}}$.
    with he eigenvectors
        $$R_{1}^{\Phi} = \left(\begin{array}{c}
     \frac{\bar{\rho}}{\bar{\gamma}\bar{p}}  \\
     -\frac{\bar{a}}{\bar{\gamma}\bar{p}} \\
     0\\
     1 \\
     0
\end{array}\right), \ \ R_{2}^{\Phi} = \left(\begin{array}{c}
     1  \\
     0 \\
     0 \\
     0 \\
     0
\end{array}\right), \ \ R_{3}^{\Phi} = \left(\begin{array}{c}
     0  \\
     0 \\
     1 \\
     0 \\
     0
\end{array}\right),\ \ R_{4}^{\Phi} = \left(\begin{array}{c}
     0  \\
     0 \\
     0 \\
     0 \\
     1
\end{array}\right) \ \ R_{4}^{\Phi} = \left(\begin{array}{c}
     \frac{\bar{\rho}}{\bar{\gamma}\bar{p}}  \\
     \frac{\bar{a}}{\bar{\gamma}\bar{p}} \\
     0\\
     1 \\
     0
\end{array}\right).$$

Note that, close to an isolated shock wave, the Roe matrix is only used to mark the cells. The convergence to
the right weak solution is guaranteed by the use of the DR method with the proper shock speed. We have now all the ingredients to apply the in-cell discontinuous reconstruction method \eqref{eq:discrete1}.
 
 \section{Numerical tests}\label{sec:4}
 Several tests will be considered for validation and application for each model. The following notation will be used in the simulations:

 \begin{itemize}
     \item O1\_Primitive+DR: first-order path-conservative in-cell discontinuous reconstruction method \eqref{eq:discretefo}.
     \item O2\_Primitive+DR: second-order path-conservative in-cell discontinuous reconstruction method \eqref{eq:discrete2}.
     \item O1\_Conservative: standard first-order fully conservative Roe method.
    \item O1\_Primitive: standard first-order path-conservative Roe method without discontinuous in-cell reconstruction.
    \item O1\_Hybrid: first-order primitive-conservative method without discontinuous in-cell reconstruction. This approach is similar  to the one in \cite{GABURRO2024117311, audebert2006hybrid}: once the cells have been marked, O1\_Conservative is used in marked cells and O1\_Primitive in the rest. The marking is relaxed in comparison with the methods based on discontinuous reconstruction: $\bar \alpha = 0.75$ is chosen in \eqref{eq:choice_k_shock}.
     \item Reference\_Solution: it is computed with the O1\_Hybrid scheme using a mesh of 5000 points.
 \end{itemize}
 In all tests the CFL number is fixed to 0.5.

 \subsection{Euler equations with an interface}

 All tests are performed in the spatial domain $I = [0,1]$ using a mesh of $N=200$ points. In all test cases $\gamma_p = 1.4$ and $\gamma_s = 1.67$. We remember $c=1$ refers to $\gamma_p$ and $c=0$ to $\gamma_s$.  In all these tests, $v \equiv 0$, and we evaluate the ability of the schemes to handle a material interface.

 \paragraph{Stationary contact discontinuity}

The initial condition is given by
 $$
 \bw_{0}(x) = (\rho_0(x), u_0(x), v_0(x), p_0(x), c_0(x))^{T} = \begin{cases}
     (1, 0, 0, 1/1.4, 1)^{T} & \text{if $x<0.5$,}\\
     (0.1819, 0, 0, 1/1.4, 0)^{T} & \text{if $x\geq0.5$.}
 \end{cases}
$$

\begin{figure}[h]
    \centering
    \subfloat[Variable $\rho$]{
        \includegraphics[width=0.45\textwidth]{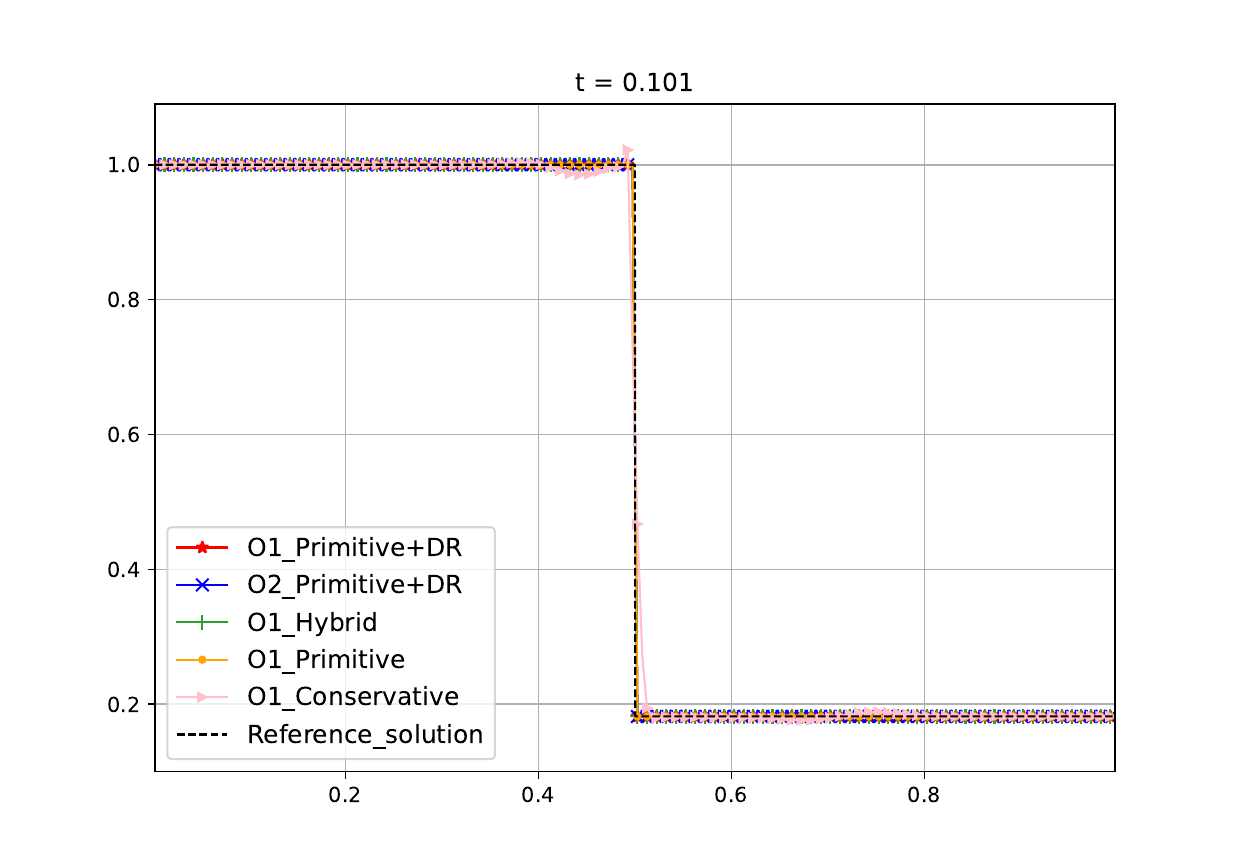}
    }
    \subfloat[Variable $u$]{
        \includegraphics[width=0.45\textwidth]{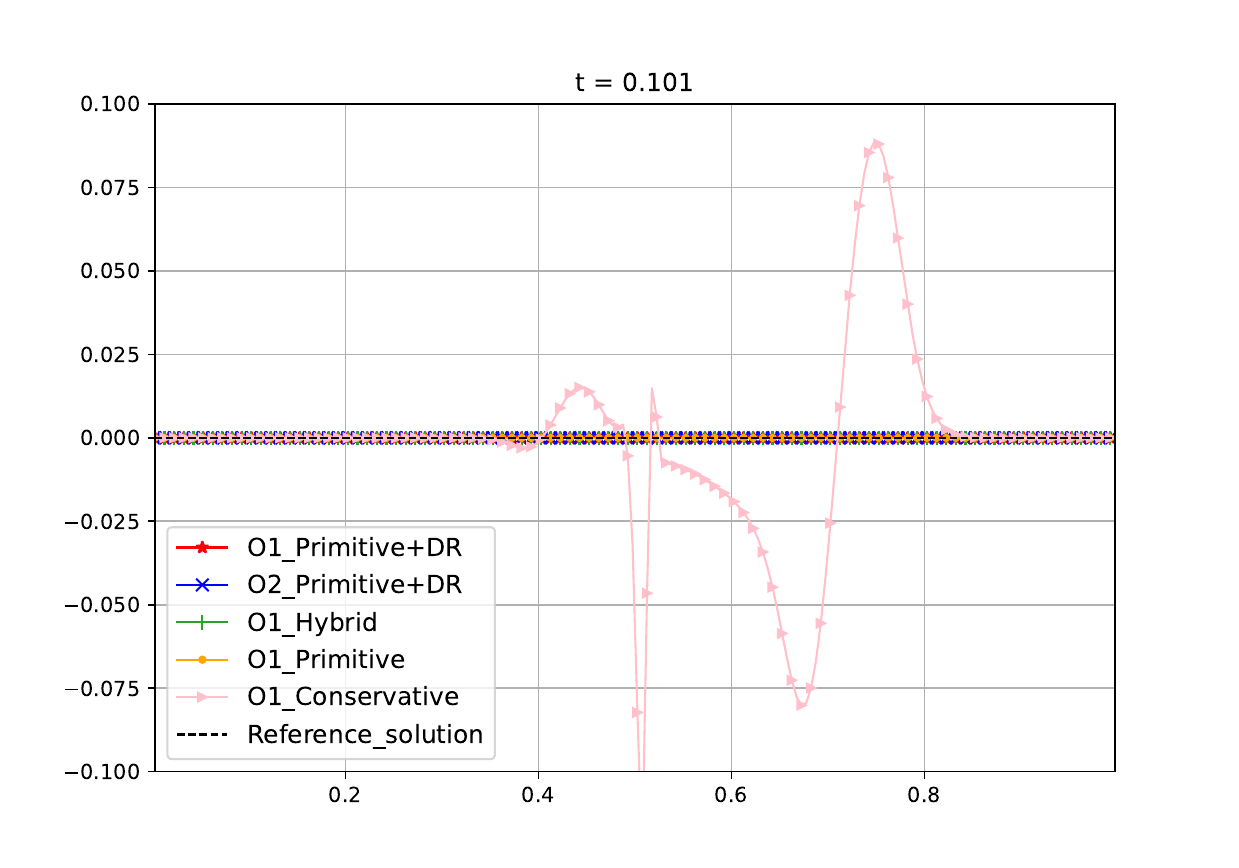}
    }
    \\
    \subfloat[Variable $p$]{
        \includegraphics[width=0.45\textwidth]{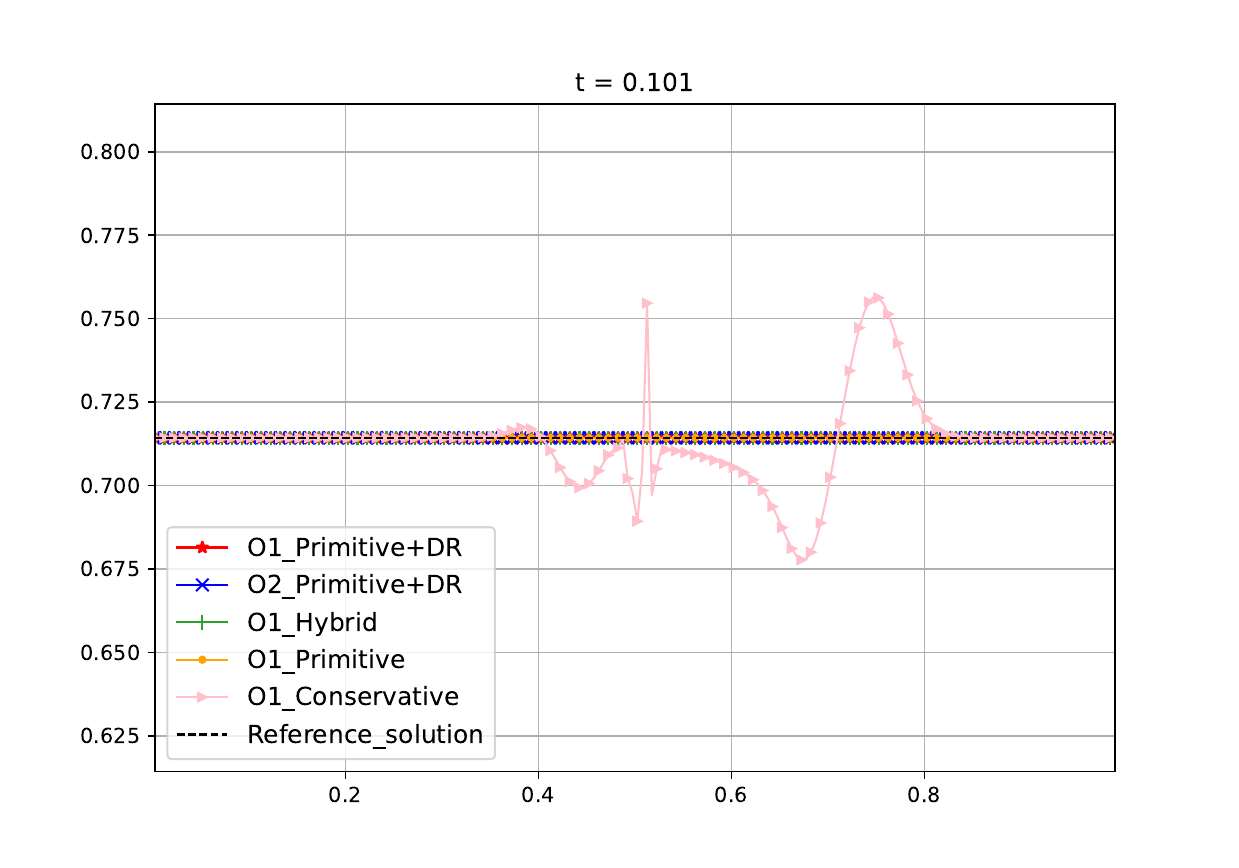}
    }
    \subfloat[Variable $c$]{
        \includegraphics[width=0.45\textwidth]{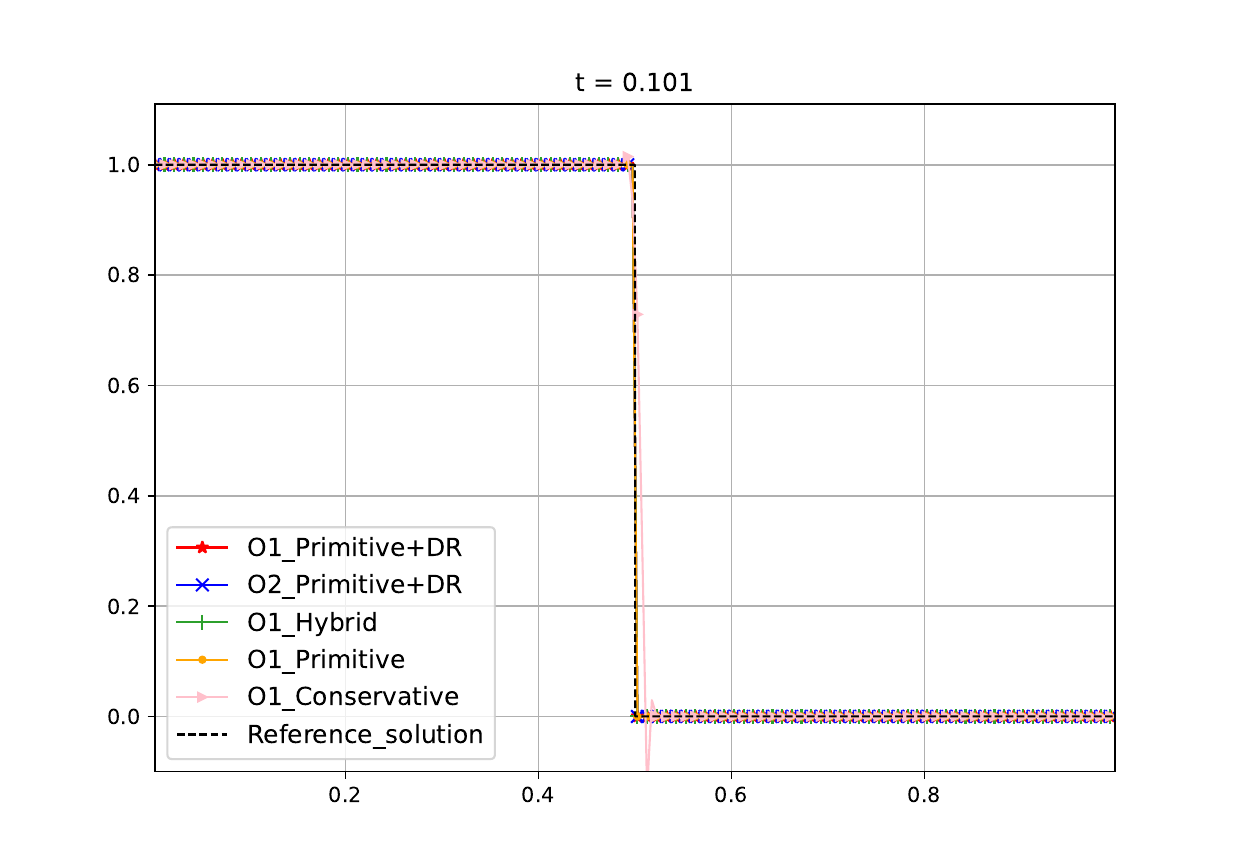}
    }

    \caption{Numerical results for Test: Stationary contact discontinuity.}
    \label{Test_Isolated_contact}
\end{figure}

Figure \ref{Test_Isolated_contact} shows the numerical results: we observe how all numerical schemes capture well the isolated contact discontinuity except the fully conservative one that introduce spurious oscillations. This is what was pointed out in \cite{Quirk1996129}. Since no shock wave is detected, in this case the numerical results provided by O1\_Primitive and O1\_Hybrid coincide. In Table \ref{tab:Error_TestStationaryContact}, we report the numerical errors obtained with the O1\_Primitive+DR, O2\_Primitive+DR and O1\_Conservative schemes for different meshes. The first two schemes achieve errors at machine precision, whereas the conservative formulation exhibits nonzero errors that decrease with mesh refinement, as expected.

\begin{table}[!h]
\centering

\begin{tabular}{|c|c|c|c|c|}
\hline
\multicolumn{5}{|c|}{O1\_Primitive+DR} \\
\hline
N &
$||\Delta \rho||_1$ &
$||\Delta u||_1$ &
$||\Delta p||_1$ &
$||\Delta c||_1$ \\
\hline
100 & 4.01e-16 & 0.00e+00 & 9.54e-16 & 3.47e-16 \\
200 & 3.04e-16 & 0.00e+00 & 2.60e-16 & 2.60e-16 \\
400 & 2.17e-17 & 0.00e+00 & 2.60e-16 & 0.00e+00 \\
800 & 9.76e-17 & 0.00e+00 & 1.73e-16 & 8.67e-17 \\
\hline
\end{tabular}

\vspace{0.4cm}

\begin{tabular}{|c|c|c|c|c|}
\hline
\multicolumn{5}{|c|}{O2\_Primitive+DR} \\
\hline
N &
$||\Delta \rho||_1$ &
$||\Delta u||_1$ &
$||\Delta p||_1$ &
$||\Delta c||_1$ \\
\hline
100 & 4.01e-16 & 0.00e+00 & 9.54e-16 & 3.47e-16 \\ 
   200 & 3.04e-16 & 0.00e+00 & 2.60e-16 & 2.60e-16 \\ 
   400 & 2.17e-17 & 0.00e+00 & 2.60e-16 & 0.00e+00 \\ 
   800 & 9.76e-17 & 0.00e+00 & 1.73e-16 & 8.67e-17 \\ 
\hline
\end{tabular}

\vspace{0.4cm}


\begin{tabular}{|c|c|c|c|c|c|c|c|c|}
\hline
\multicolumn{9}{|c|}{O1\_Conservative} \\
\hline
N &
$||\Delta \rho||_1$ & Order &
$||\Delta u||_1$ & Order &
$||\Delta p||_1$ & Order &
$||\Delta c||_1$ & Order \\
\hline
100 & 7.27e-03 & 0.00 & 2.68e-02 & 0.00 & 1.40e-02 & 0.00 & 1.45e-02 & 0.00 \\
200 & 5.06e-03 & 0.52 & 1.36e-02 & 0.98 & 6.67e-03 & 1.07 & 6.42e-03 & 1.18 \\
400 & 3.14e-03 & 0.69 & 6.52e-03 & 1.06 & 3.11e-03 & 1.10 & 3.93e-03 & 0.71 \\
800 & 1.98e-03 & 0.67 & 3.18e-03 & 1.04 & 1.51e-03 & 1.05 & 2.71e-03 & 0.54 \\
\hline
\end{tabular}

\caption{Test \textit{Stationary contact discontinuity}: order of accuracy for O1\_Primitive+DR, O2\_Primitive+DR and O1\_Conservative: $L^{1}$ errors $||\Delta \cdot||_1$ at time $t=0.1$.}
\label{tab:Error_TestStationaryContact}
\end{table}

\paragraph{Isolated shock wave}
  The initial condition is given by
 $$
 \bw_{0}(x) = (\rho_0(x), u_0(x), v_0(x), p_0(x), c_0(x))^{T} = \begin{cases}
     (\rho_l, u_l, 0, p_l, 1)^{T} & \text{if $x<0.5$,}\\
     (\rho_r, u_r, 0, p_r, 1)^{T} & \text{if $x\geq0.5$,}
 \end{cases}
 $$
 where 
 $$\rho_l = 1.225, \ \ u_l = 0.0, \ \ p_l = 101325, \ \ \gamma_l = \gamma_p, \ \ M_s = 0.8,$$
 $$\rho_r = \frac{\rho_l(\gamma_l+1)M_s^2}{2+(\gamma_l-1)M_s^2}, \ u_r = \frac{2}{\gamma_l+1}(M_s-1/M_s)\sqrt{\gamma_lp_l/\rho_l}, \ p_r = p_l\left(\frac{2\gamma_l M_s^2}{\gamma_l+1}-\frac{\gamma_l-1}{\gamma_l+1}\right)$$

\begin{figure}[h]
    \centering

    \subfloat[Variable $\rho$]{
        \includegraphics[width=0.30\textwidth]{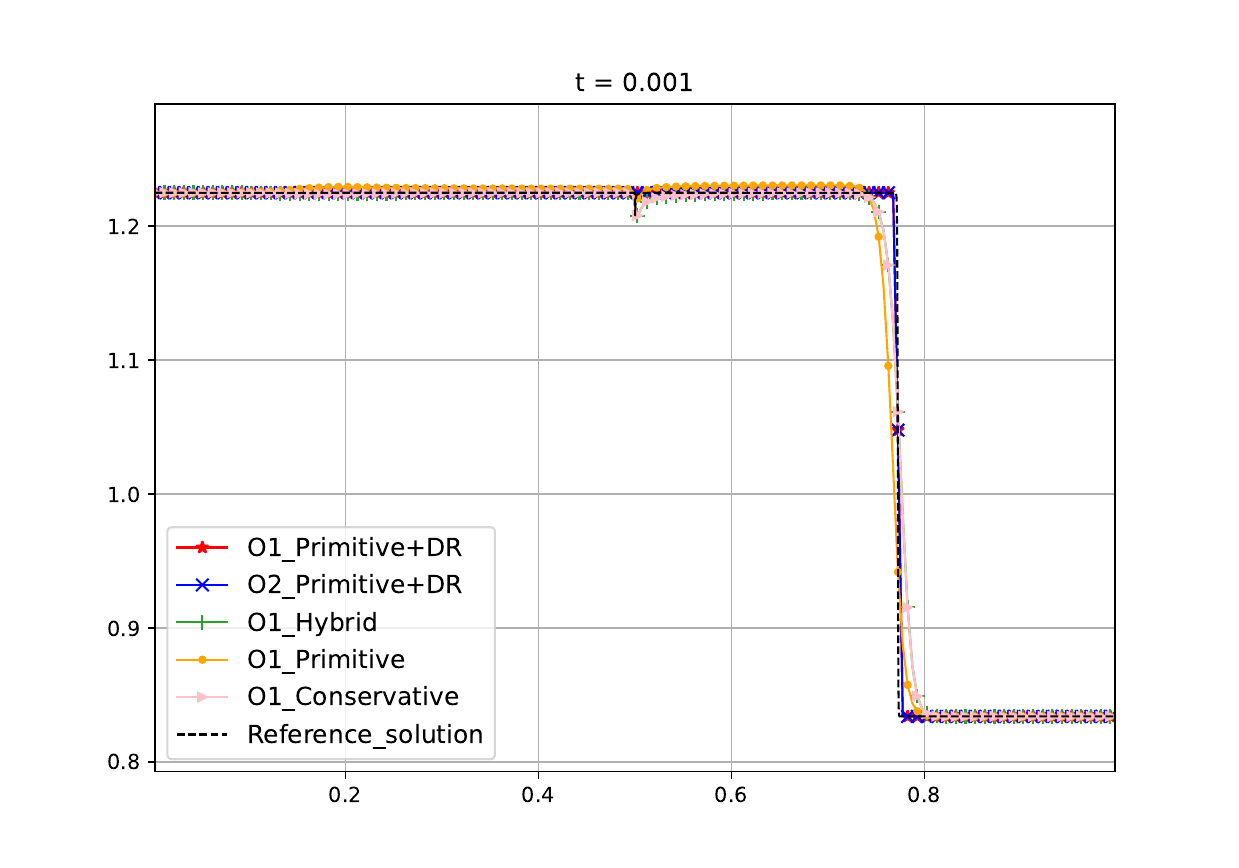}
    }
    \subfloat[Variable $\rho$. Zoom shock]{
        \includegraphics[width=0.30\textwidth]{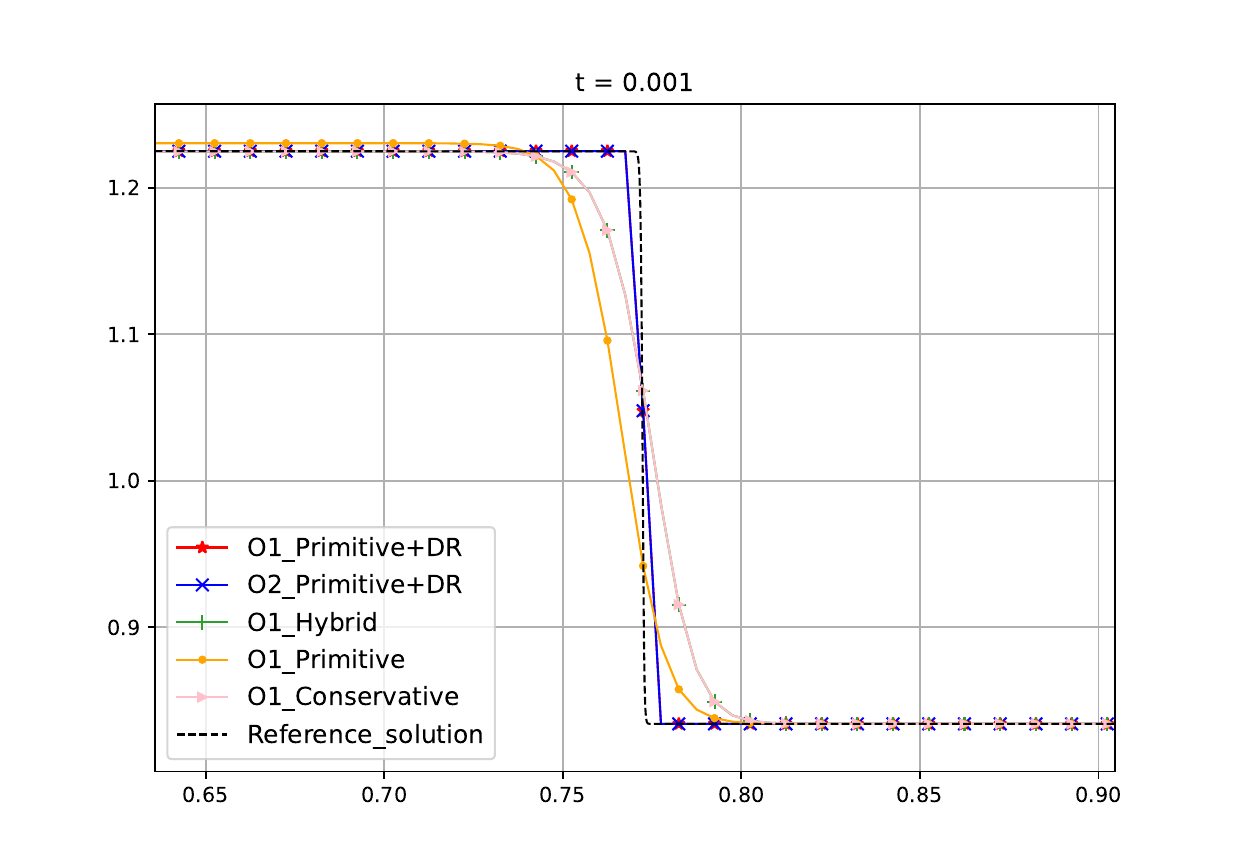}
    }
    \subfloat[Variable $\rho$. Zoom left state]{
        \includegraphics[width=0.30\textwidth]{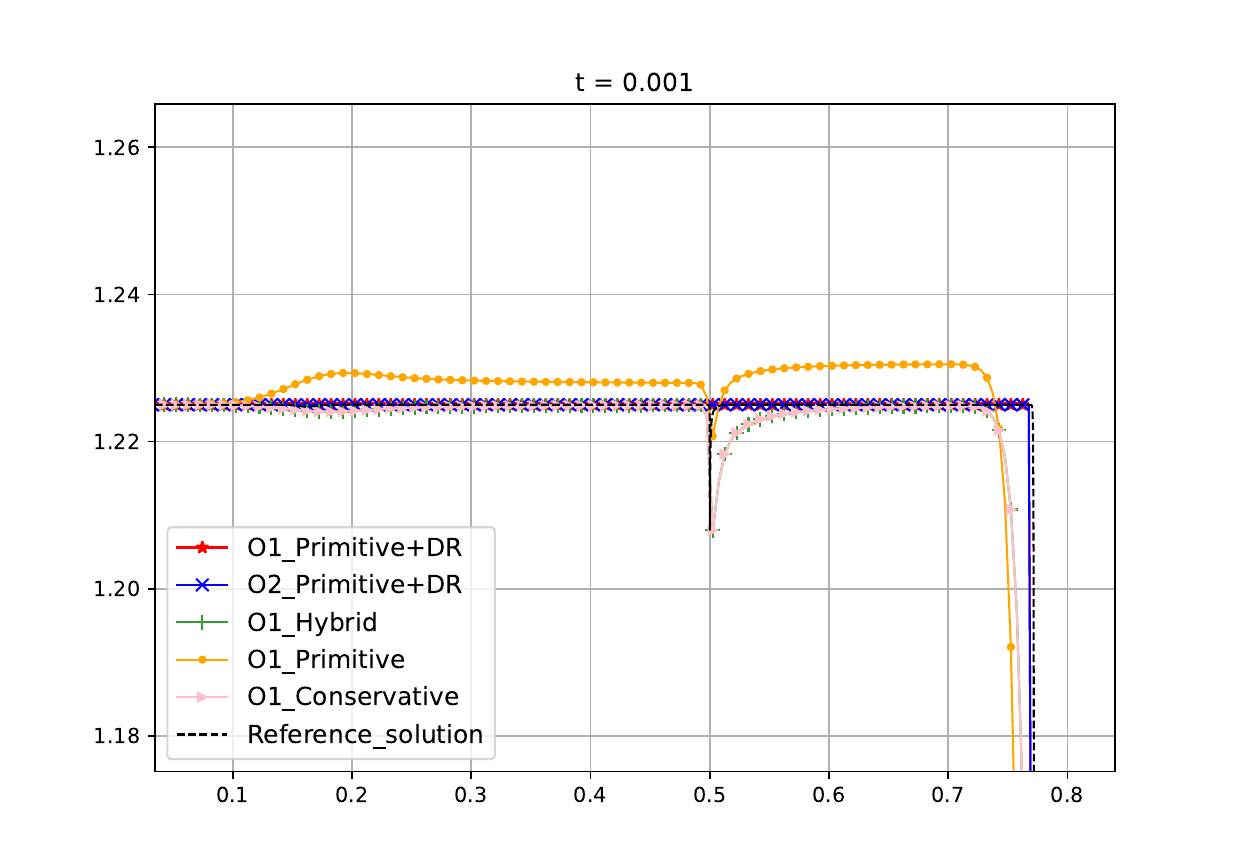}
    }

    \subfloat[Variable $u$]{
        \includegraphics[width=0.30\textwidth]{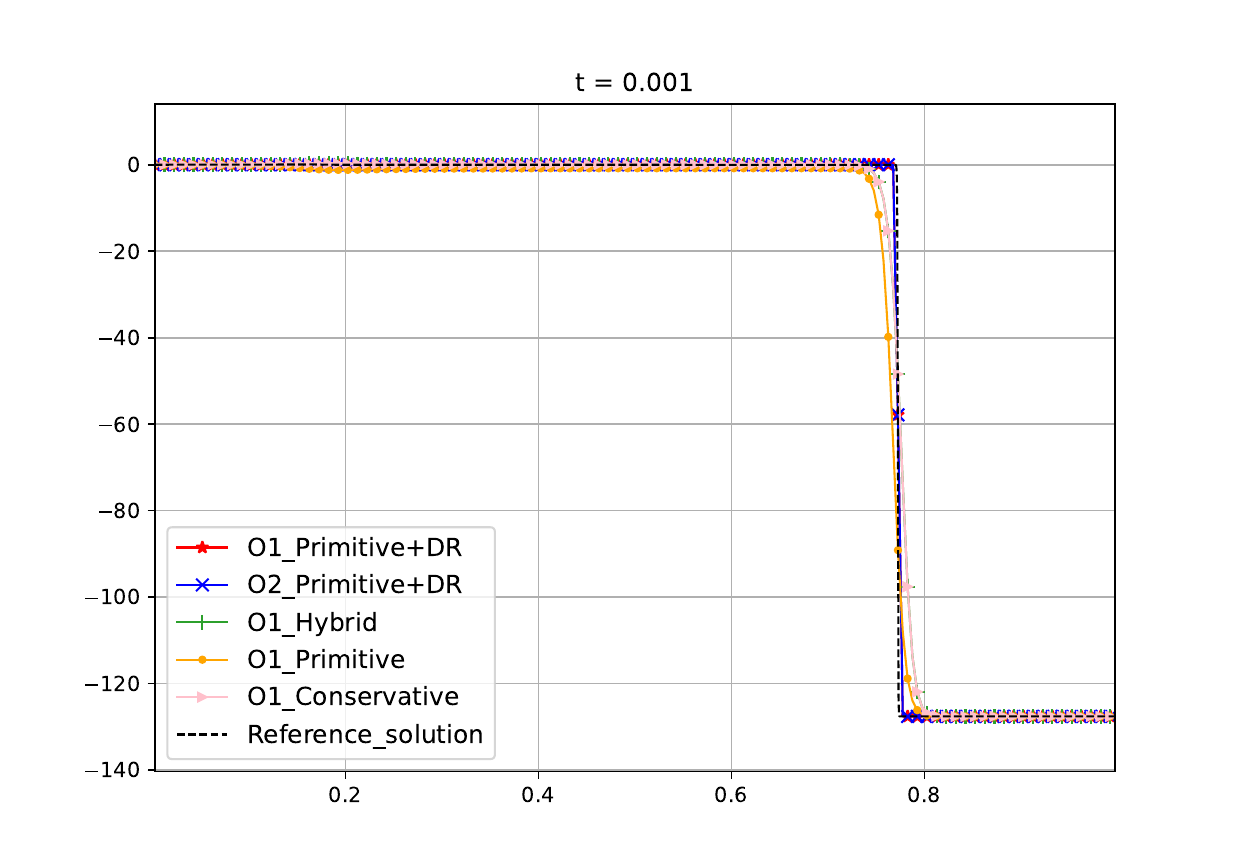}
    }
    \subfloat[Variable $u$. Zoom shock]{
        \includegraphics[width=0.30\textwidth]{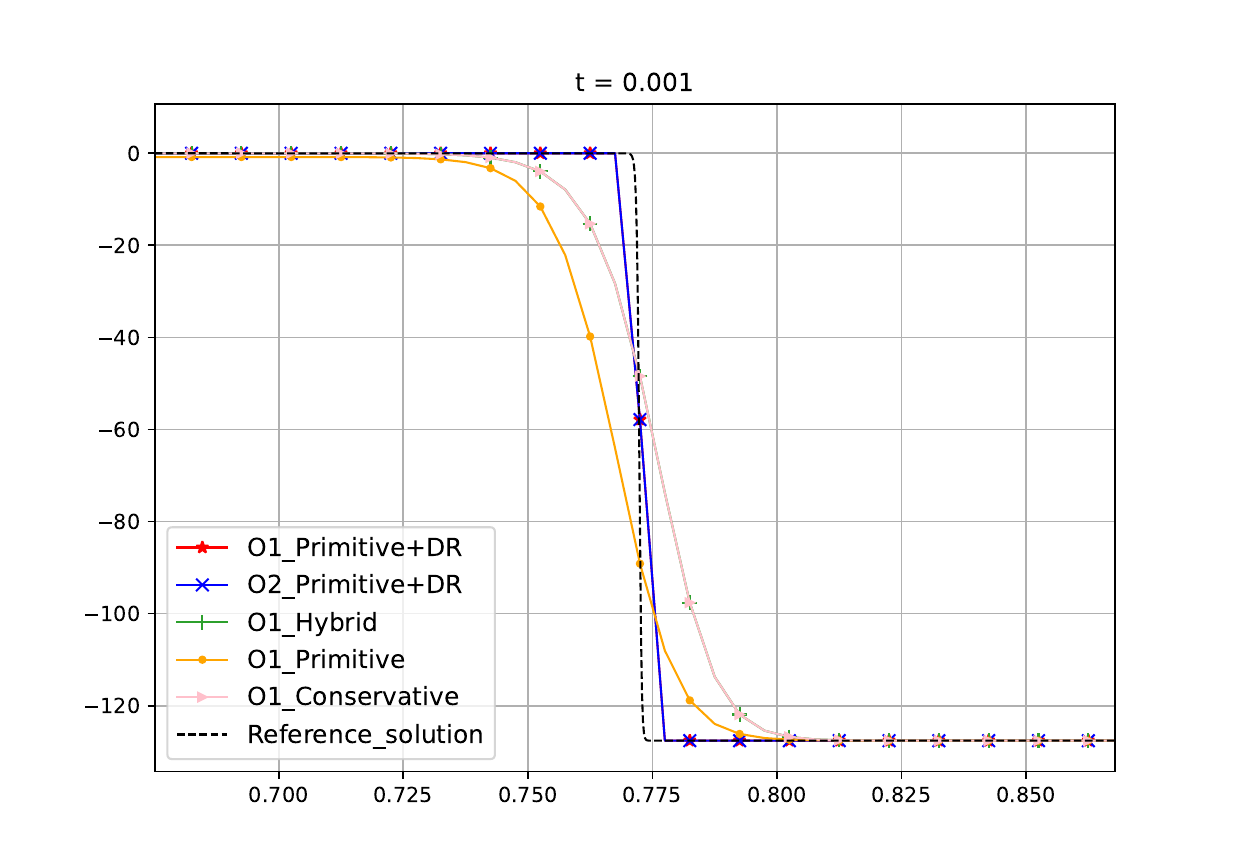}
    }
    \subfloat[Variable $u$. Zoom left state]{
        \includegraphics[width=0.30\textwidth]{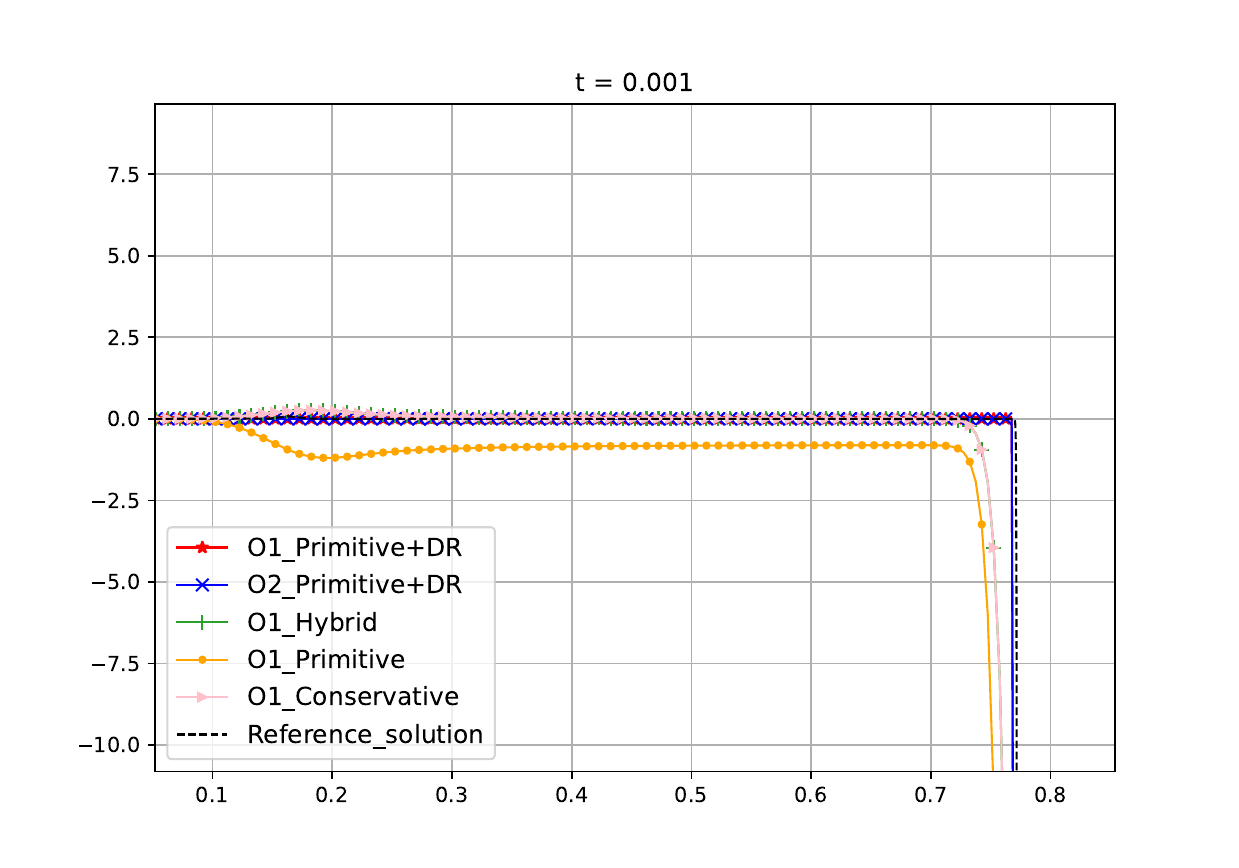}
    }

    \subfloat[Variable $p$]{
        \includegraphics[width=0.30\textwidth]{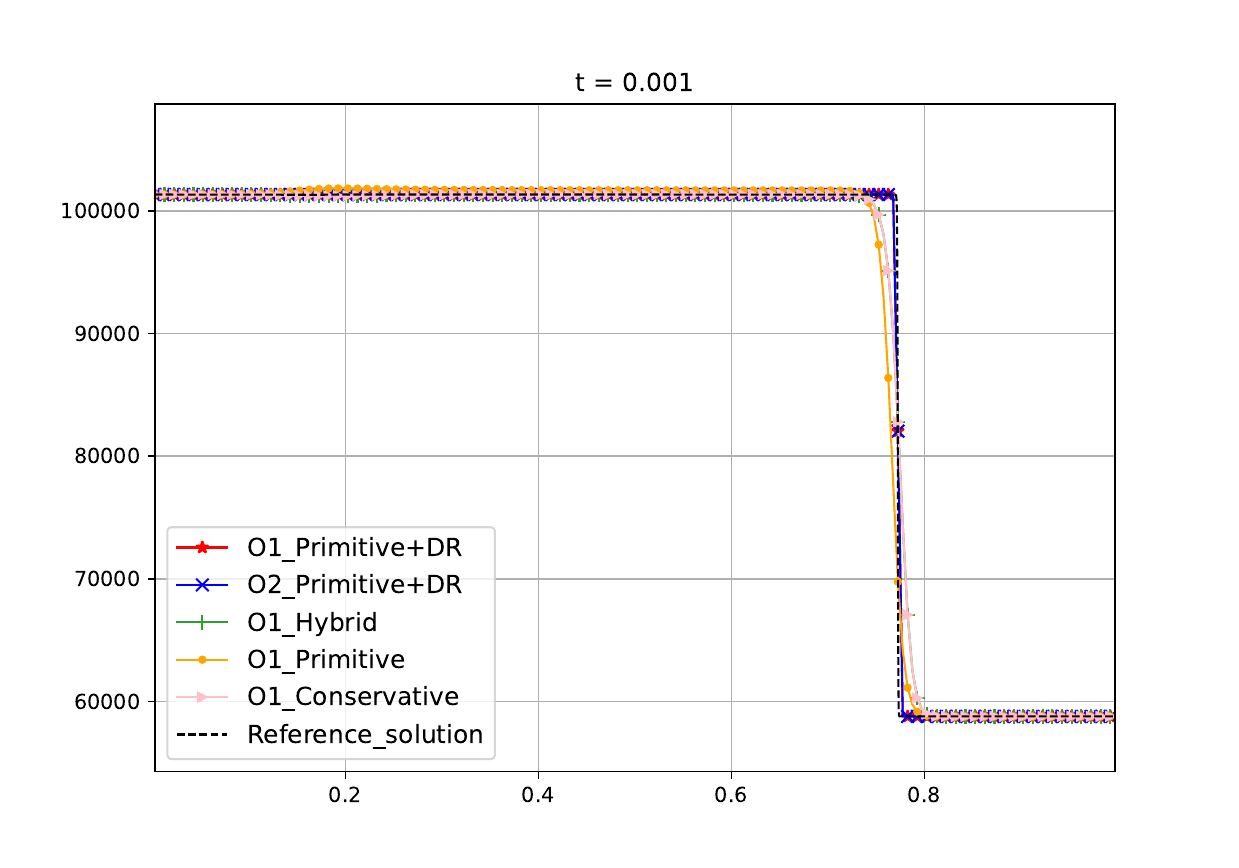}
    }
    \subfloat[Variable $p$. Zoom shock]{
        \includegraphics[width=0.30\textwidth]{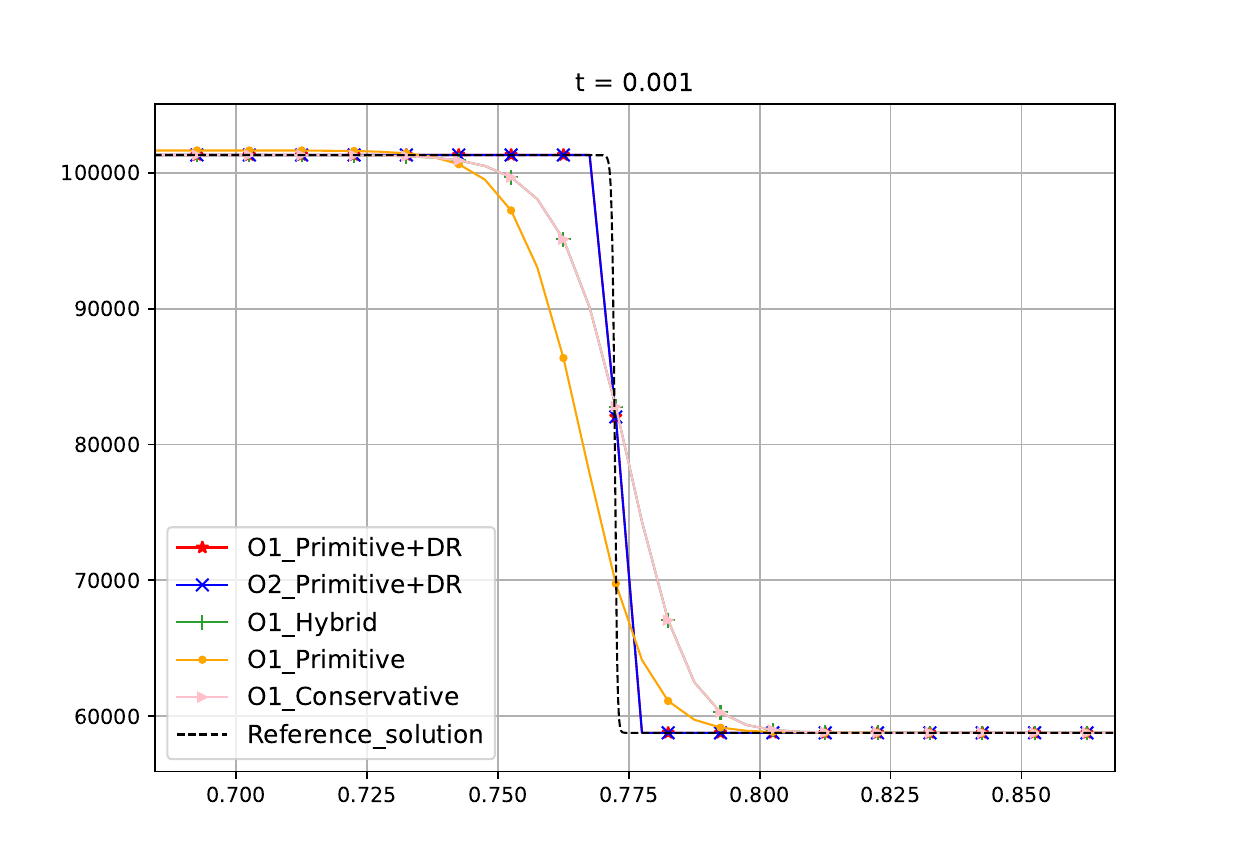}
    }
    \subfloat[Variable $p$. Zoom left state]{
        \includegraphics[width=0.30\textwidth]{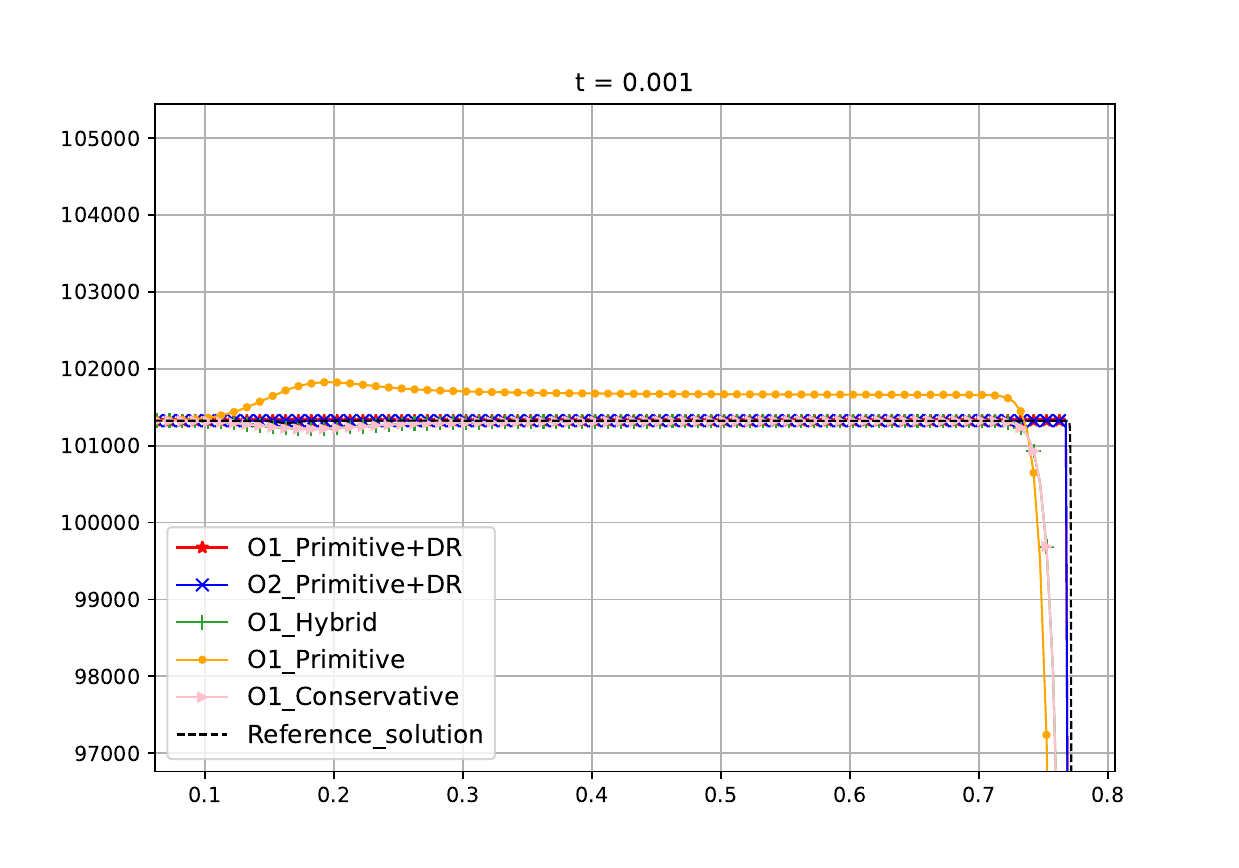}
    }
    \\
    \subfloat[Variable $c$]{
        \includegraphics[width=0.30\textwidth]{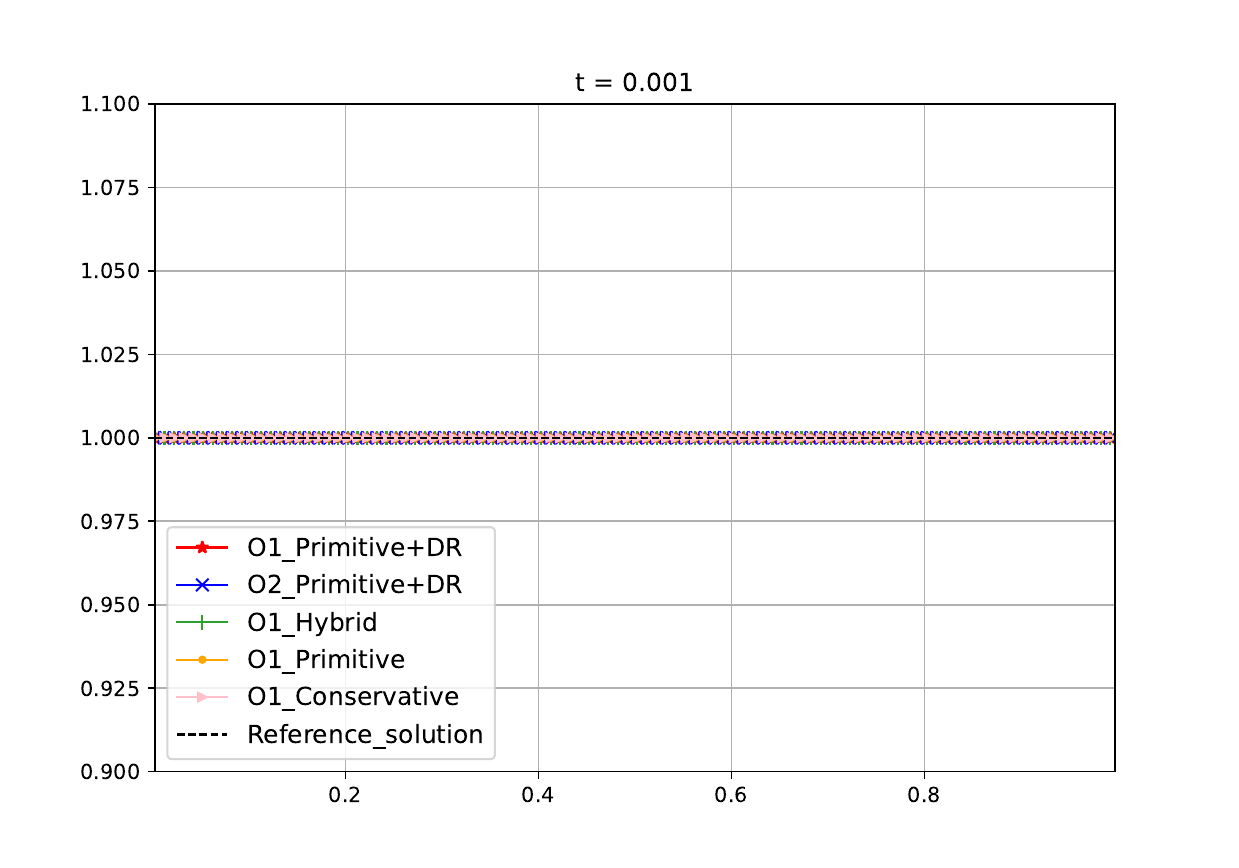}
    }

    \caption{Numerical results for Test: Isolated shock.}
    \label{Test_Isolated_shock}
\end{figure}

 Only one material (air $\gamma_p = 1.4$) is considered in this test. See Figure \ref{Test_Isolated_shock} for numerical results. In this case all the schemes capture the correct shock wave but O1\_Primitive converge as expected. In this case, the numerical results corresponding to O1\_Conservative and O1\_Hybrid coincide. Observe that the use of the discontinuous reconstruction makes the numerical method to capture the isolated shock wave much more sharply.

\begin{figure}[h]
    \centering

    \subfloat[Variable $\rho$]{
        \includegraphics[width=0.23\textwidth]{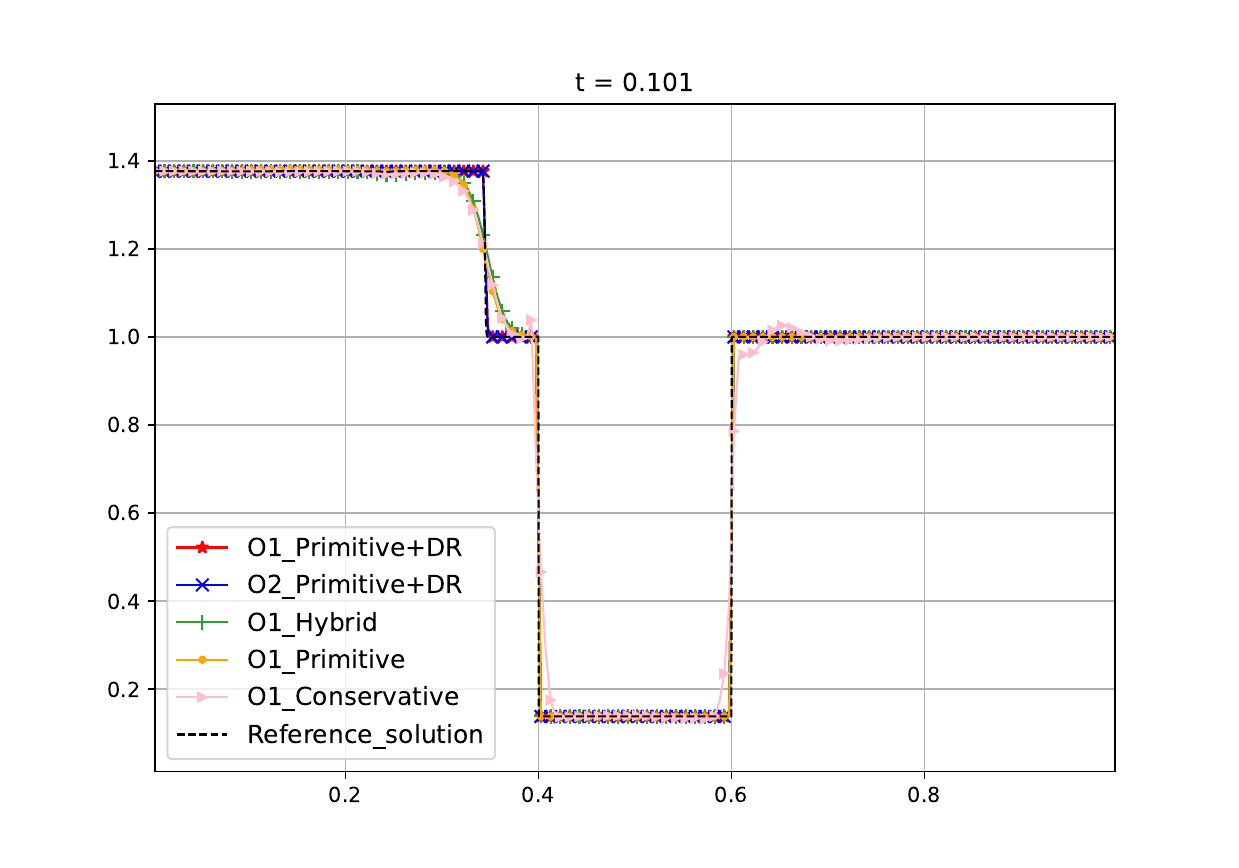}
    }
    \subfloat[Variable $\rho$. Zoom]{
        \includegraphics[width=0.23\textwidth]{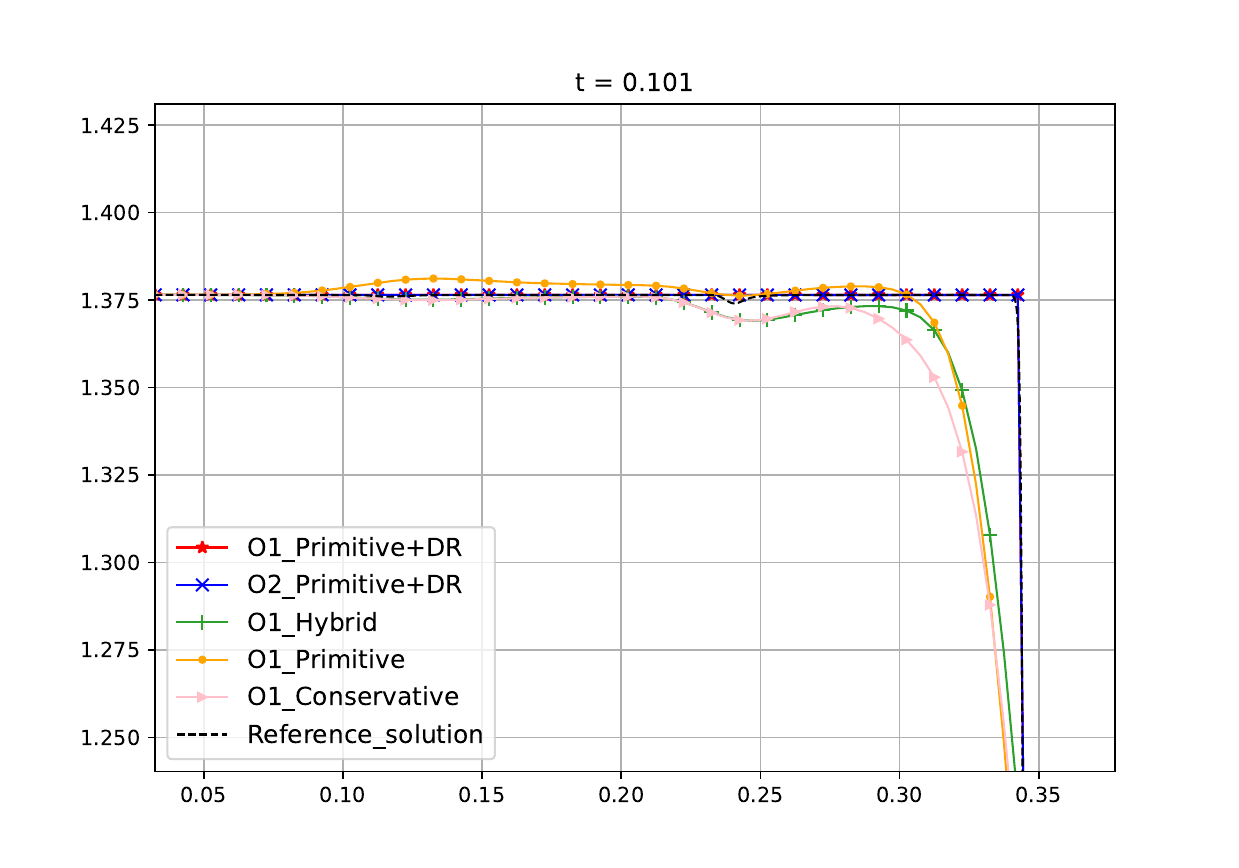}
    }
    \subfloat[Variable $\rho$]{
        \includegraphics[width=0.23\textwidth]{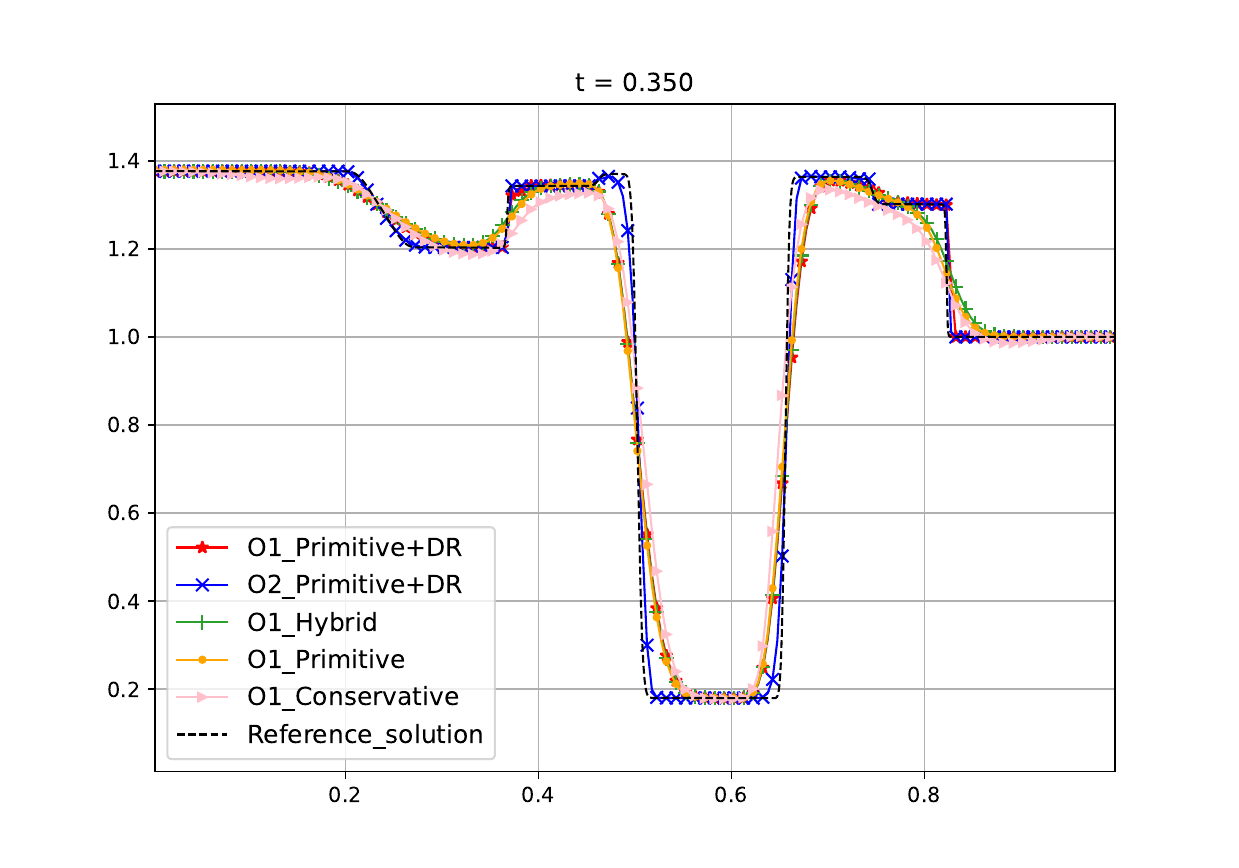}
    }
    \subfloat[Variable $\rho$. Zoom]{
        \includegraphics[width=0.23\textwidth]{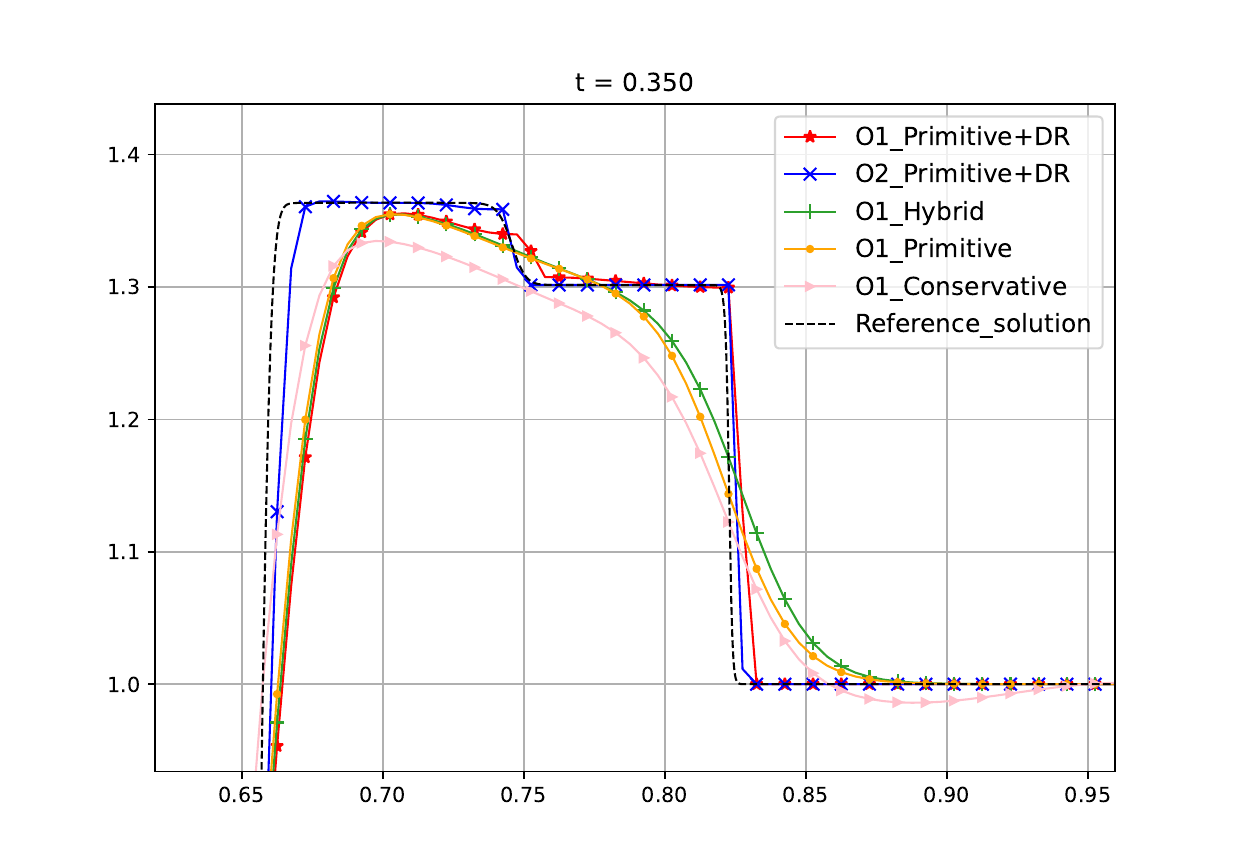}
    }

    \subfloat[Variable $u$]{
        \includegraphics[width=0.23\textwidth]{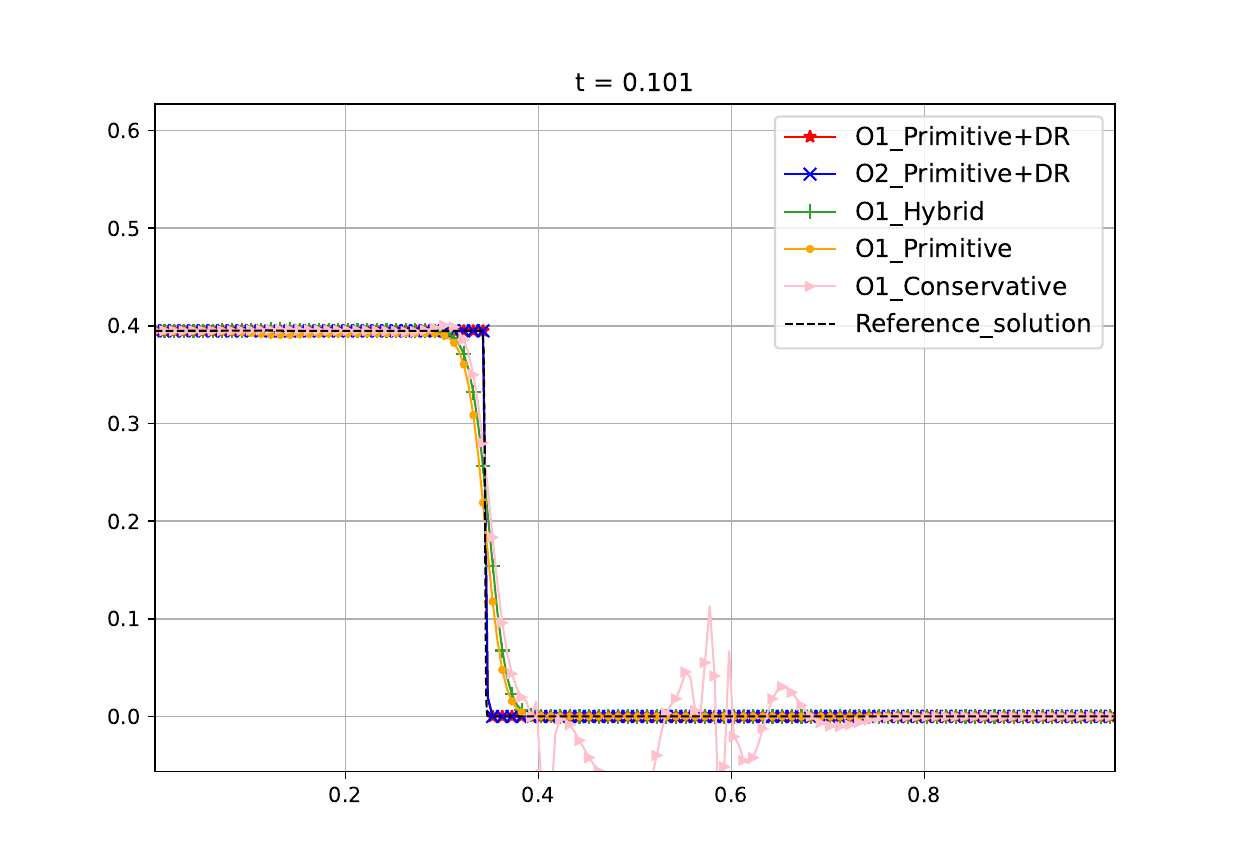}
    }
    \subfloat[Variable $u$. Zoom]{
        \includegraphics[width=0.23\textwidth]{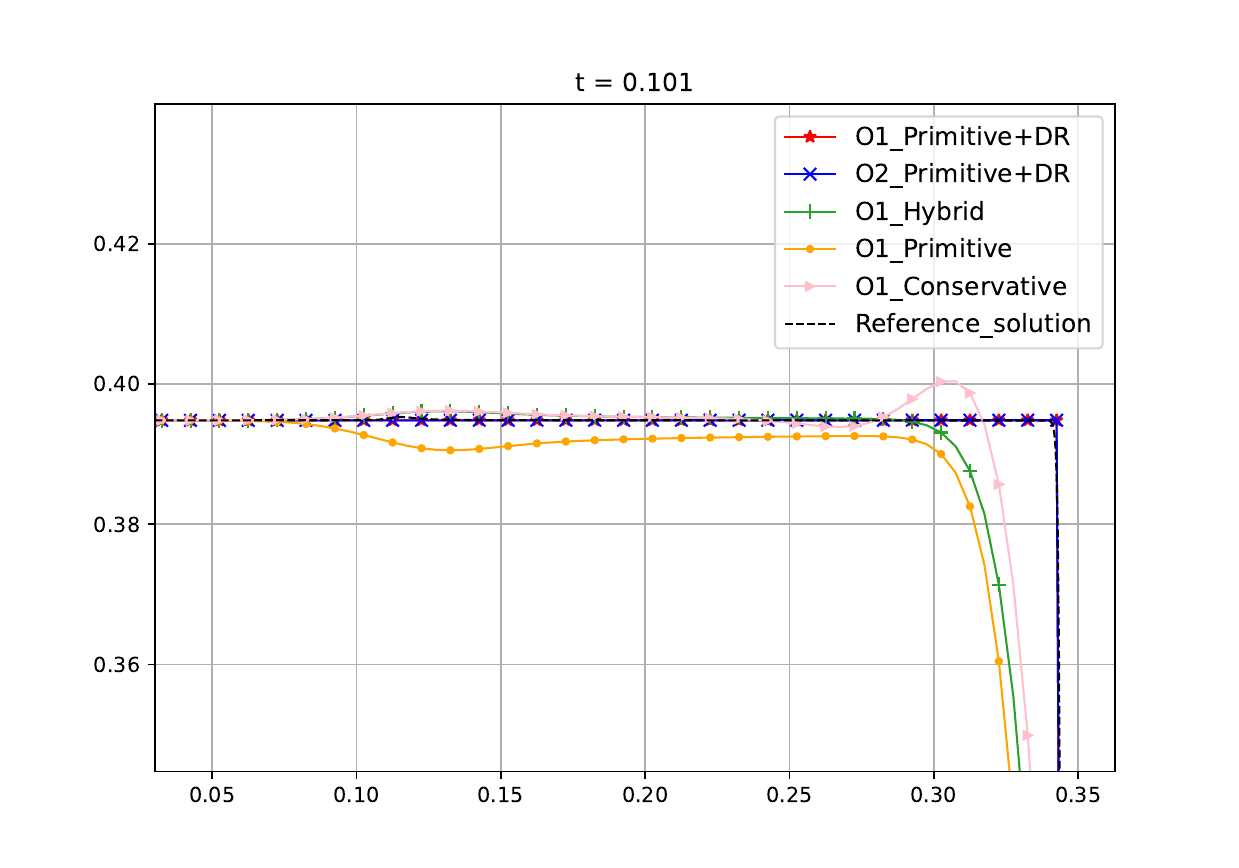}
    }
    \subfloat[Variable $u$]{
        \includegraphics[width=0.23\textwidth]{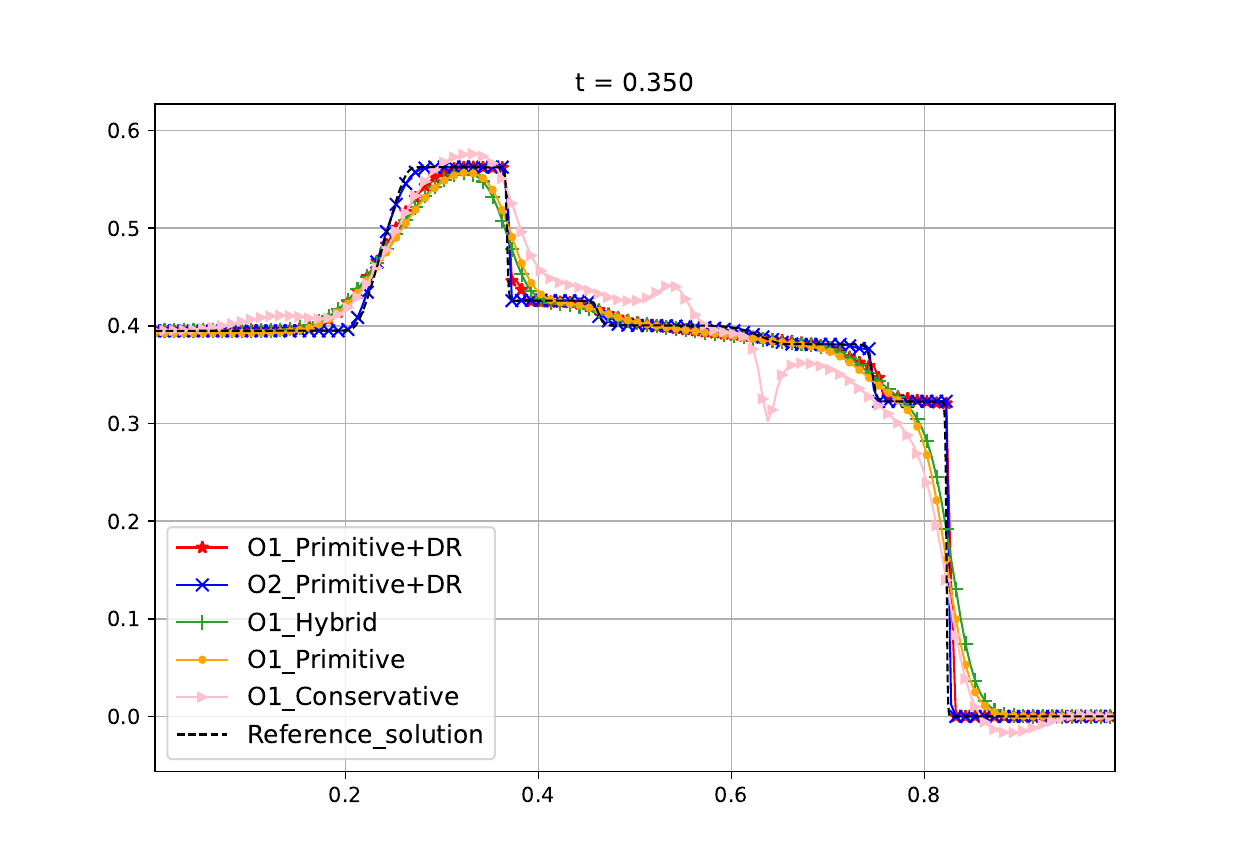}
    }
    \subfloat[Variable $u$. Zoom]{
        \includegraphics[width=0.23\textwidth]{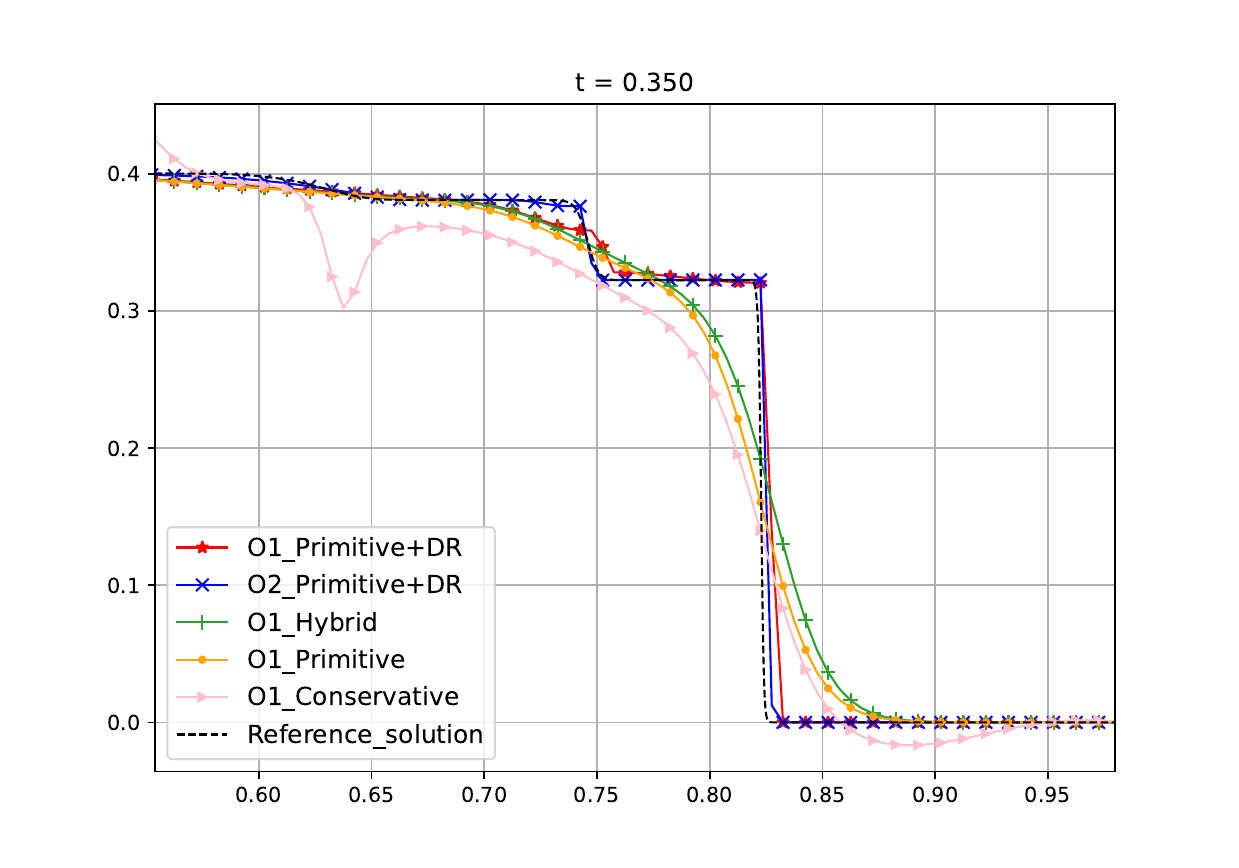}
    }

    \subfloat[Variable $p$]{
        \includegraphics[width=0.23\textwidth]{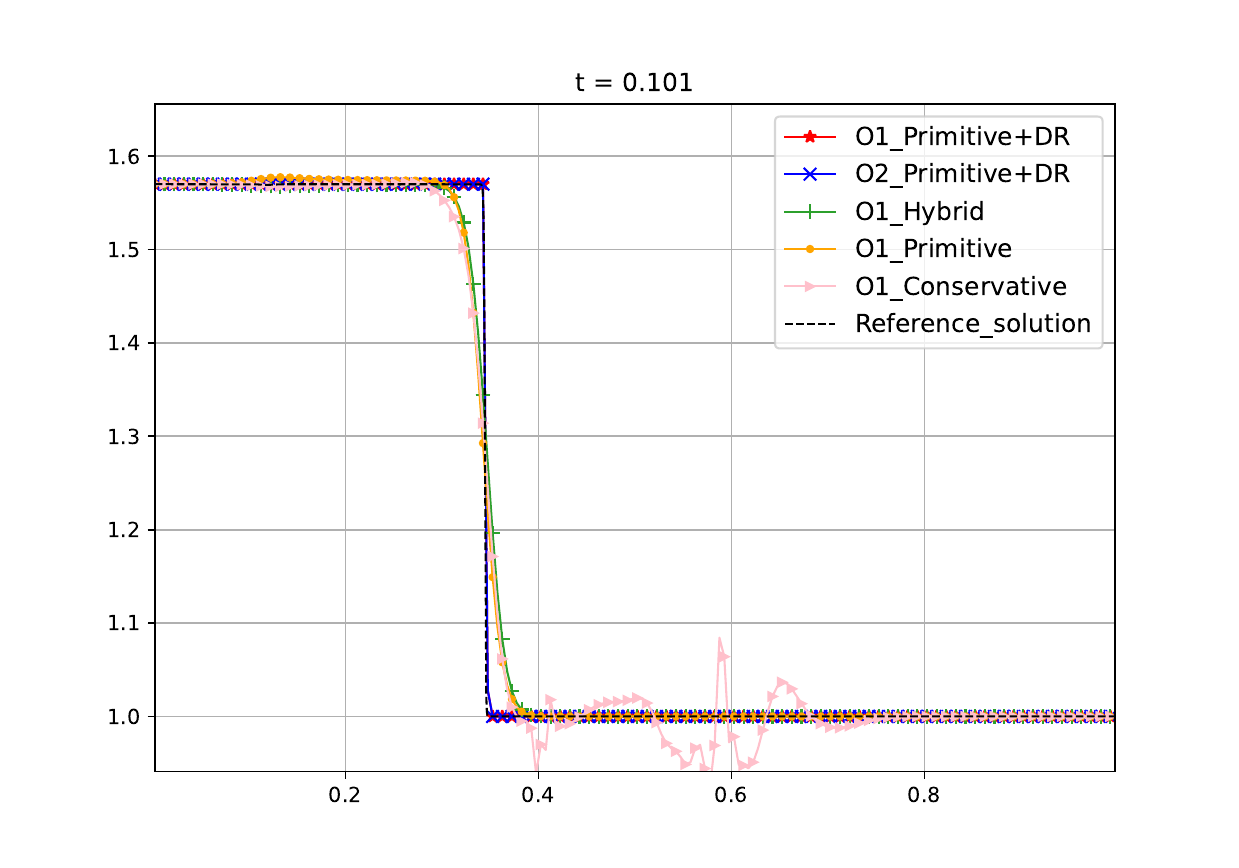}
    }
    \subfloat[Variable $p$. Zoom]{
        \includegraphics[width=0.23\textwidth]{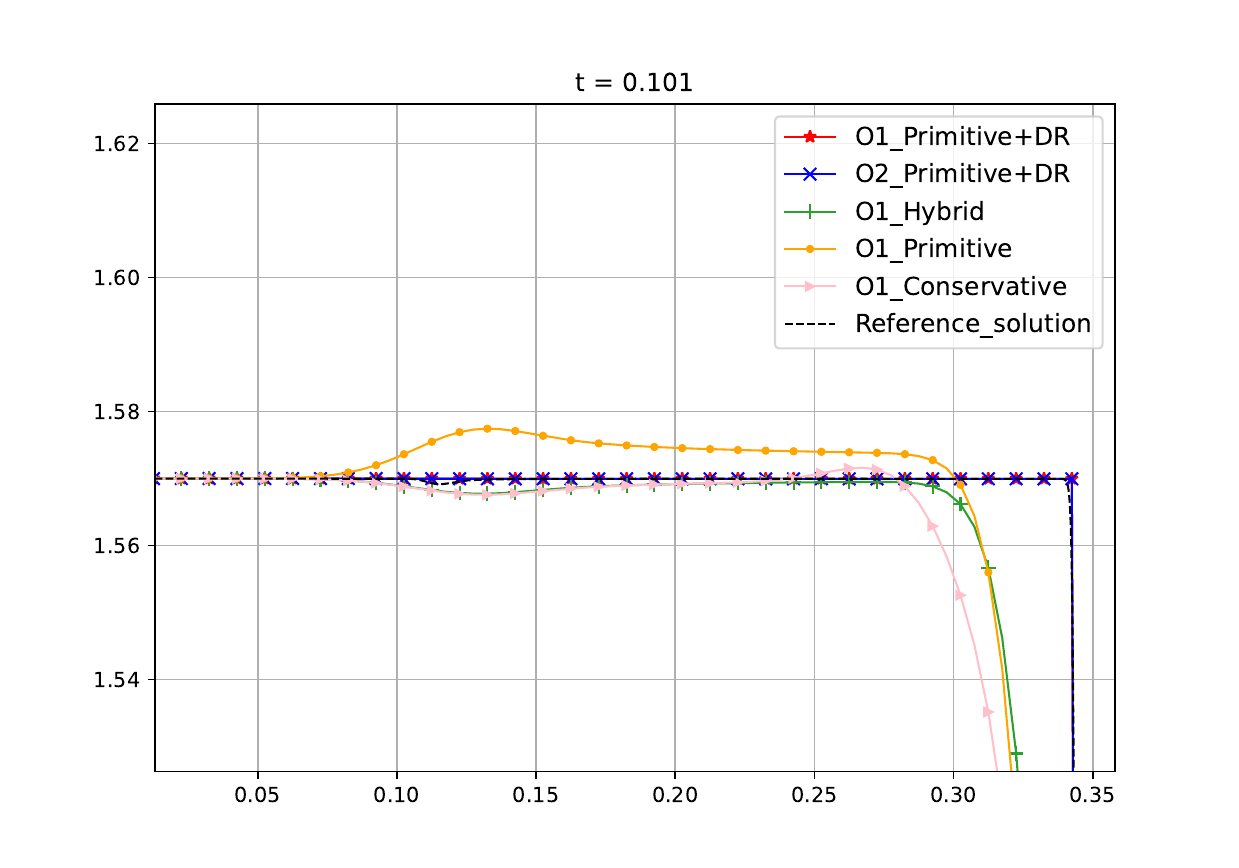}
    }
    \subfloat[Variable $p$]{
        \includegraphics[width=0.23\textwidth]{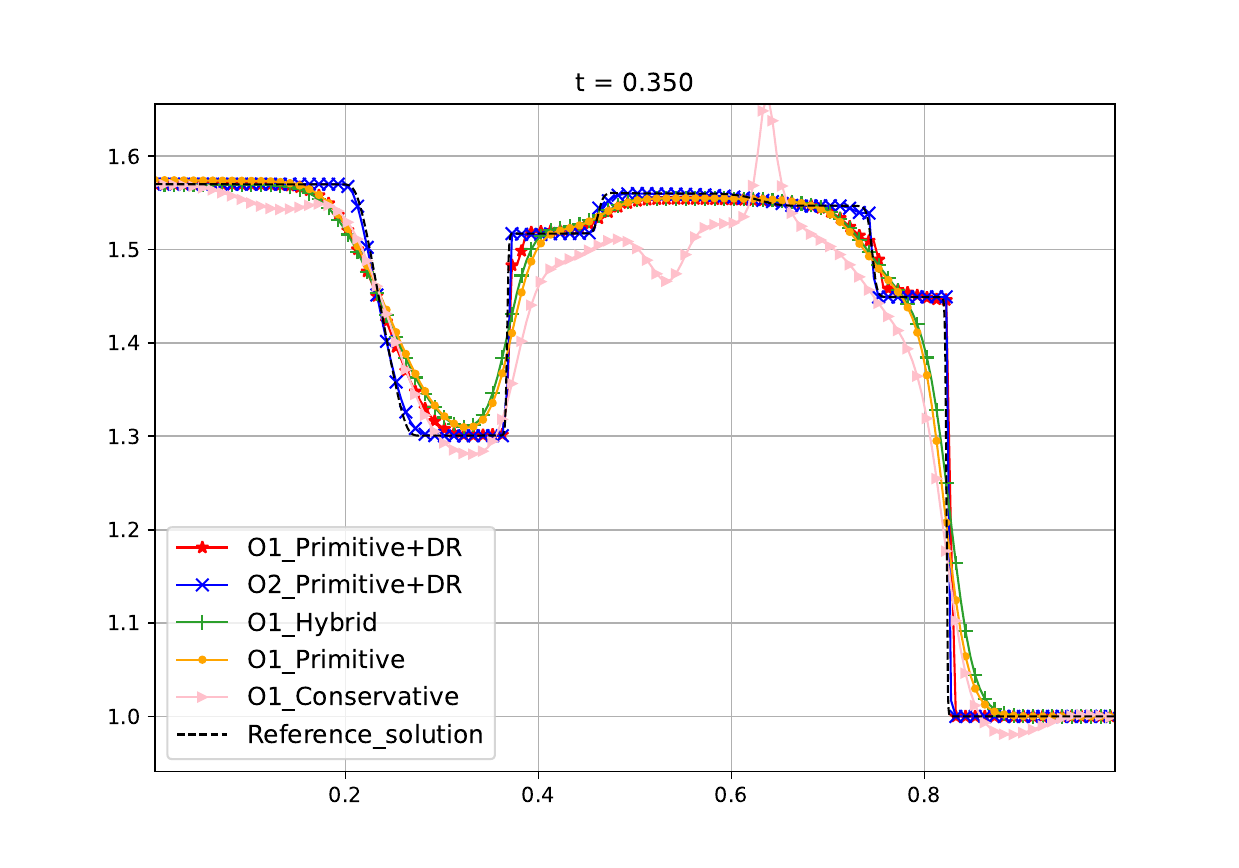}
    }
    \subfloat[Variable $p$. Zoom]{
        \includegraphics[width=0.23\textwidth]{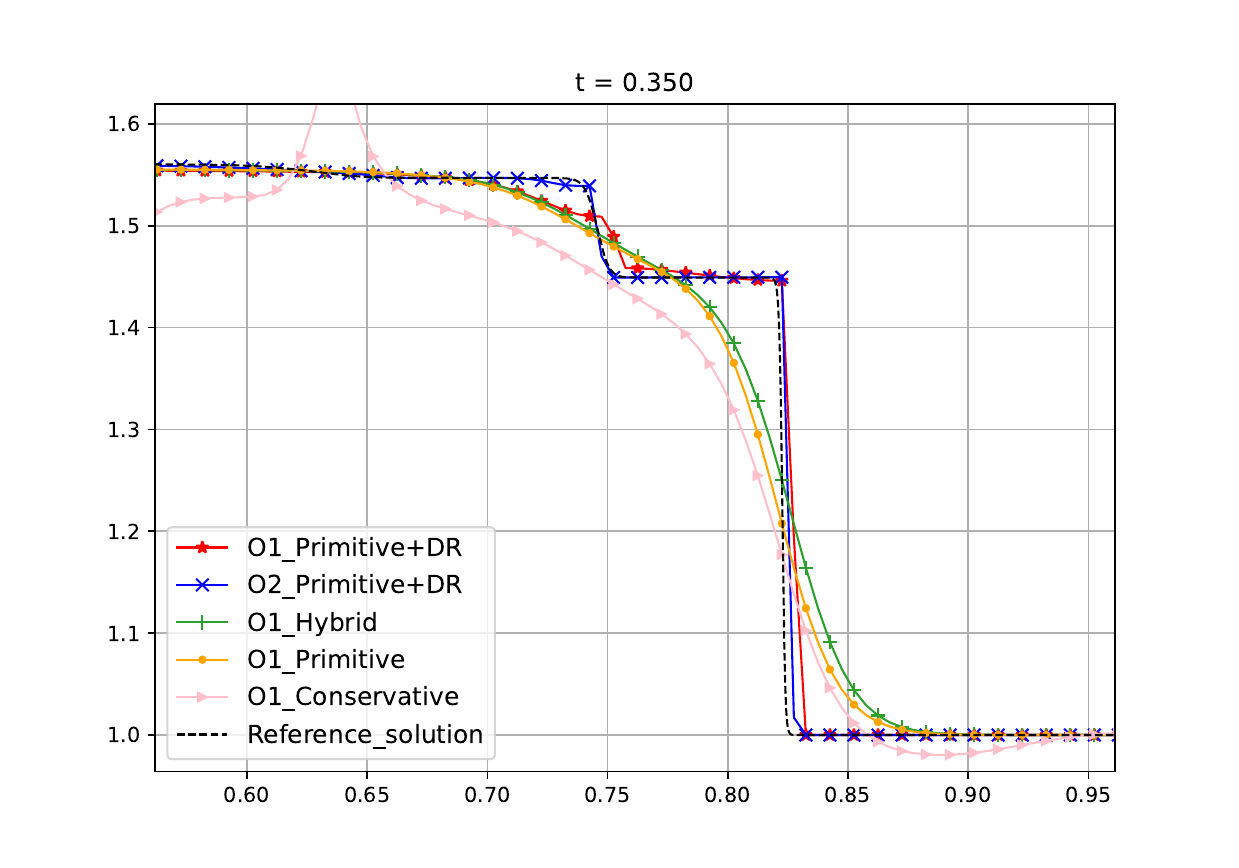}
    }
    \\
    \subfloat[Variable $c$]{
        \includegraphics[width=0.3\textwidth]{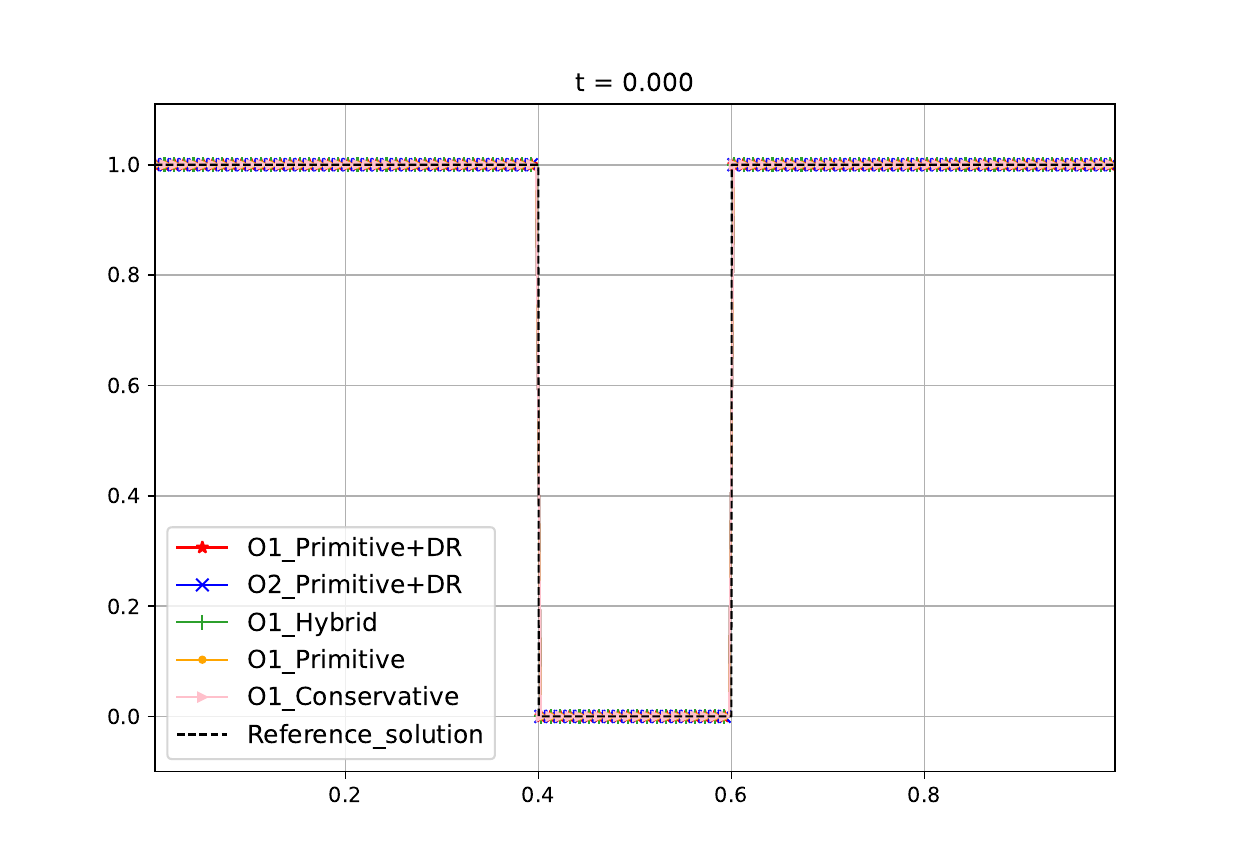}
    }
    \subfloat[Variable $c$]{
        \includegraphics[width=0.3\textwidth]{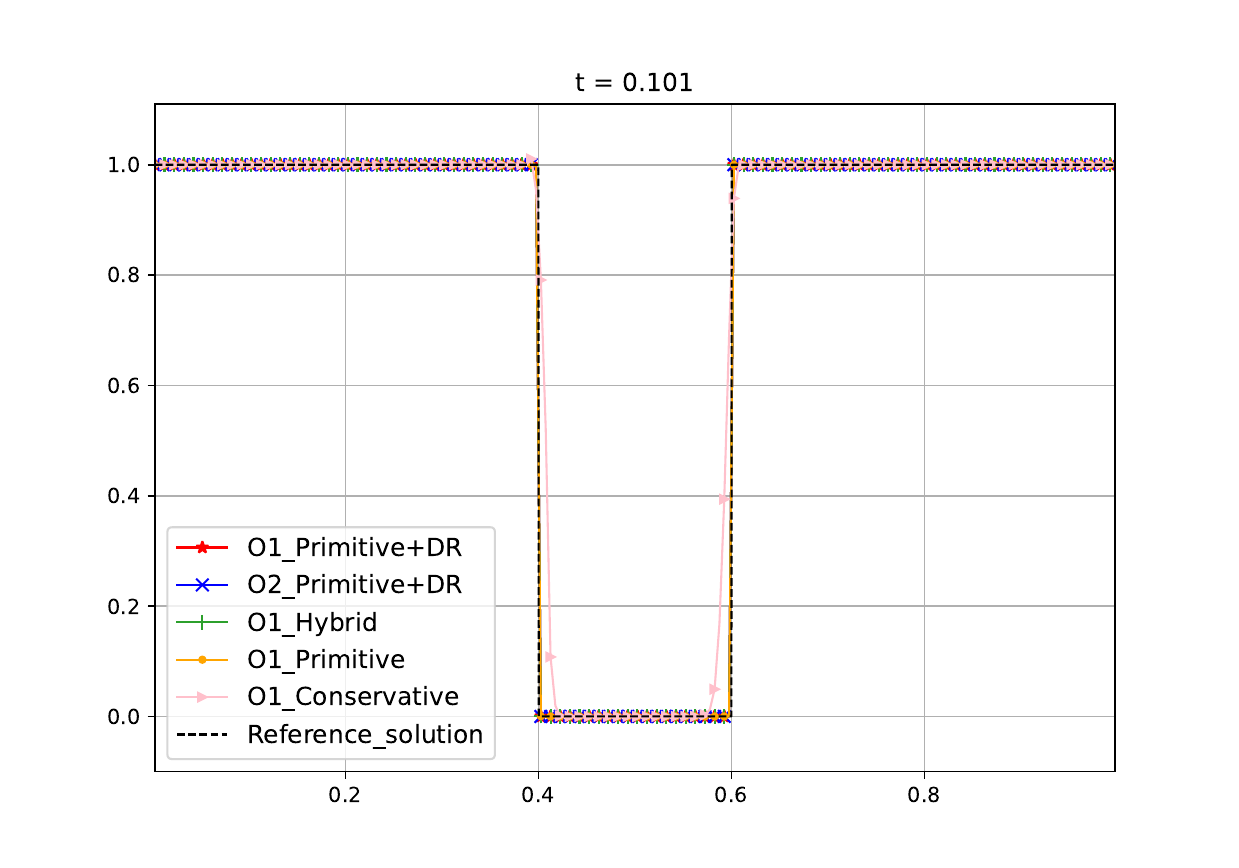}
    }
    \subfloat[Variable $c$]{
        \includegraphics[width=0.3\textwidth]{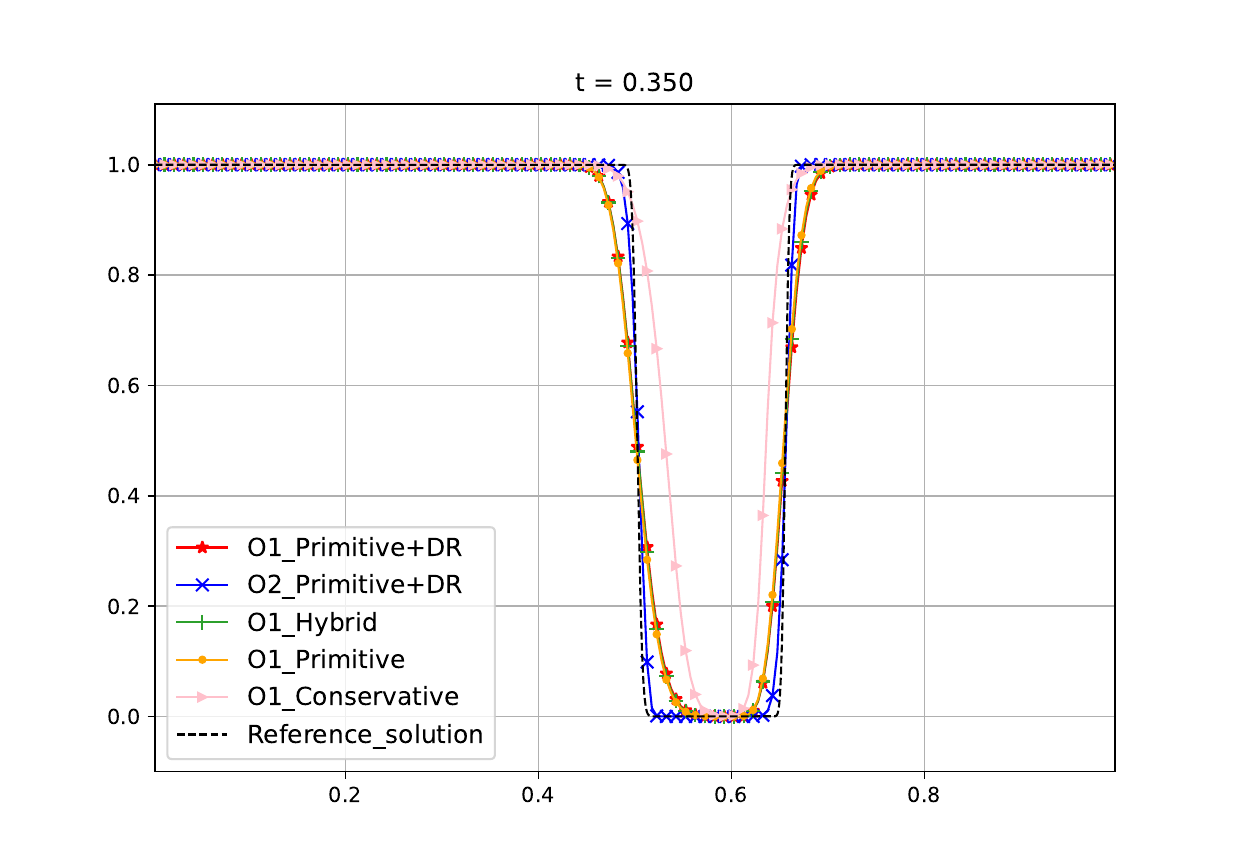}
    }

    \caption{Numerical results for Test: Quirk and Karni shock helium-bubble interaction.}
    \label{Test_Karni}
\end{figure}

\begin{figure}[h]
    \centering

    \subfloat[Variable $\rho$]{
        \includegraphics[width=0.45\textwidth]{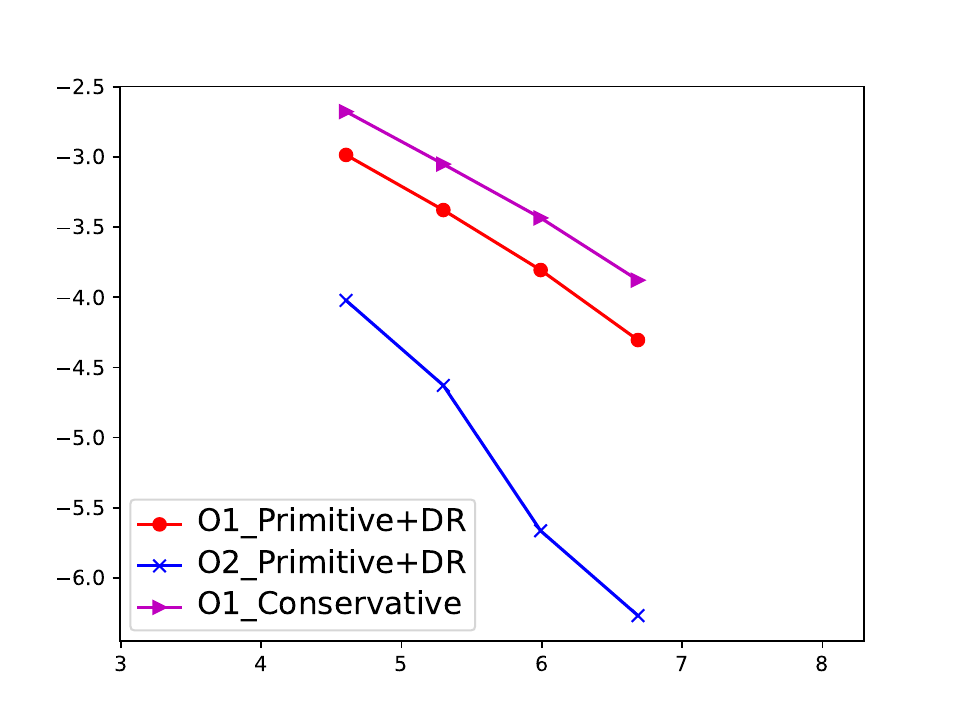}
    }
    \subfloat[Variable $u$]{
        \includegraphics[width=0.45\textwidth]{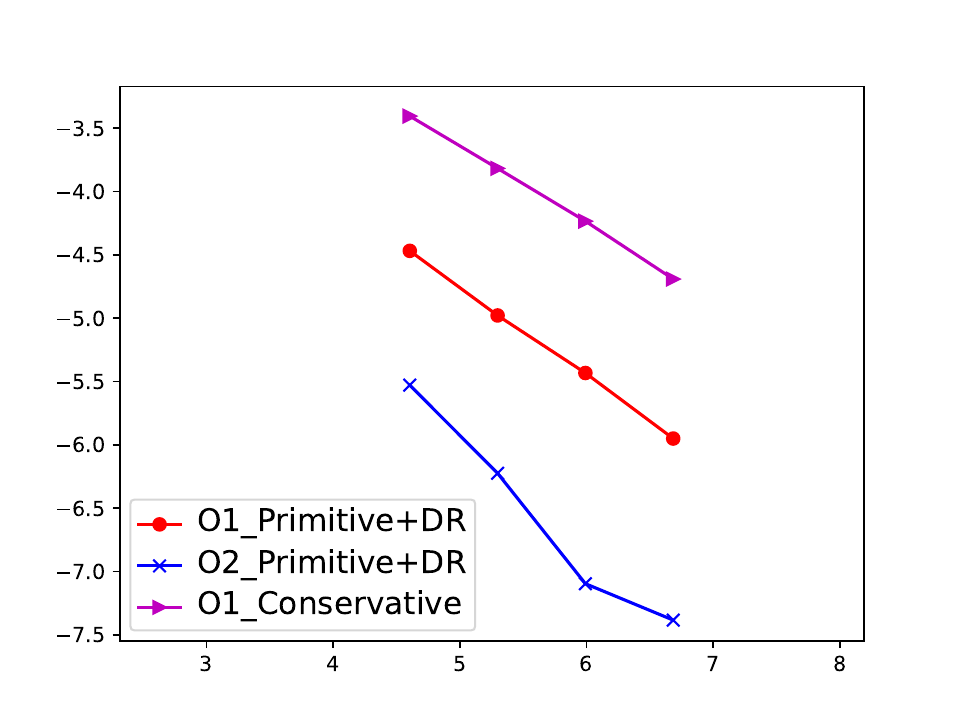}
    }
    \\

    \subfloat[Variable $p$]{
        \includegraphics[width=0.45\textwidth]{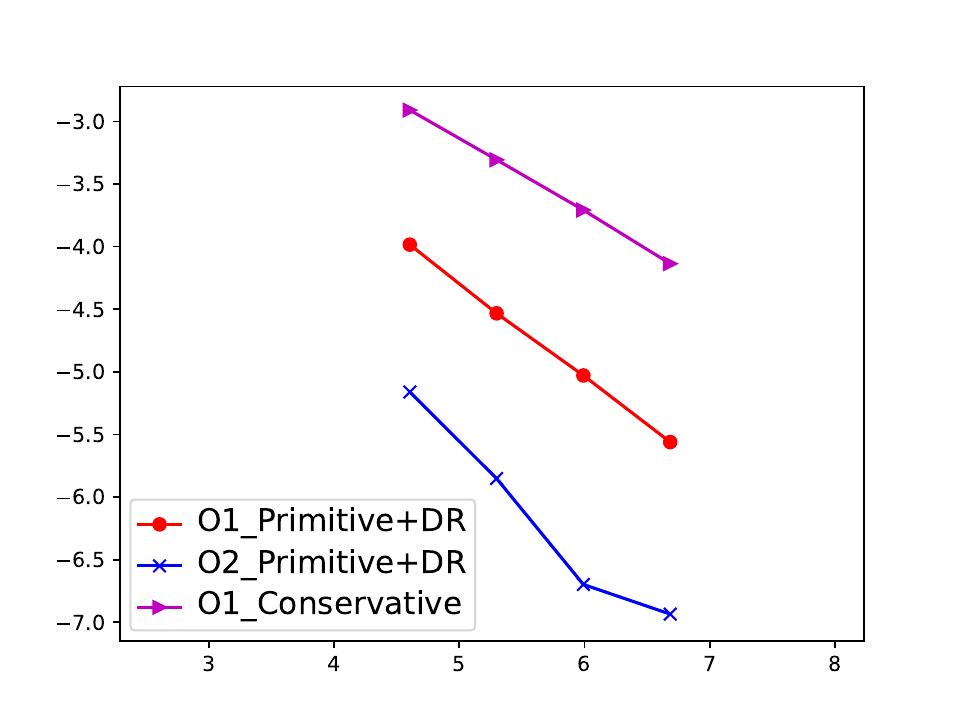}
    }
    \subfloat[Variable $c$]{
        \includegraphics[width=0.45\textwidth]{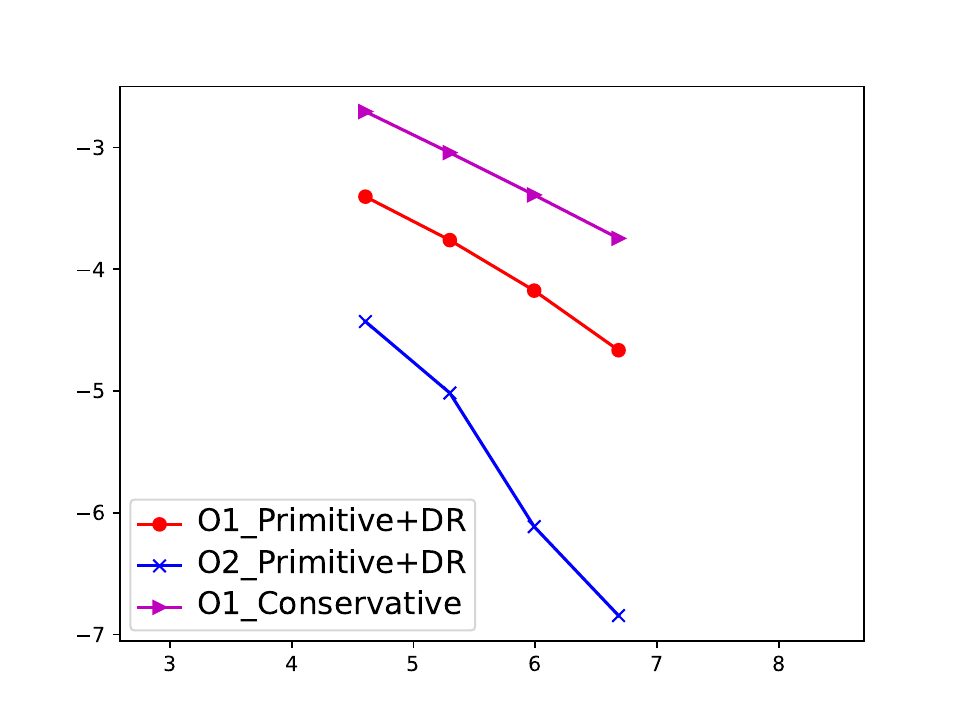}
    }

    \caption{Convergence for Test \textit{Quirk and Karni shock helium-bubble interaction}: Numerical errors obtained in log scale for O1\_Primitive+DR, O2\_Primitive+DR, O1\_Conservative at time $t=0.35$.}
    \label{Test_Error_rates_Quirk}
\end{figure}

 \paragraph{Quirk and Karni shock helium-bubble interaction}
  This test is taken from Quirk and Karni in \cite{Quirk1996129}. The initial condition is given by
 $$
 \begin{aligned}
 \bw_{0}(x) & =  (\rho_0(x), u_0(x), v_0(x), p_0(x), c_0(x))^{T} \\ & =  \begin{cases}
     (1.3765, 0.3948, 0, 1.57, 1)^{T} & \text{if $x<0.2$,}\\
     (1, 0, 0, 1, 1)^{T} & \text{if $0.2\leq x<0.4$,}\\
     (0.138, 0, 0, 1, 0)^{T} & \text{if $0.4\leq x\leq 0.6$,}\\
     (1, 0, 0, 1, 1)^{T} & \text{if $x>0.6$.}
 \end{cases}  
 \end{aligned}
 $$
 This test corresponds to the interaction between a shock and a helium bubble. Figure \ref{Test_Karni} shoes the numerical results. We observe that the poorest results are produced by the fully conservative scheme, due to spurious oscillations arising in the vicinity of contact discontinuities. The standard path-conservative method, when used without the in-cell reconstruction operator, fails to correctly capture the plateau near the left and shocks. In the remaining cases, all schemes appear to converge to the correct solution; however, those that incorporate the in-cell discontinuous reconstruction operator exhibit superior performance in shock regions. Moreover, as expected, the second-order scheme gives better results in rarefaction waves. In \ref{tab:Error_TestQuirk} we observe the errors obtained with O1\_Primitive+DR, O2\_Primitive+DR, O1\_Conservative. All seem to converge to the correct solution, but together with Figure \ref{Test_Error_rates_Quirk}, the convergence is much faster in the case of the in-cell discontinuous reconstruction path-conservative methods.

\begin{table}[!h]
  	\centering
  	\begin{tabular}{|c|c|c|c|c|c|c|c|c|}
  	    \hline
         \multicolumn{9}{|c|}{O1\_Primitive+DR} \\
 \hline
  		N & $||\Delta \rho||_1$ & Order & $||\Delta u||_1$ & Order & $||\Delta p||_1$ & Order & $||\Delta c||_1$ & Order\\  
  		\hline 
   100 & 5.05e-02 & 0.00 & 1.15e-02 & 0.00 & 1.86e-02 & 0.00 & 3.32e-02 & 0.00 \\ 
   200 & 3.41e-02 & 0.57 & 6.88e-03 & 0.74 & 1.08e-02 & 0.79 & 2.33e-02 & 0.52 \\ 
   400 & 2.22e-02 & 0.62 & 4.37e-03 & 0.66 & 6.55e-03 & 0.72 & 1.54e-02 & 0.60 \\ 
   800 & 1.35e-02 & 0.72 & 2.61e-03 & 0.75 & 3.85e-03 & 0.77 & 9.42e-03 & 0.71 \\ 
  		\hline 
  		\hline
         \multicolumn{9}{|c|}{O2\_Primitive+DR} \\
 \hline
  		N & $||\Delta \rho||_1$ & Order & $||\Delta u||_1$ & Order & $||\Delta p||_1$ & Order & $||\Delta c||_1$ & Order\\  
  		\hline 
   100 & 1.79e-02 & 0.00 & 3.97e-03 & 0.00 & 5.74e-03 & 0.00 & 1.19e-02 & 0.00 \\ 
   200 & 9.77e-03 & 0.87 & 1.98e-03 & 1.00 & 2.87e-03 & 1.00 & 6.63e-03 & 0.85 \\ 
   400 & 3.47e-03 & 1.49 & 8.28e-04 & 1.26 & 1.23e-03 & 1.22 & 2.21e-03 & 1.58 \\ 
   800 & 1.90e-03 & 0.87 & 6.21e-04 & 0.41 & 9.74e-04 & 0.34 & 1.06e-03 & 1.05 \\
  		\hline
        \hline
         \multicolumn{9}{|c|}{O1\_Conservative} \\
 \hline
  		N & $||\Delta \rho||_1$ & Order & $||\Delta u||_1$ & Order & $||\Delta p||_1$ & Order & $||\Delta c||_1$ & Order\\  
  		\hline 
   100 & 6.87e-02 & 0.00 & 3.32e-02 & 0.00 & 5.45e-02 & 0.00 & 6.69e-02 & 0.00 \\ 
   200 & 4.73e-02 & 0.54 & 2.20e-02 & 0.60 & 3.66e-02 & 0.57 & 4.77e-02 & 0.49 \\ 
   400 & 3.22e-02 & 0.55 & 1.45e-02 & 0.60 & 2.45e-02 & 0.58 & 3.37e-02 & 0.50 \\ 
   800 & 2.07e-02 & 0.64 & 9.18e-03 & 0.66 & 1.60e-02 & 0.62 & 2.36e-02 & 0.51 \\
  		\hline
  	\end{tabular} 
	  	\caption{Test \textit{Quirk and Karni shock helium-bubble interaction}: convergence of O1\_Primitive+DR, O2\_Primitive+DR, O1\_Conservative: $L^{1}$ errors $||\Delta \cdot||_1$ at time $t=0.35$.}

  	\label{tab:Error_TestQuirk}
\end{table}

\begin{figure}[h]
    \centering

    \subfloat[Variable $\rho$]{
        \includegraphics[width=0.23\textwidth]{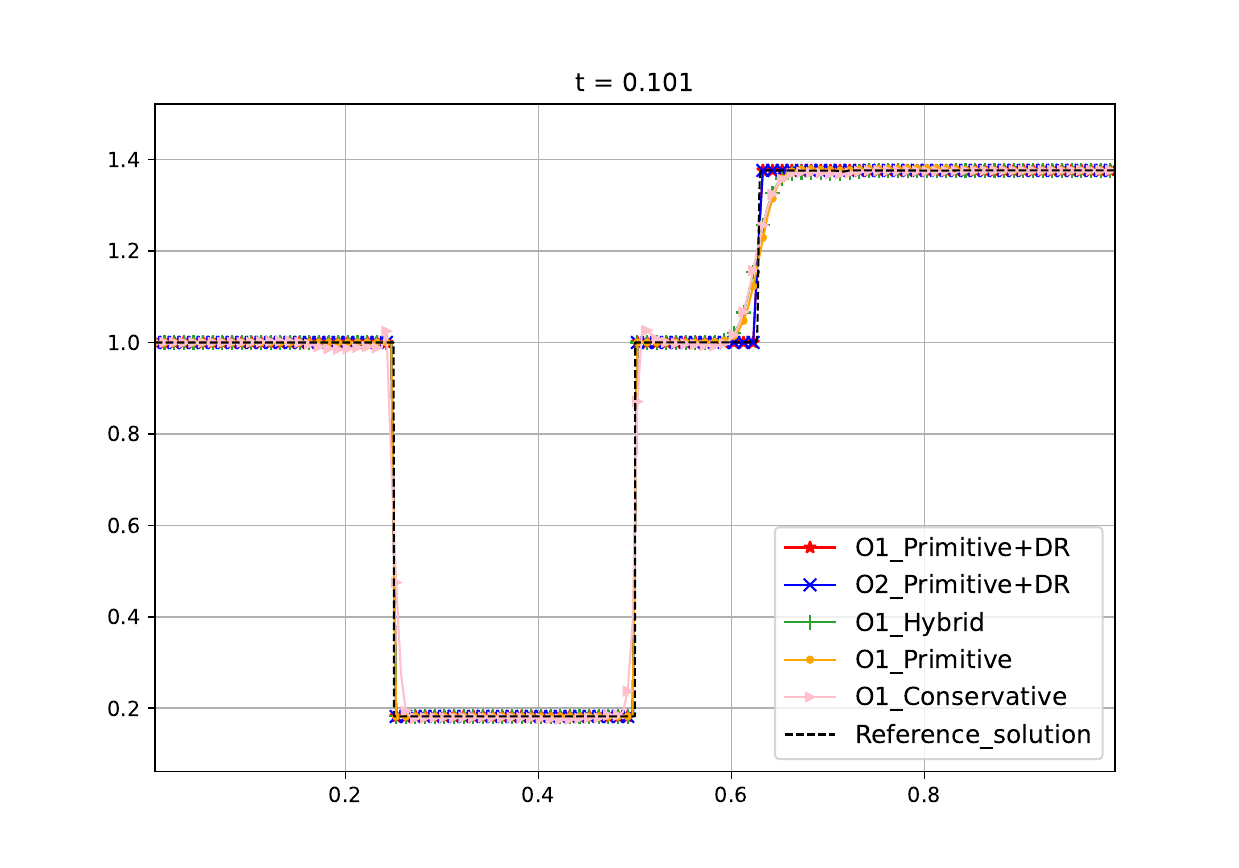}
    }
    \subfloat[Variable $\rho$. Zoom]{
        \includegraphics[width=0.23\textwidth]{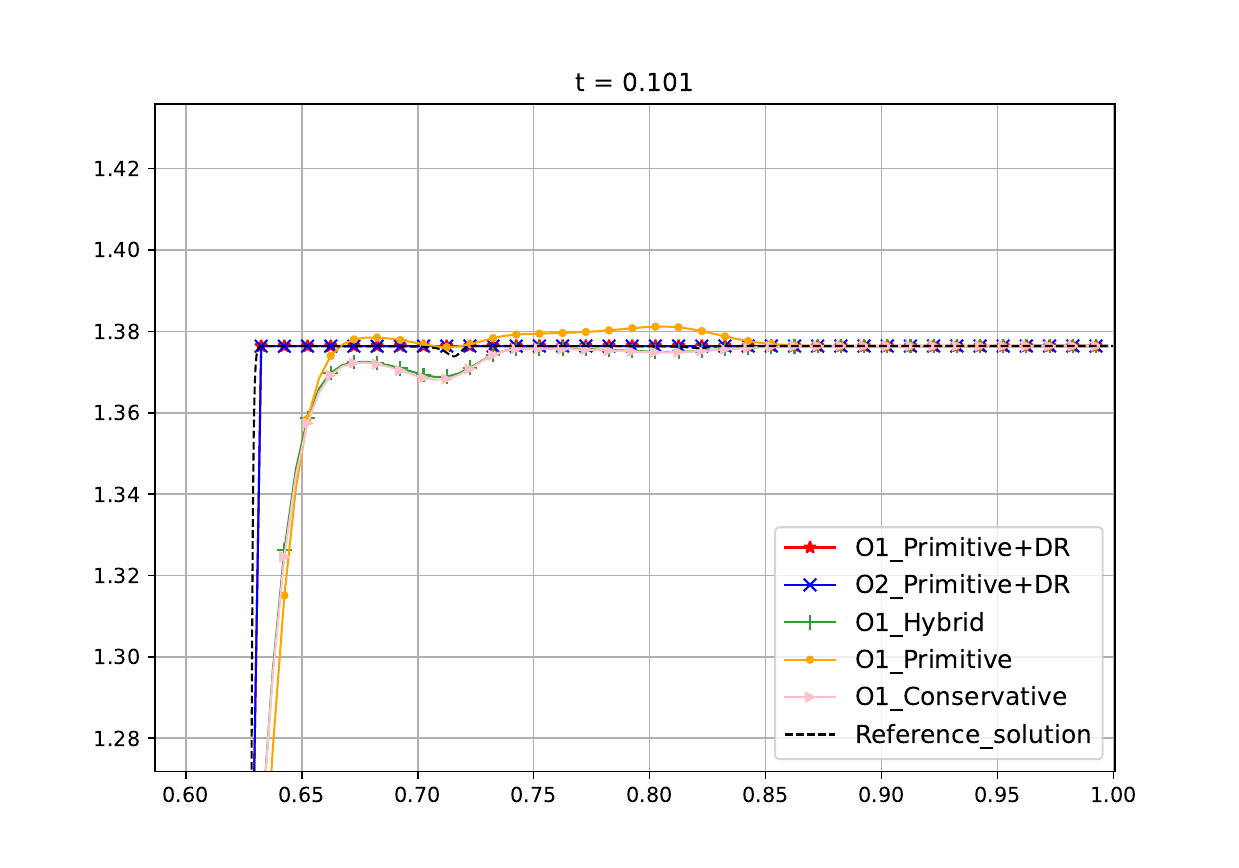}
    }
    \subfloat[Variable $\rho$]{
        \includegraphics[width=0.23\textwidth]{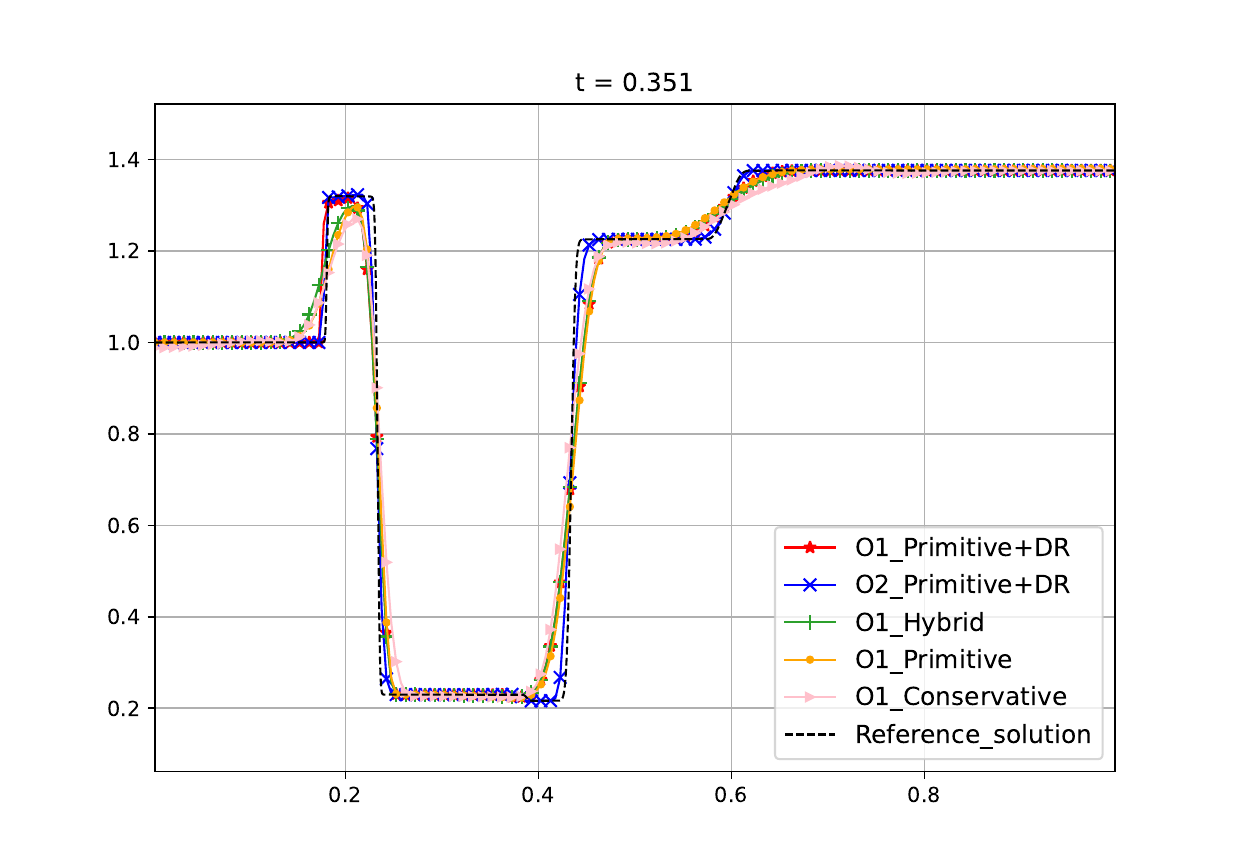}
    }
    \subfloat[Variable $\rho$. Zoom]{
        \includegraphics[width=0.23\textwidth]{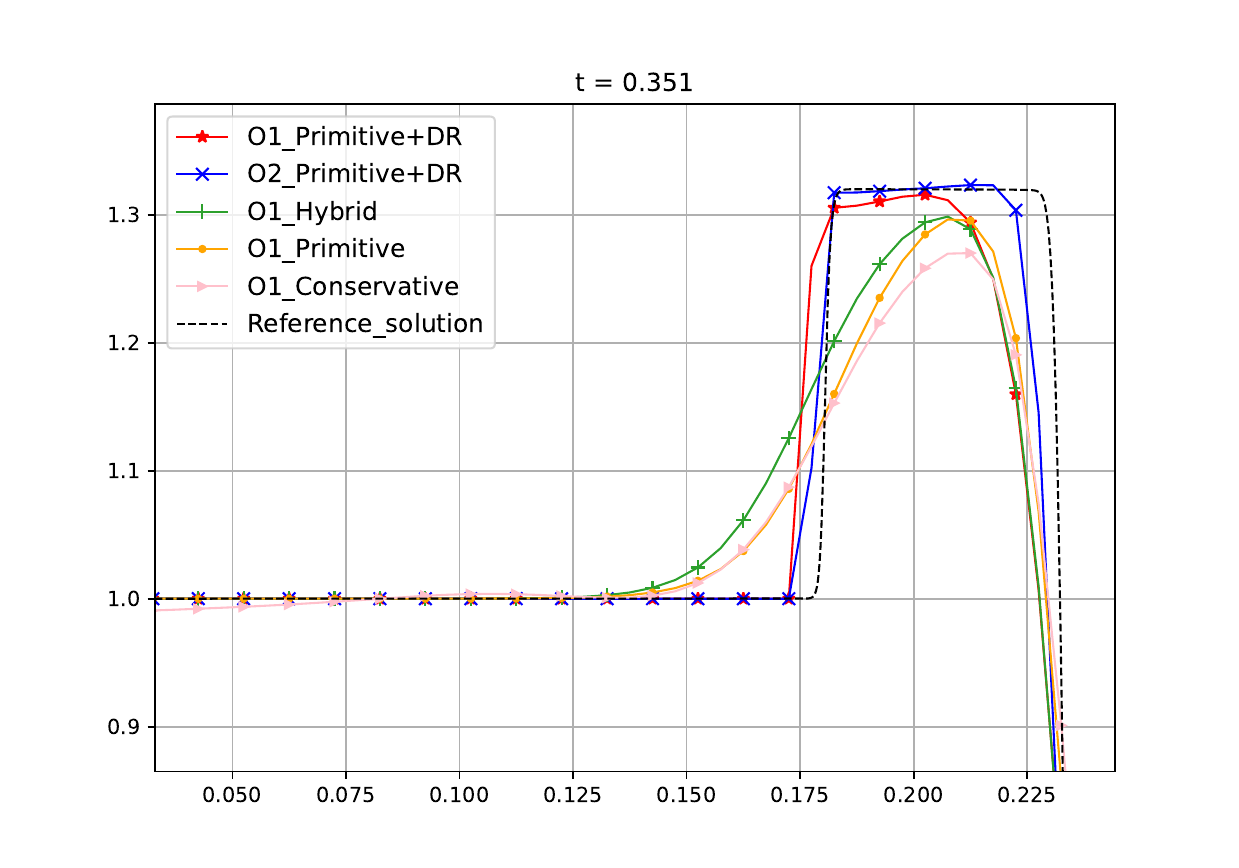}
    }

    \subfloat[Variable $u$]{
        \includegraphics[width=0.23\textwidth]{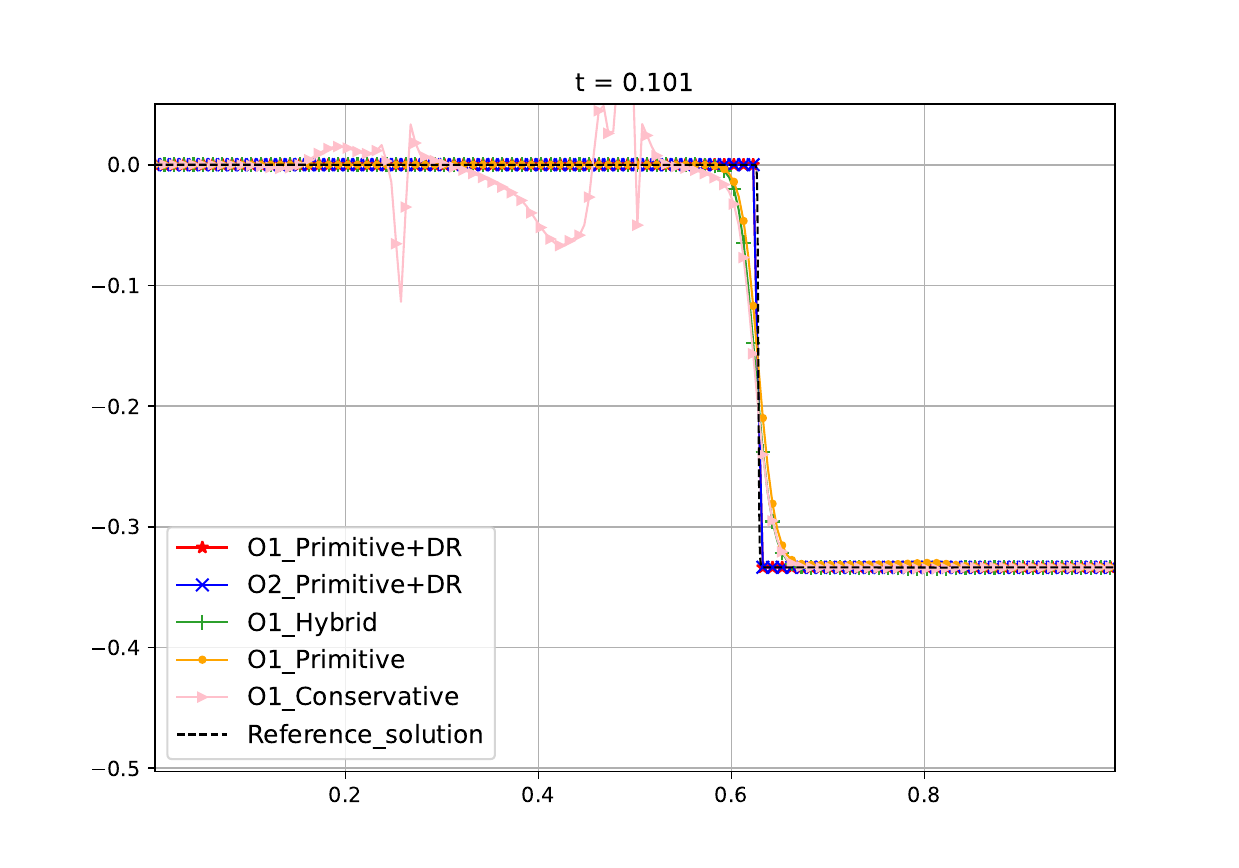}
    }
    \subfloat[Variable $u$. Zoom]{
        \includegraphics[width=0.23\textwidth]{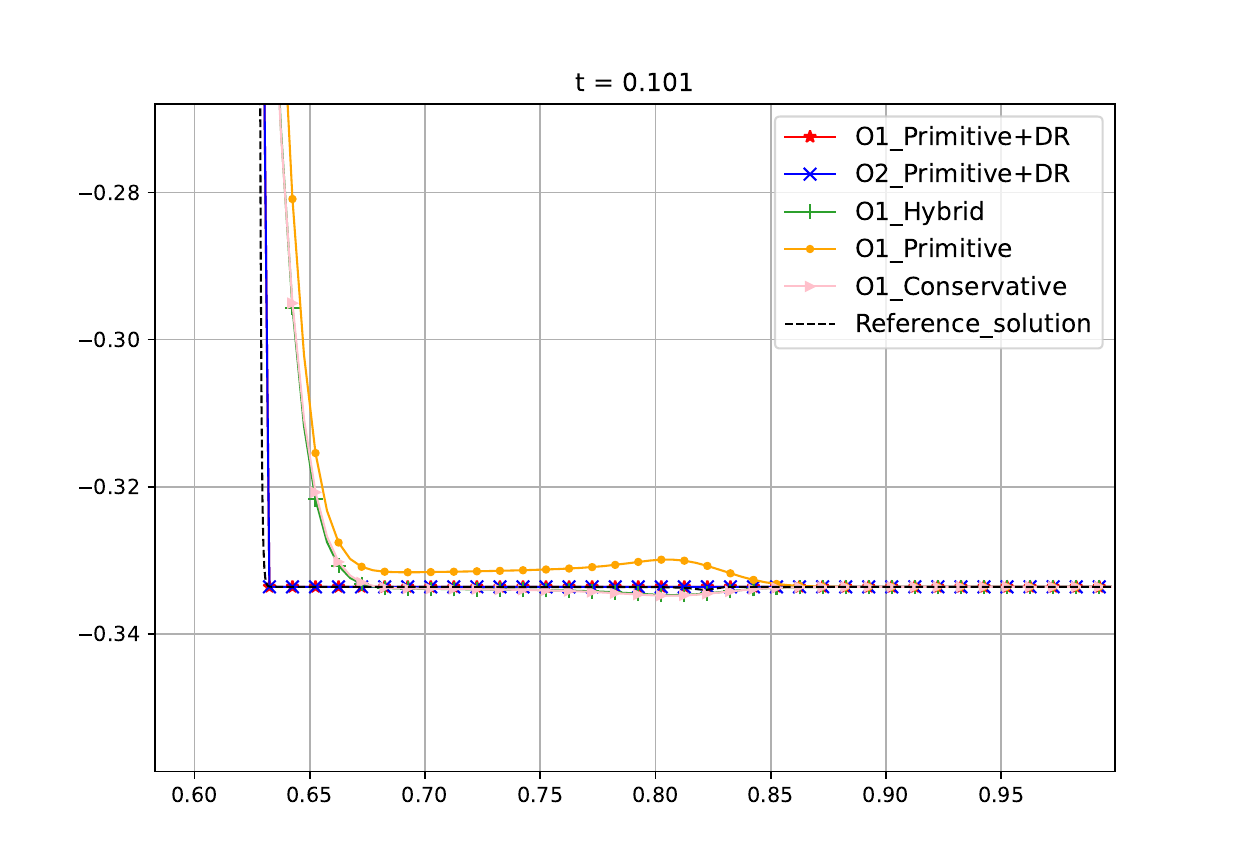}
    }
    \subfloat[Variable $u$]{
        \includegraphics[width=0.23\textwidth]{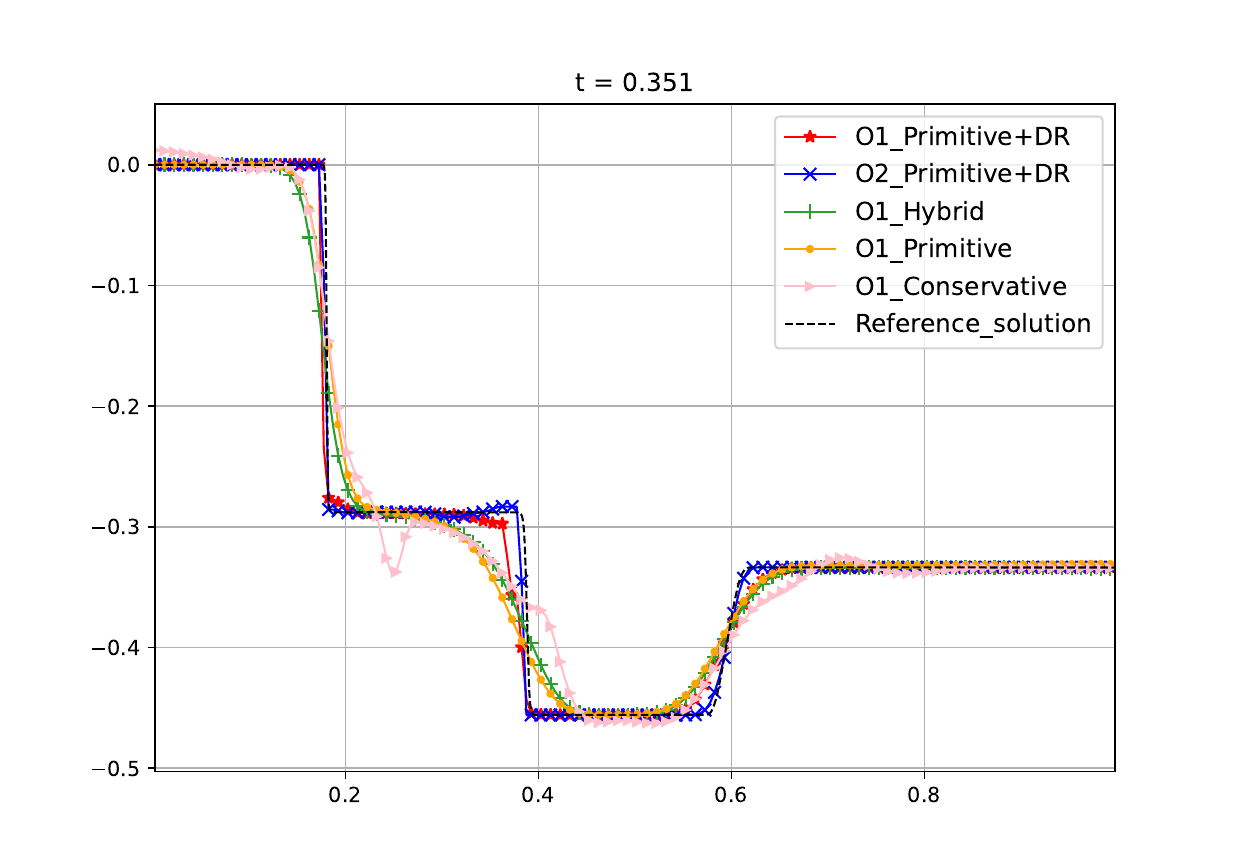}
    }
    \subfloat[Variable $u$. Zoom]{
        \includegraphics[width=0.23\textwidth]{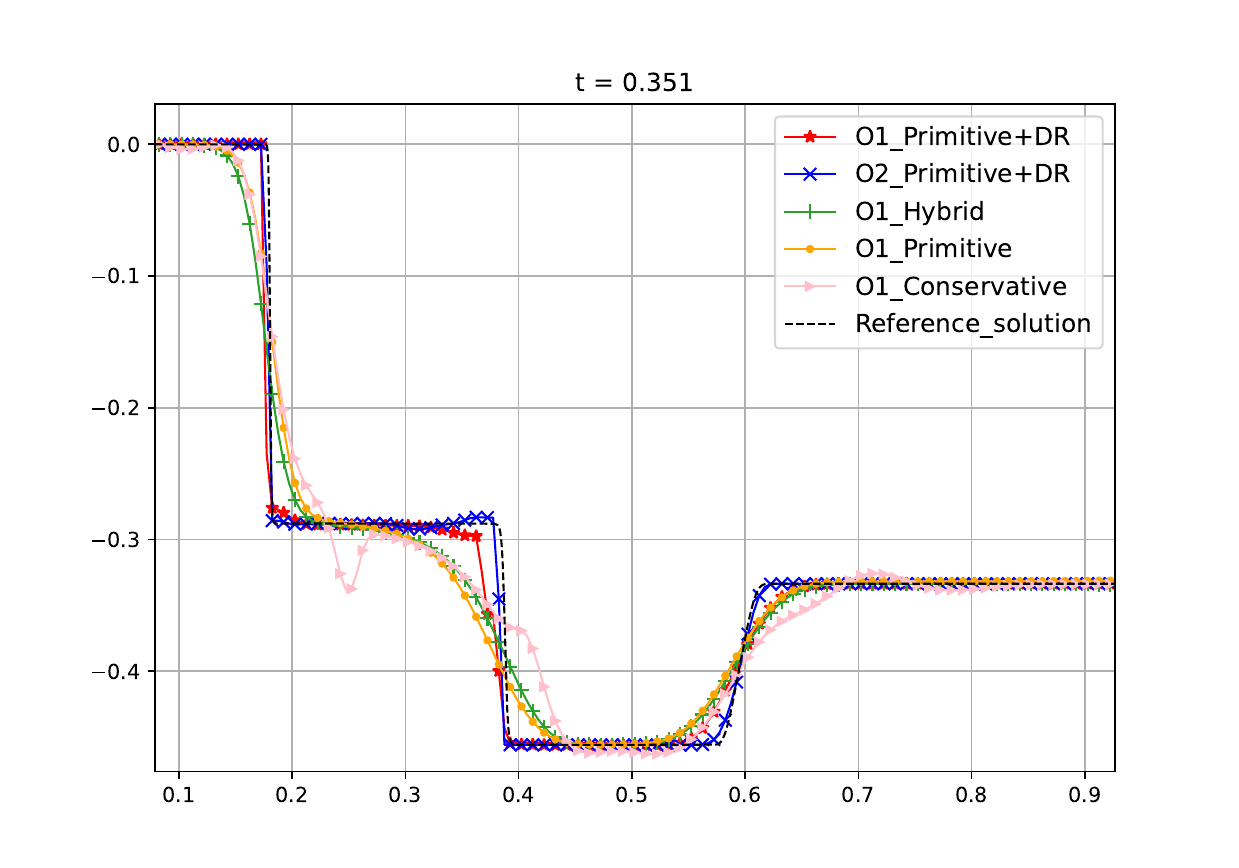}
    }

    \subfloat[Variable $p$]{
        \includegraphics[width=0.23\textwidth]{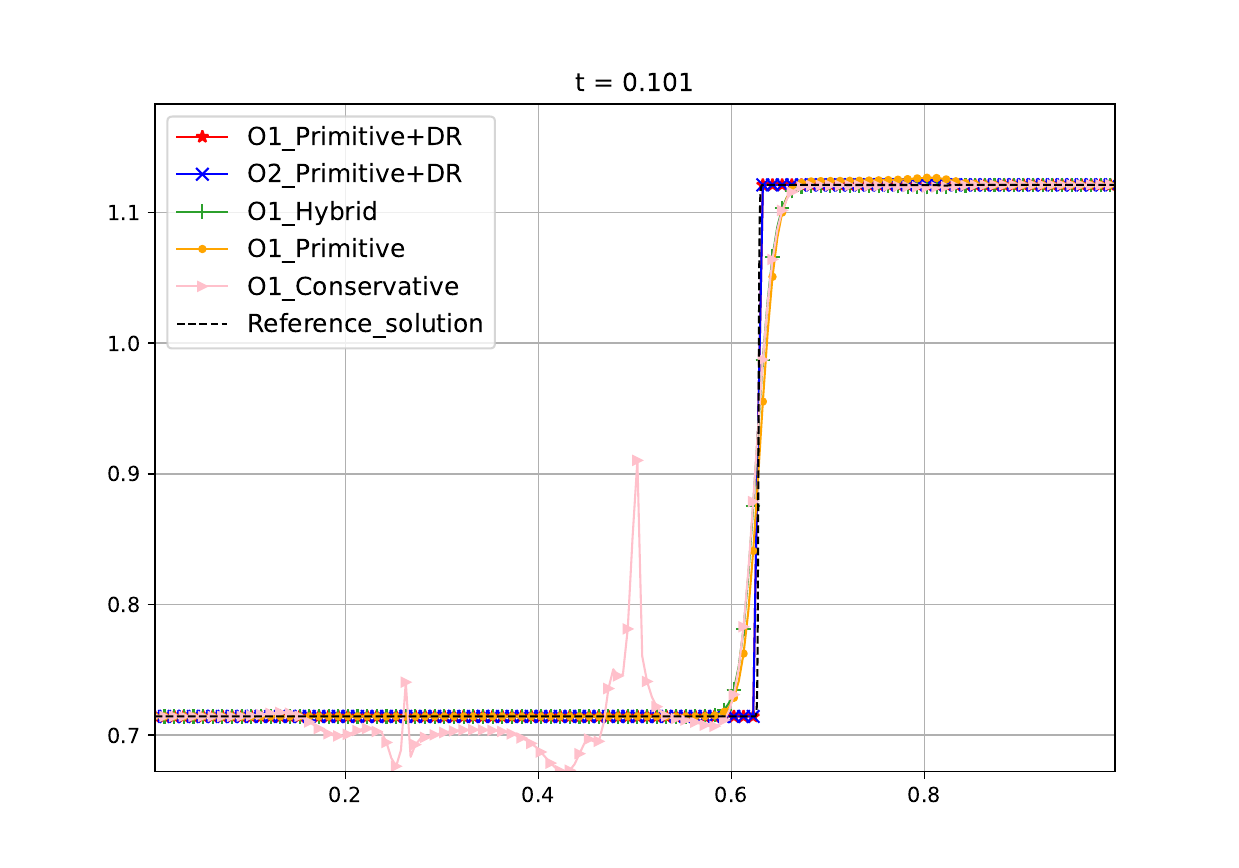}
    }
    \subfloat[Variable $p$. Zoom]{
        \includegraphics[width=0.23\textwidth]{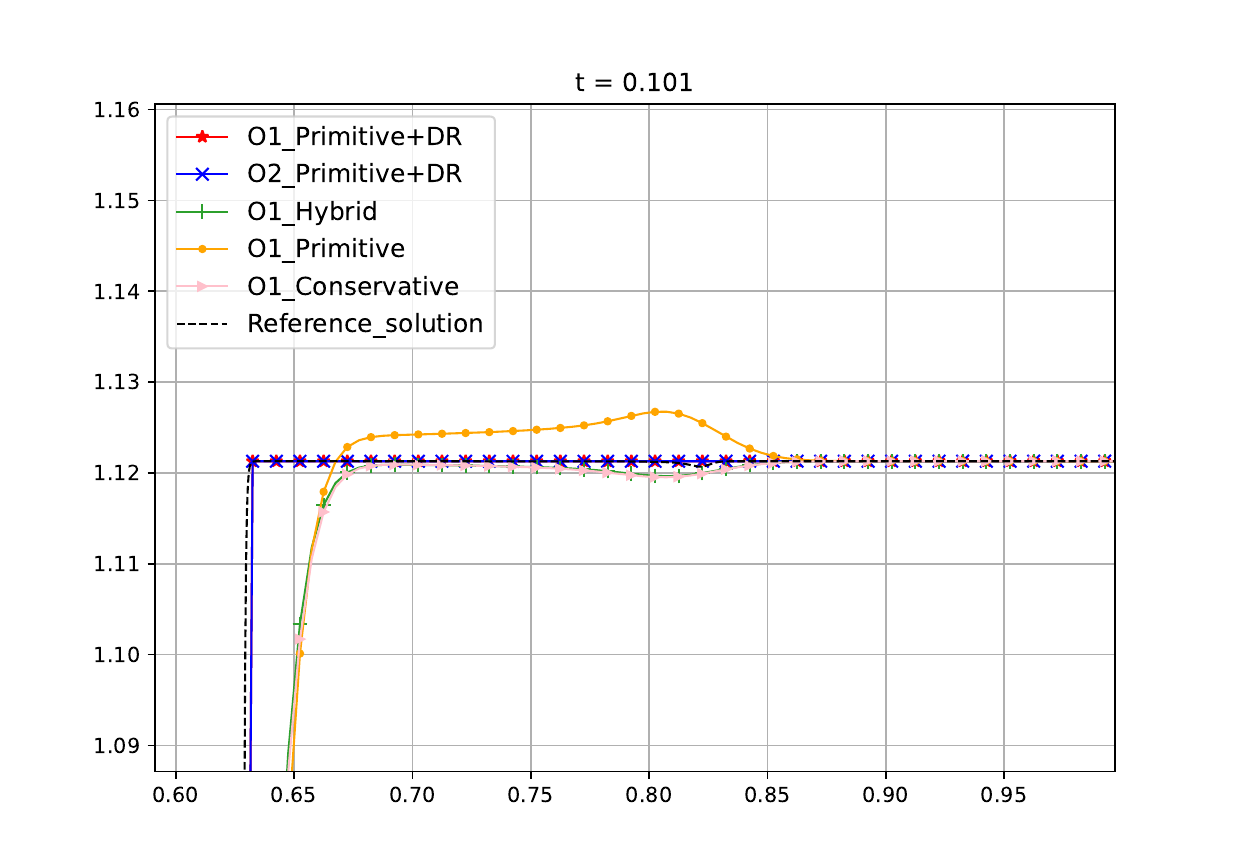}
    }
    \subfloat[Variable $p$]{
        \includegraphics[width=0.23\textwidth]{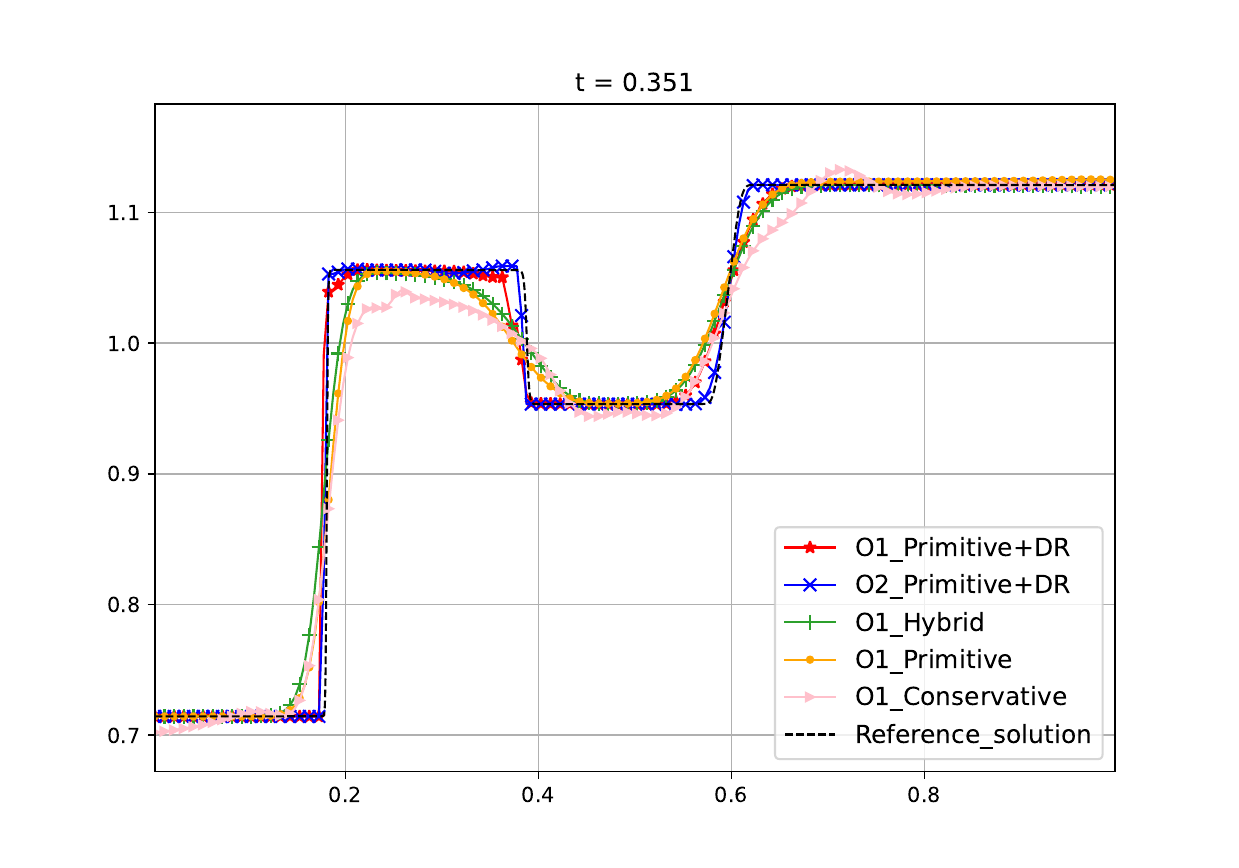}
    }
    \subfloat[Variable $p$. Zoom]{
        \includegraphics[width=0.23\textwidth]{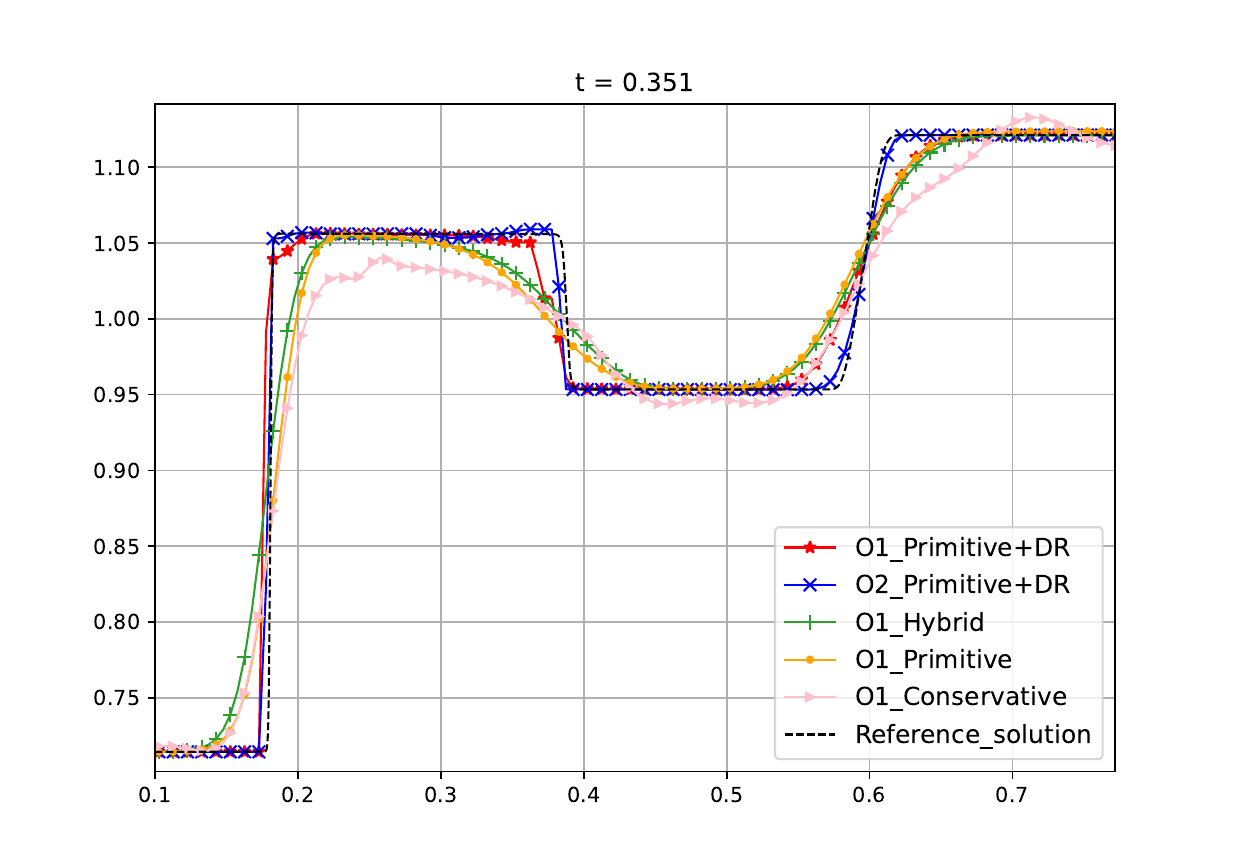}
    }
    \\
    \subfloat[Variable $c$]{
        \includegraphics[width=0.30\textwidth]{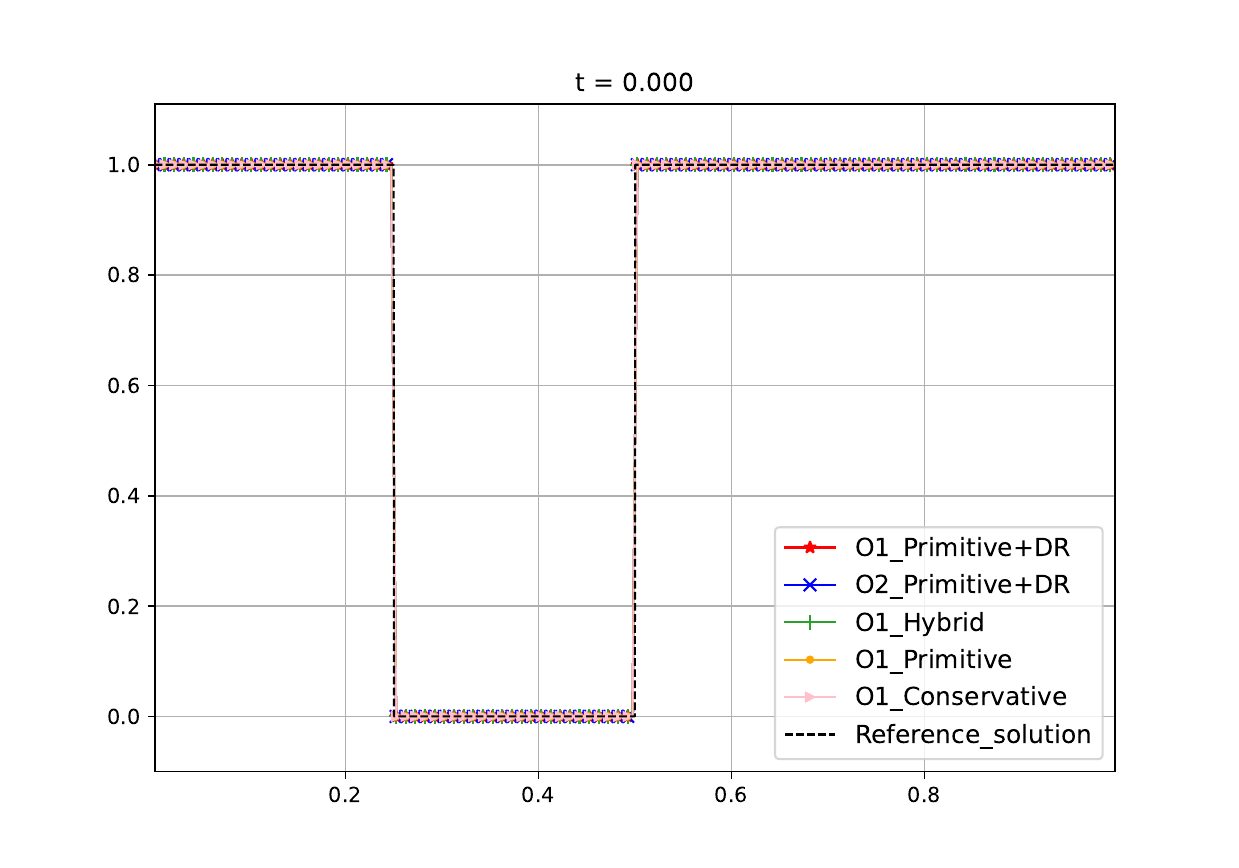}
    }
    \subfloat[Variable $c$]{
        \includegraphics[width=0.30\textwidth]{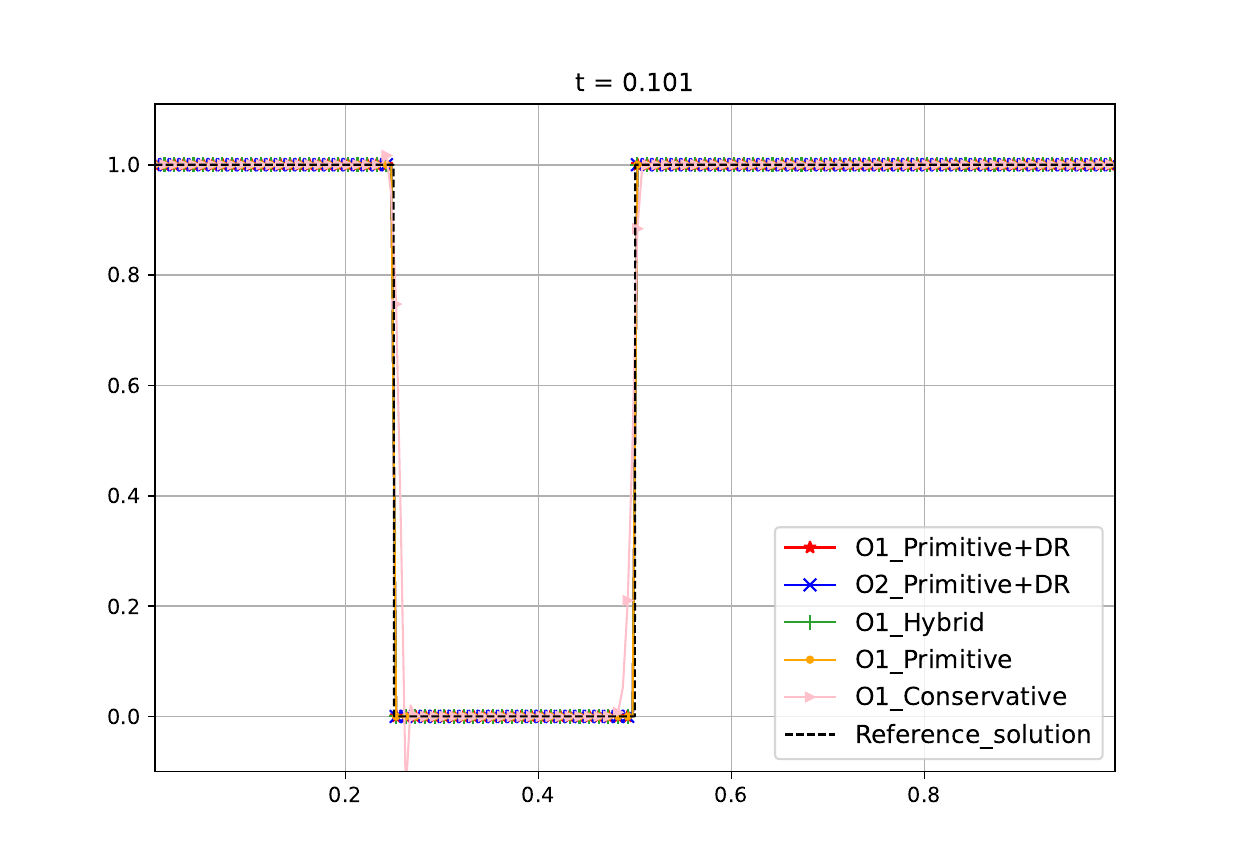}
    }
    \subfloat[Variable $c$]{
        \includegraphics[width=0.30\textwidth]{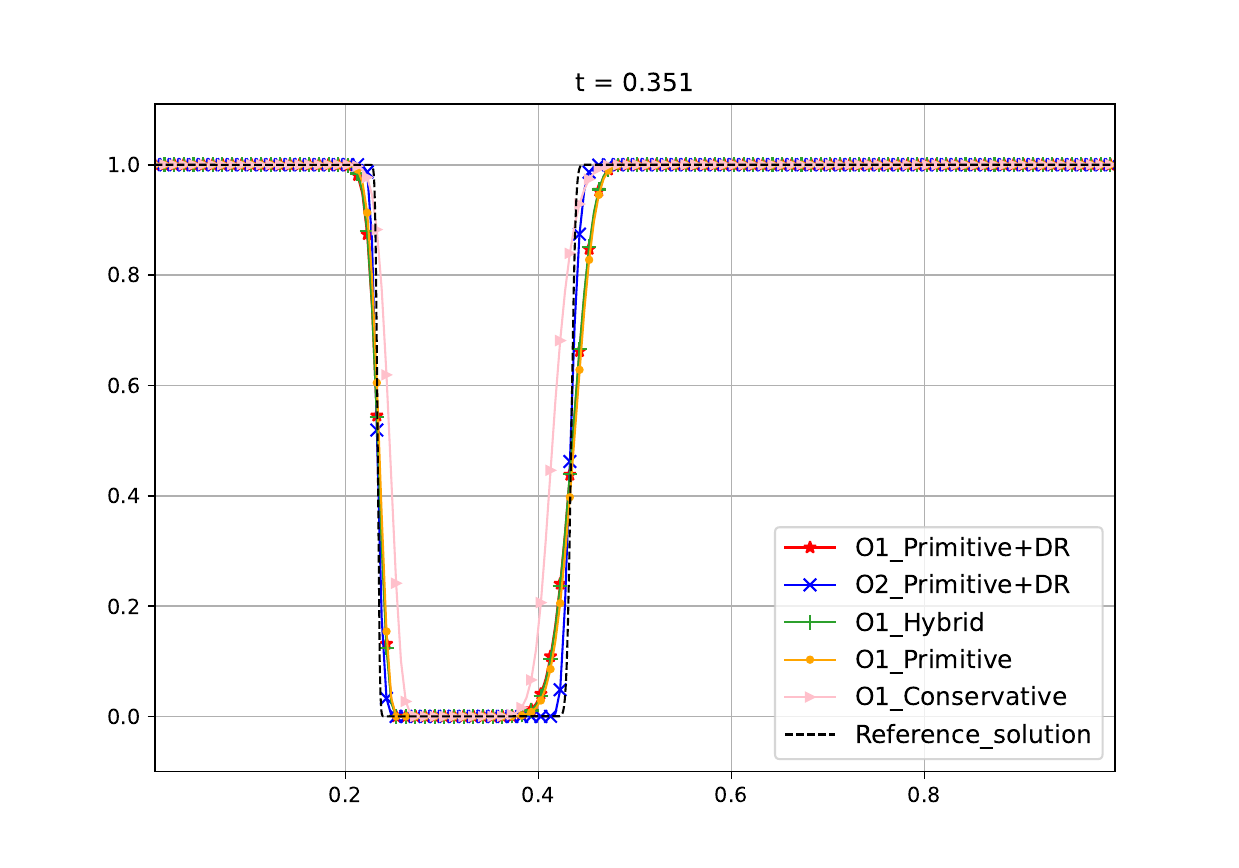}
    }

    \caption{Numerical results for Test: Shock helium-bubble interaction}
    \label{Test_Shock_Helium_bubble_interaction}
\end{figure}

 \paragraph{Gaburro et. al. shock helium-bubble interaction}
  This test is taken from \cite{GABURRO2024117311}. The initial condition is given by
 $$
 \begin{aligned}
 \bw_{0}(x) & = (\rho_0(x), u_0(x), v_0(x), p_0(x), c_0(x))^{T}  \\ & = \begin{cases}
     (1, 0, 0, 1/1.4, 1)^{T} & \text{if $x<0.25$,}\\
     (0.181, 0, 0, 1/1.4, 0)^{T} & \text{if $0.25\leq x<0.5$,}\\
     (1, 0, 0, 1/1.4, 1)^{T} & \text{if $0.5\leq x\leq 0.75$,}\\
     (1.3764, -0.3336, 0, 1.5698/1.4, 1)^{T} & \text{if $x>0.75$.}
 \end{cases}
 \end{aligned}
 $$
This test corresponds to the interaction between a shock and a helium bubble too. Figure \ref{Test_Shock_Helium_bubble_interaction} show the numerical results. The same conclusions apply as in the previous test.

\subsection{Euler equations with shear}

All tests are performed in the spatial domain $I = [0,1]$ using a mesh of $N=200$ points. In all test cases we take $\gamma = 1.4$ and the variable $c$ is set to 1, so that there is only one material (air) and we evaluate the ability of the schemes to handle a shear interface.

\paragraph{Stationary contact discontinuity}

The initial condition is given by
 $$
 \bw_{0}(x) = (\rho_0(x), u_0(x), v_0(x), p_0(x), c_0(x))^{T} = \begin{cases}
     (1, 0, 0, 1/1.4, 1)^{T} & \text{if $x<0.5$,}\\
     (0.1819, 0, 100, 1/1.4, 1)^{T} & \text{if $x\geq0.5$,}
 \end{cases}
$$

\begin{figure}[h]
    \centering

    \subfloat[Variable $\rho$]{
        \includegraphics[width=0.45\textwidth]{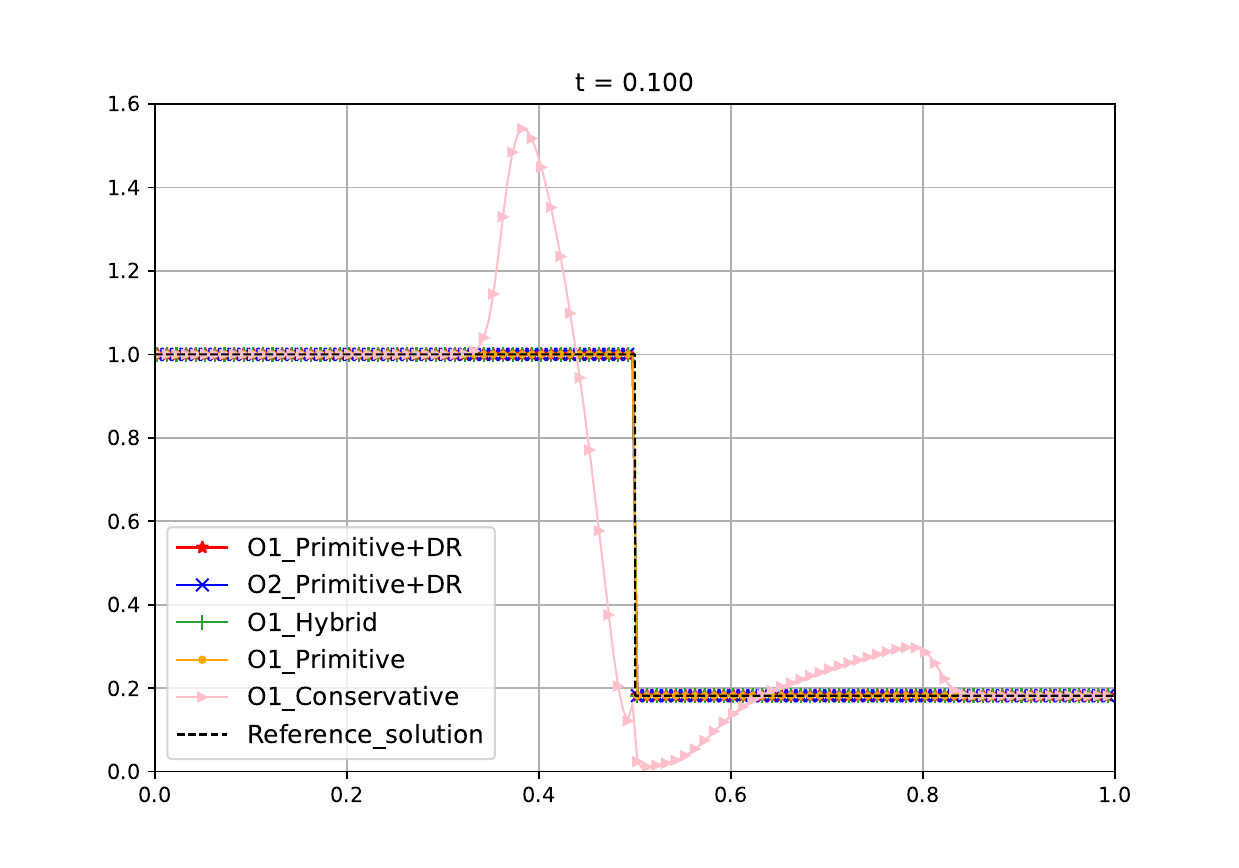}
    }
    \subfloat[Variable $u$]{
        \includegraphics[width=0.45\textwidth]{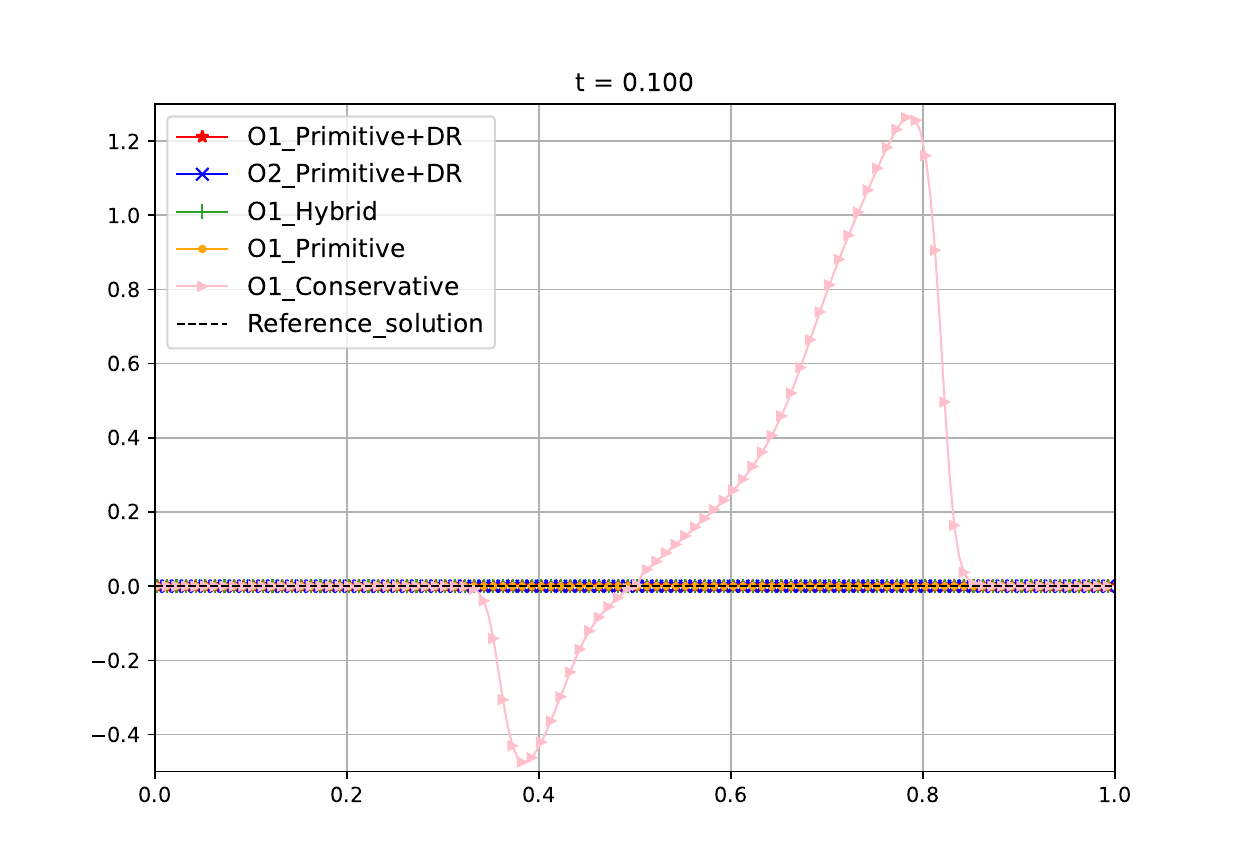}
    }
    \\

    \subfloat[Variable $p$]{
        \includegraphics[width=0.45\textwidth]{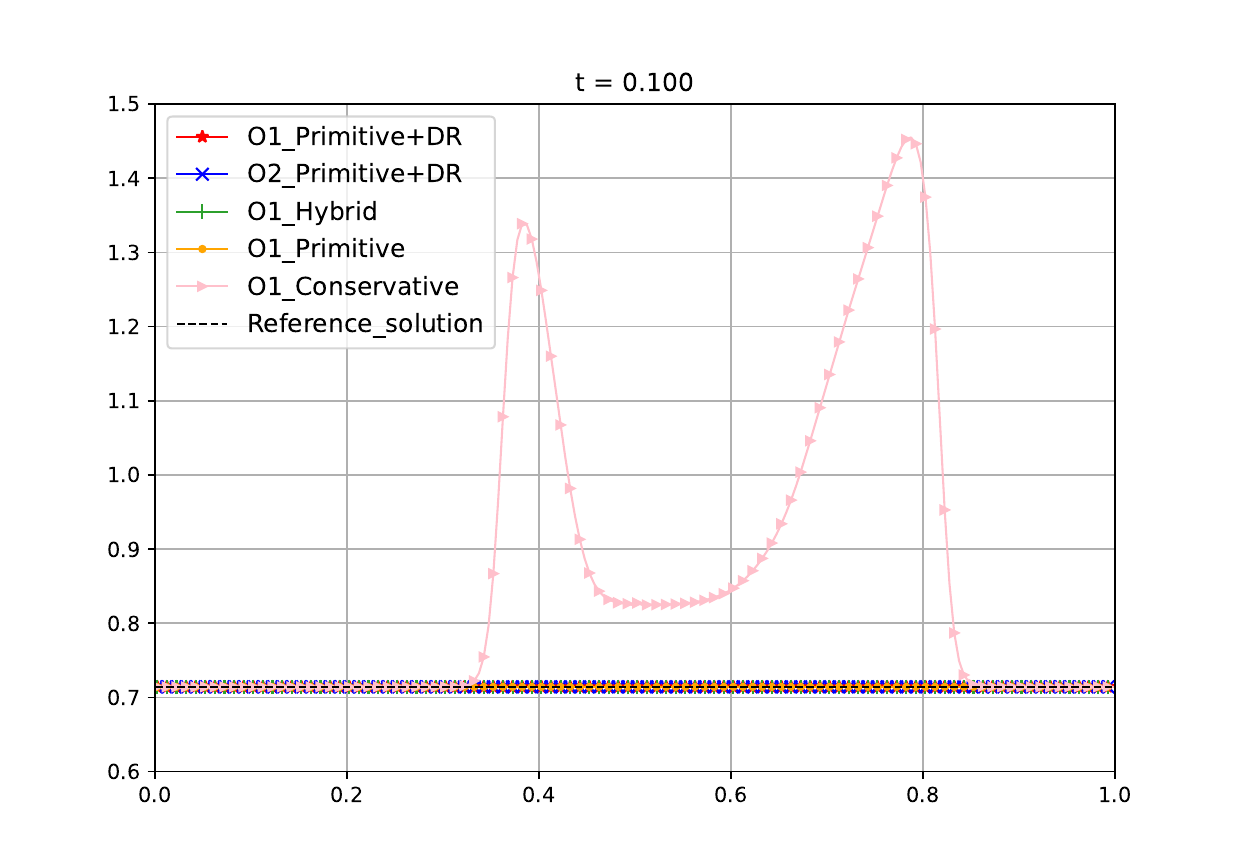}
    }
    \subfloat[Variable $v$]{
        \includegraphics[width=0.45\textwidth]{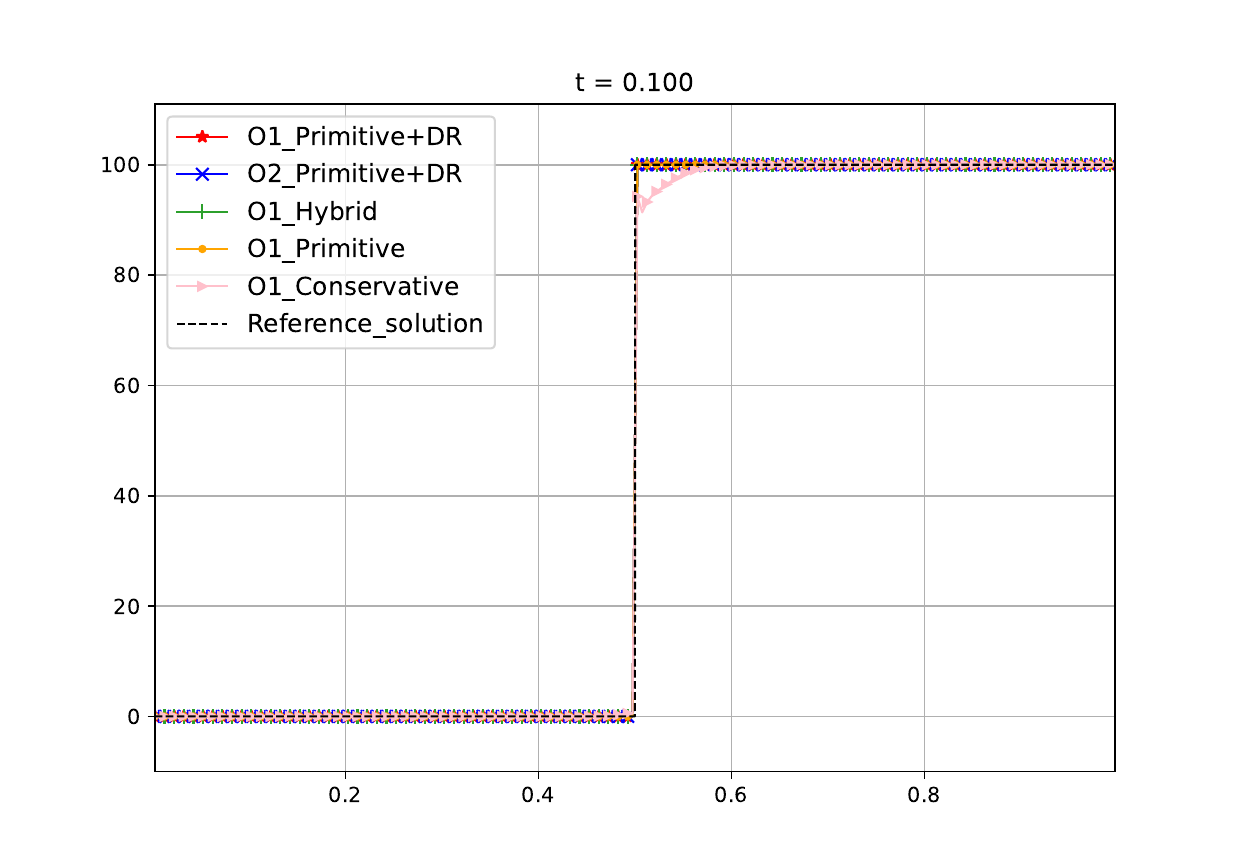}
    }

    \caption{Numerical results for Test: Stationary contact discontinuity.}
    \label{Test_Multispeed_Isolated_contact}
\end{figure}

 Figure \ref{Test_Multispeed_Isolated_contact} shows the numerical results: as in the two-material interface, we observe how all numerical schemes capture well the isolated contact discontinuity except O1\_Conservative that introduce spurious oscillations. Again, the numerical results given by O1\_Primitive and O1\_Hybrid coincide.

\paragraph{Modified Quirk and Karni small shock-shear interaction}

\begin{figure}[h]
    \centering

    \subfloat[Variable $\rho$]{
        \includegraphics[width=0.23\textwidth]{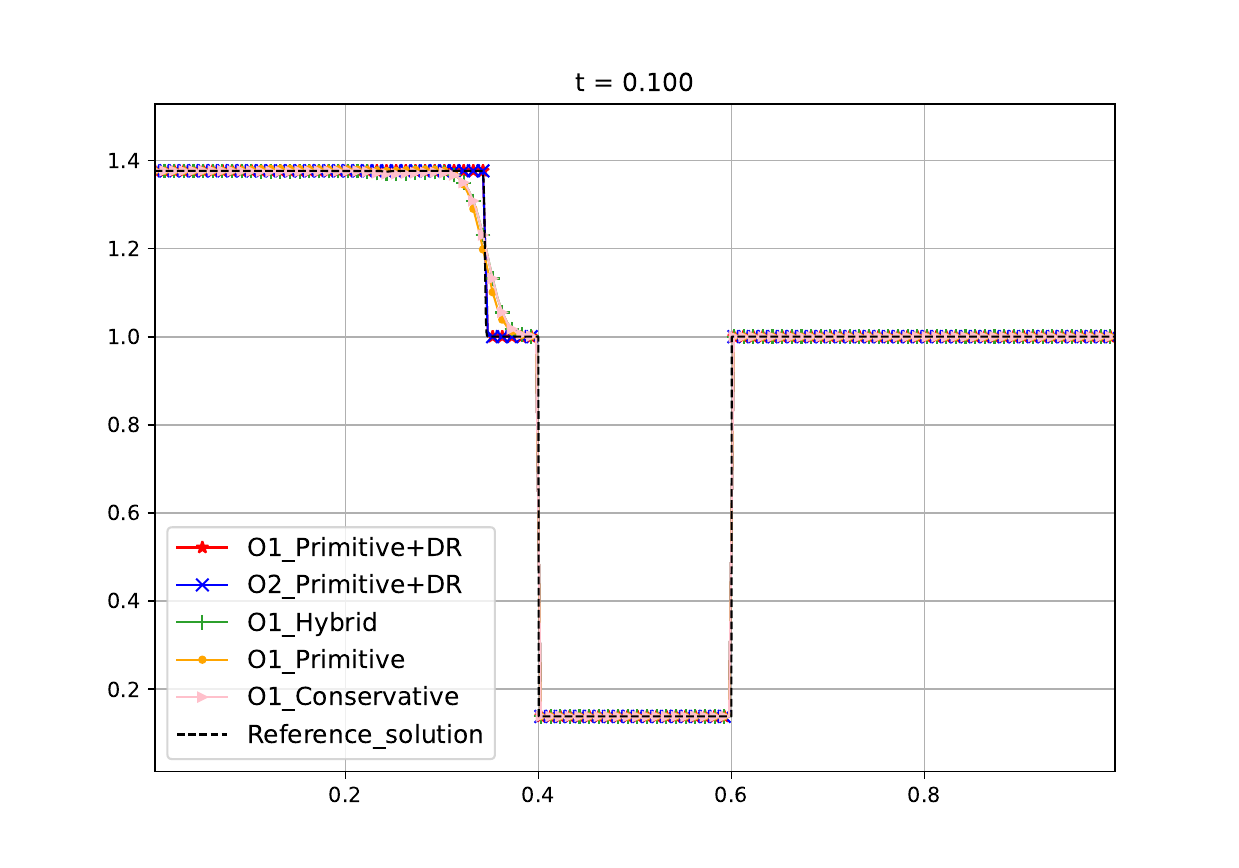}
    }
    \subfloat[Variable $\rho$. Zoom]{
        \includegraphics[width=0.23\textwidth]{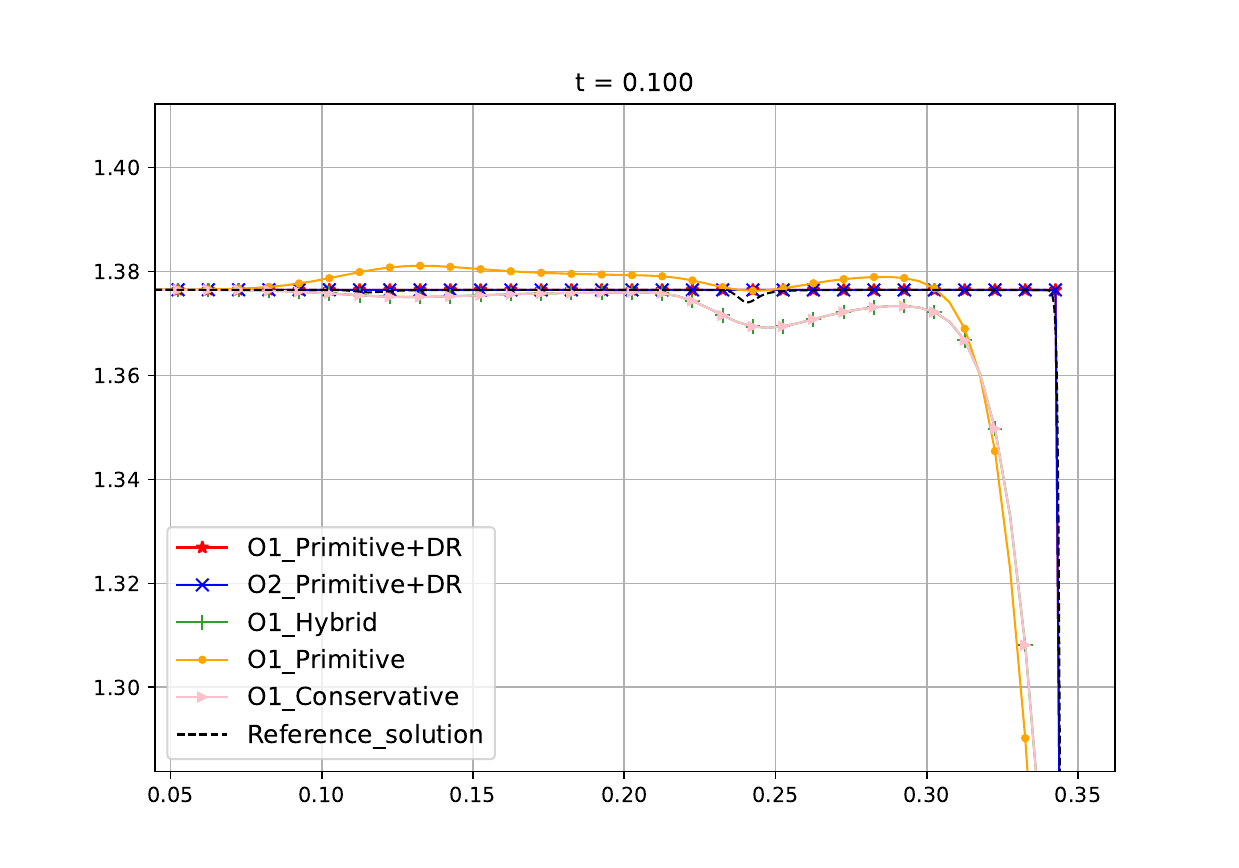}
    }
    \subfloat[Variable $\rho$]{
        \includegraphics[width=0.23\textwidth]{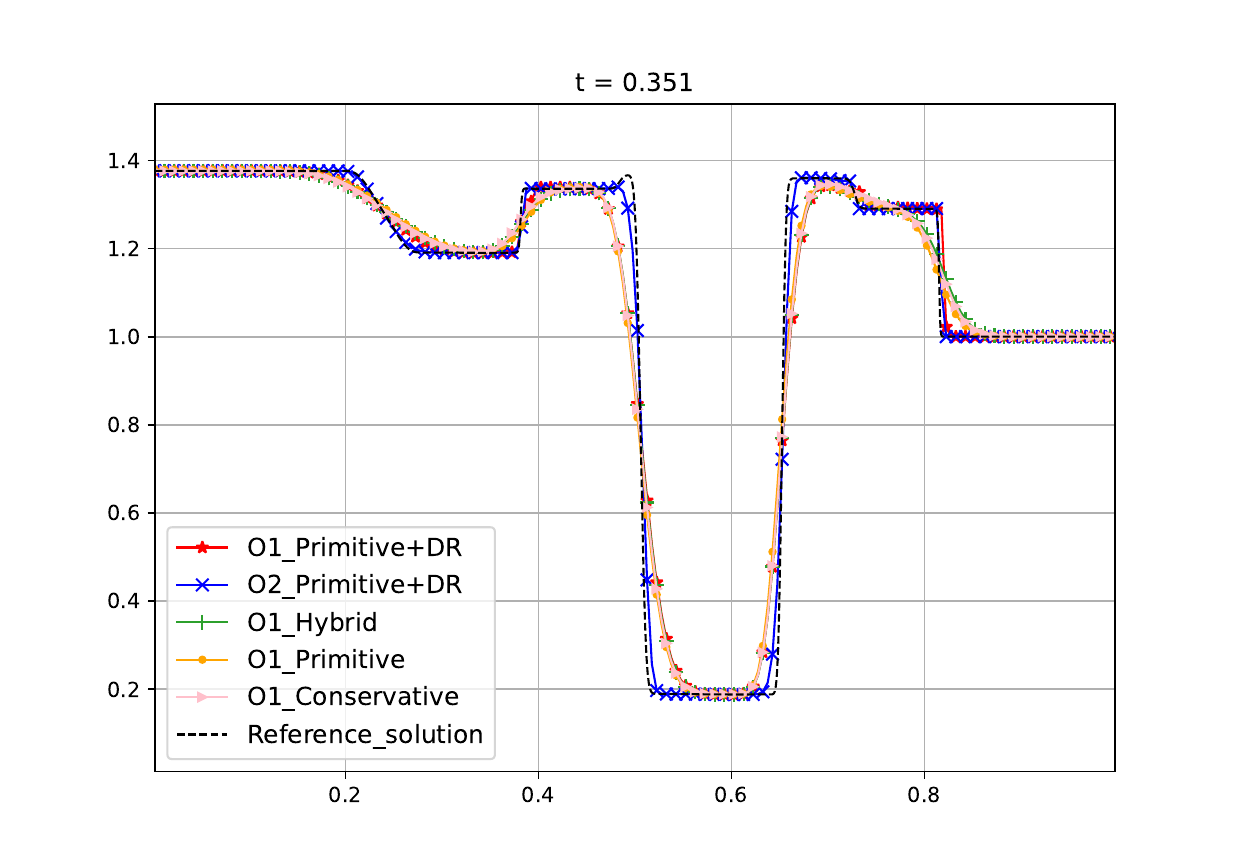}
    }
    \subfloat[Variable $\rho$. Zoom]{
        \includegraphics[width=0.23\textwidth]{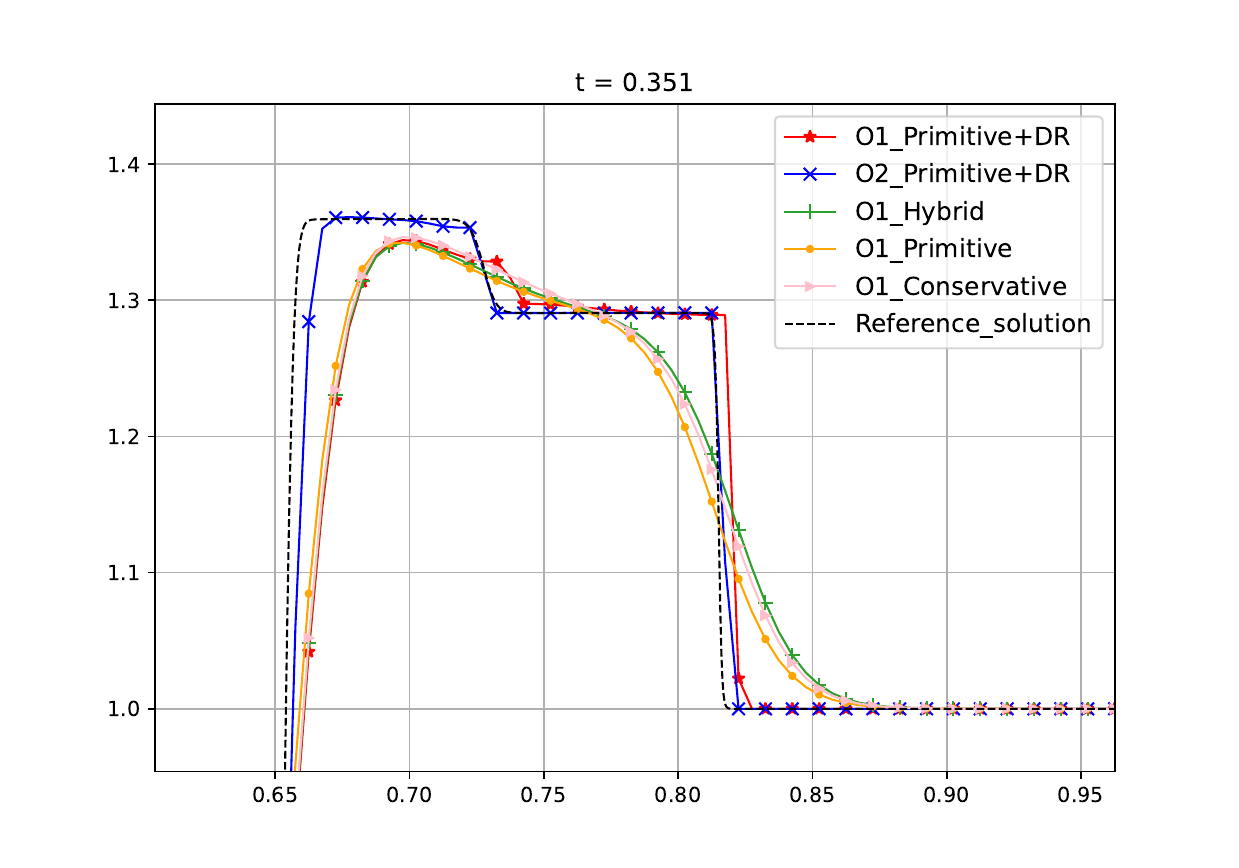}
    }

    \subfloat[Variable $u$]{
        \includegraphics[width=0.23\textwidth]{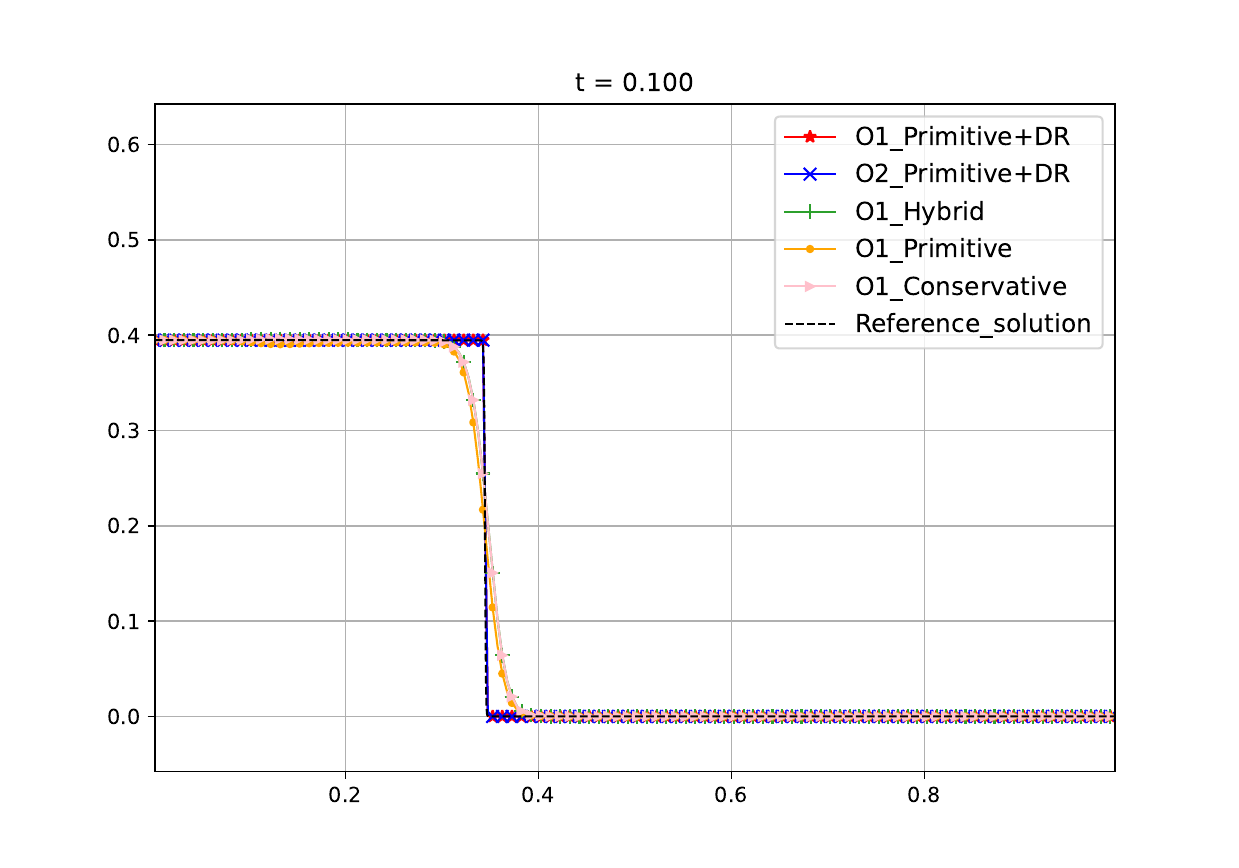}
    }
    \subfloat[Variable $u$. Zoom]{
        \includegraphics[width=0.23\textwidth]{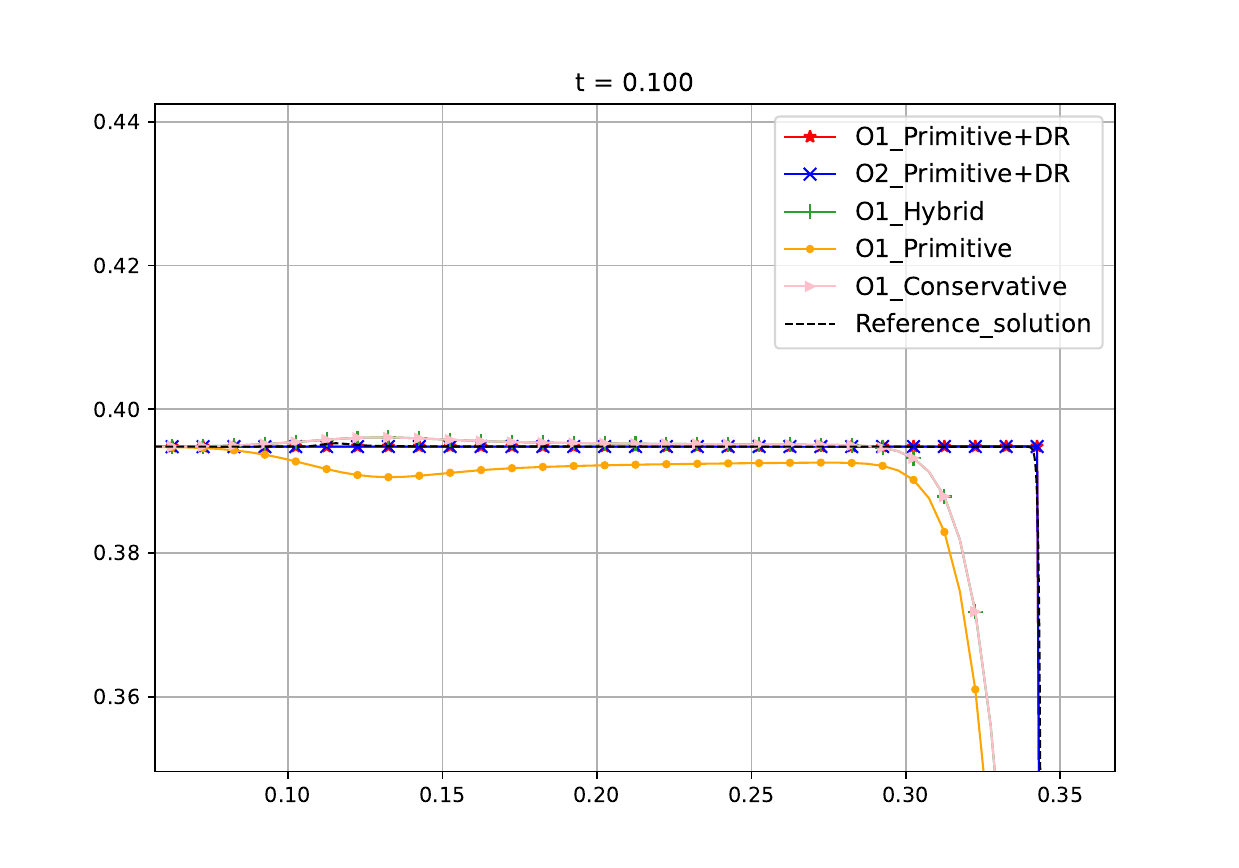}
    }
    \subfloat[Variable $u$]{
        \includegraphics[width=0.23\textwidth]{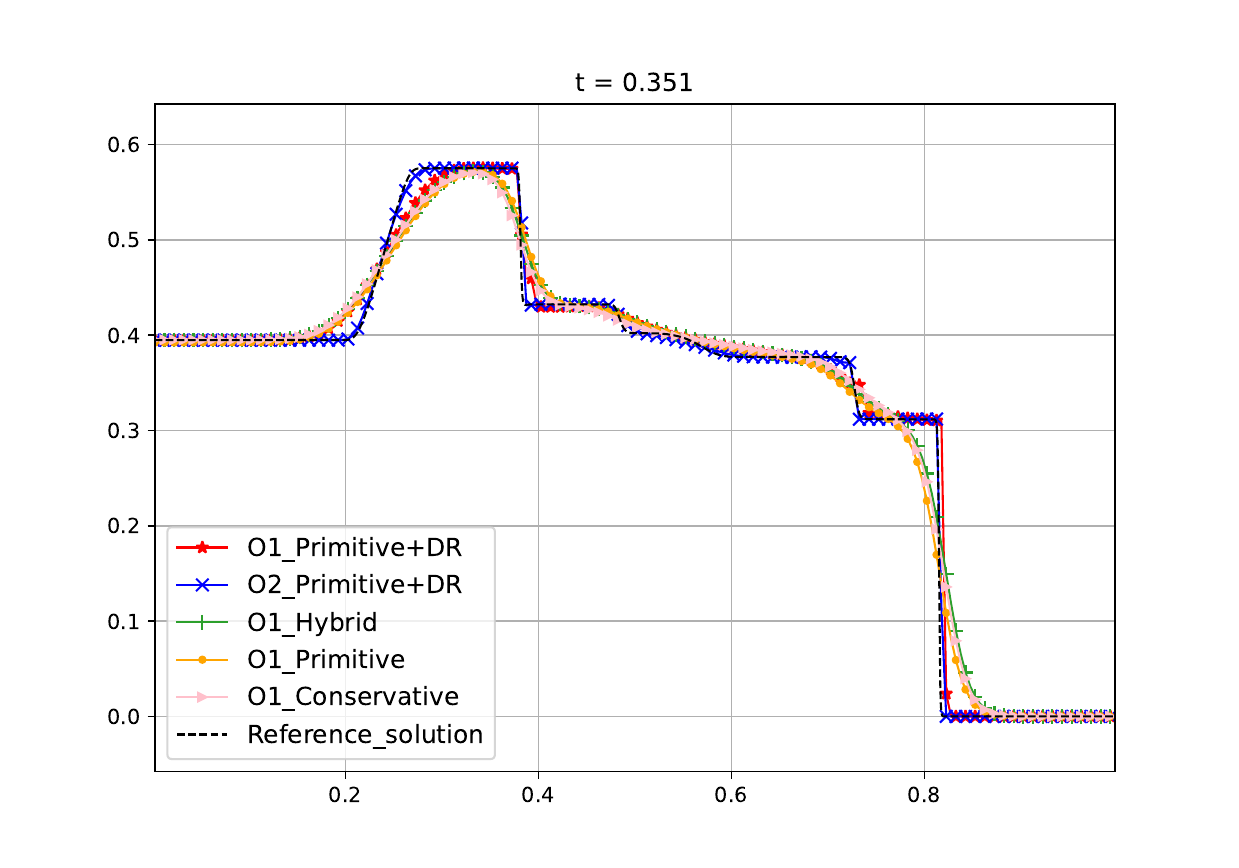}
    }
    \subfloat[Variable $u$. Zoom]{
        \includegraphics[width=0.23\textwidth]{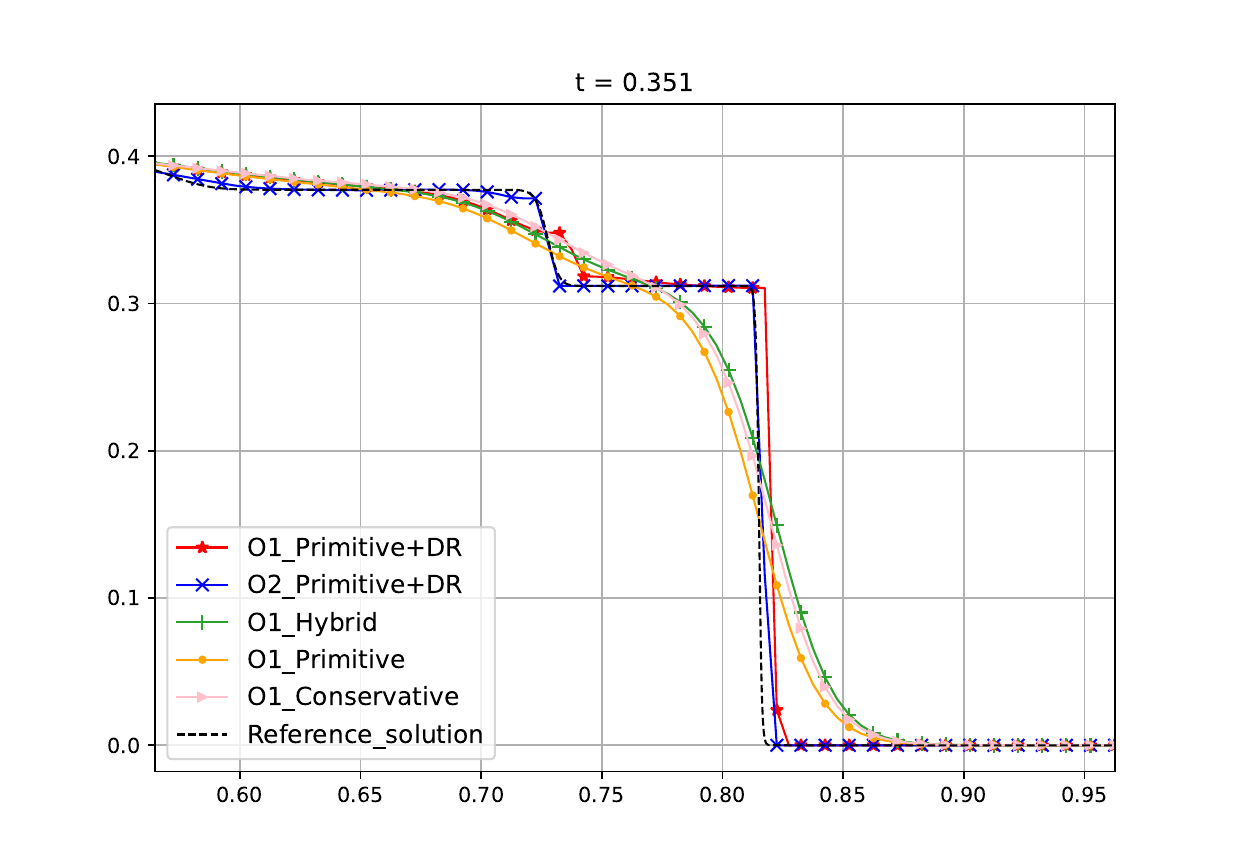}
    }

    \subfloat[Variable $p$]{
        \includegraphics[width=0.23\textwidth]{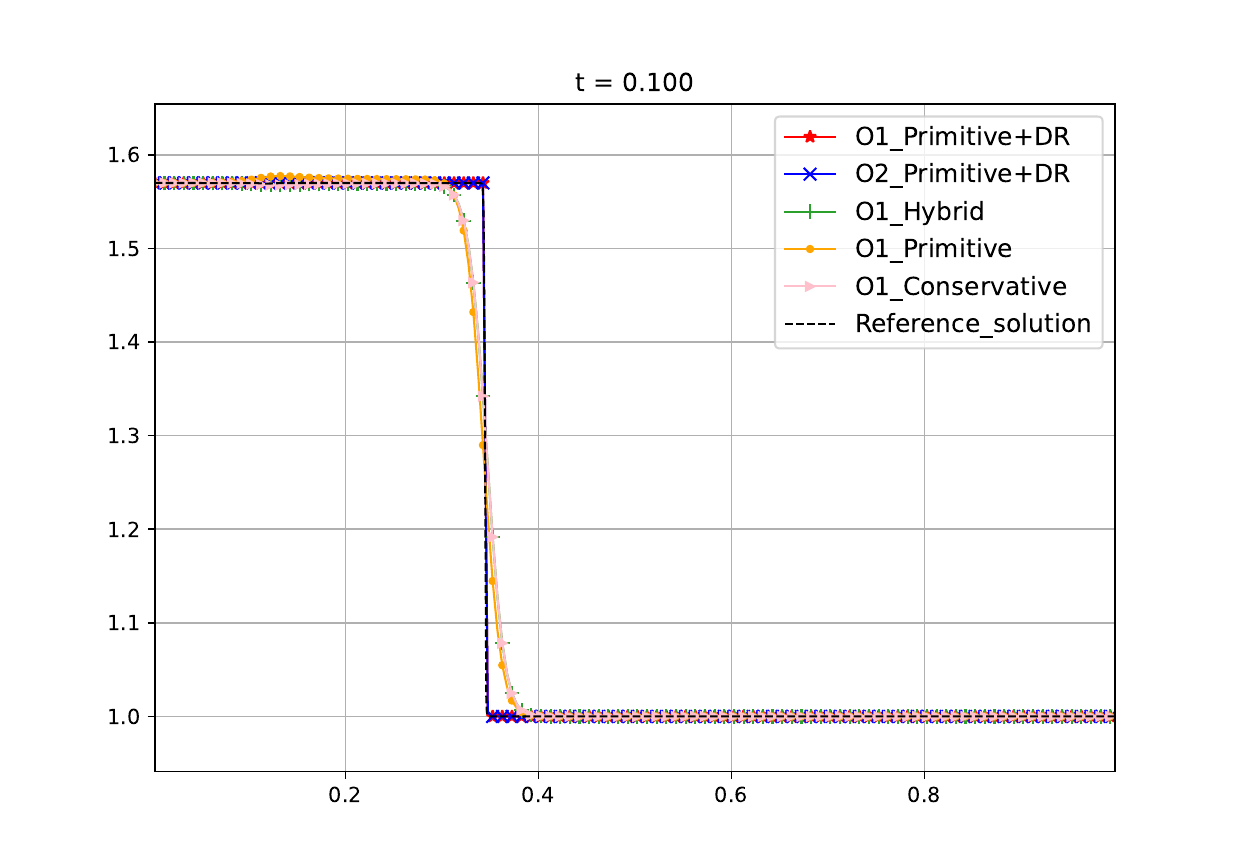}
    }
    \subfloat[Variable $p$. Zoom]{
        \includegraphics[width=0.23\textwidth]{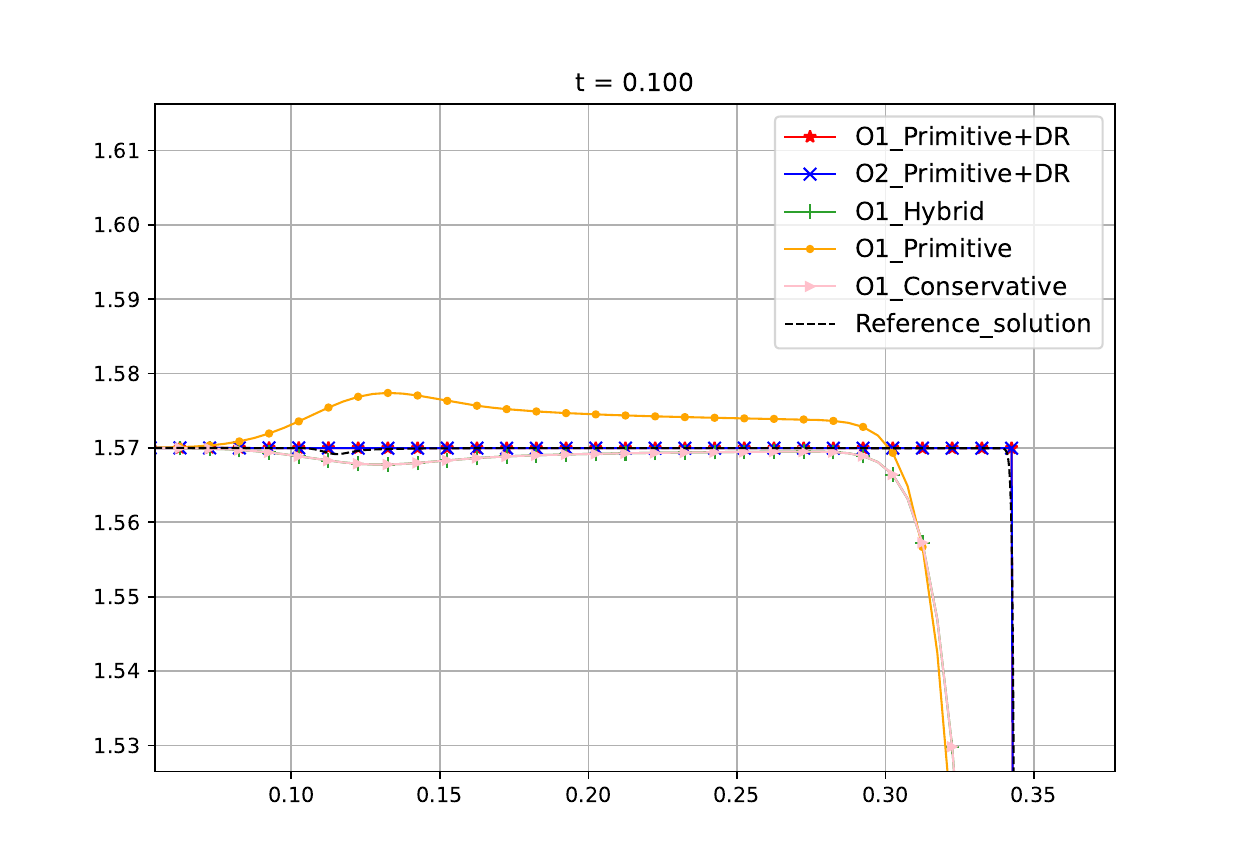}
    }
    \subfloat[Variable $p$]{
        \includegraphics[width=0.23\textwidth]{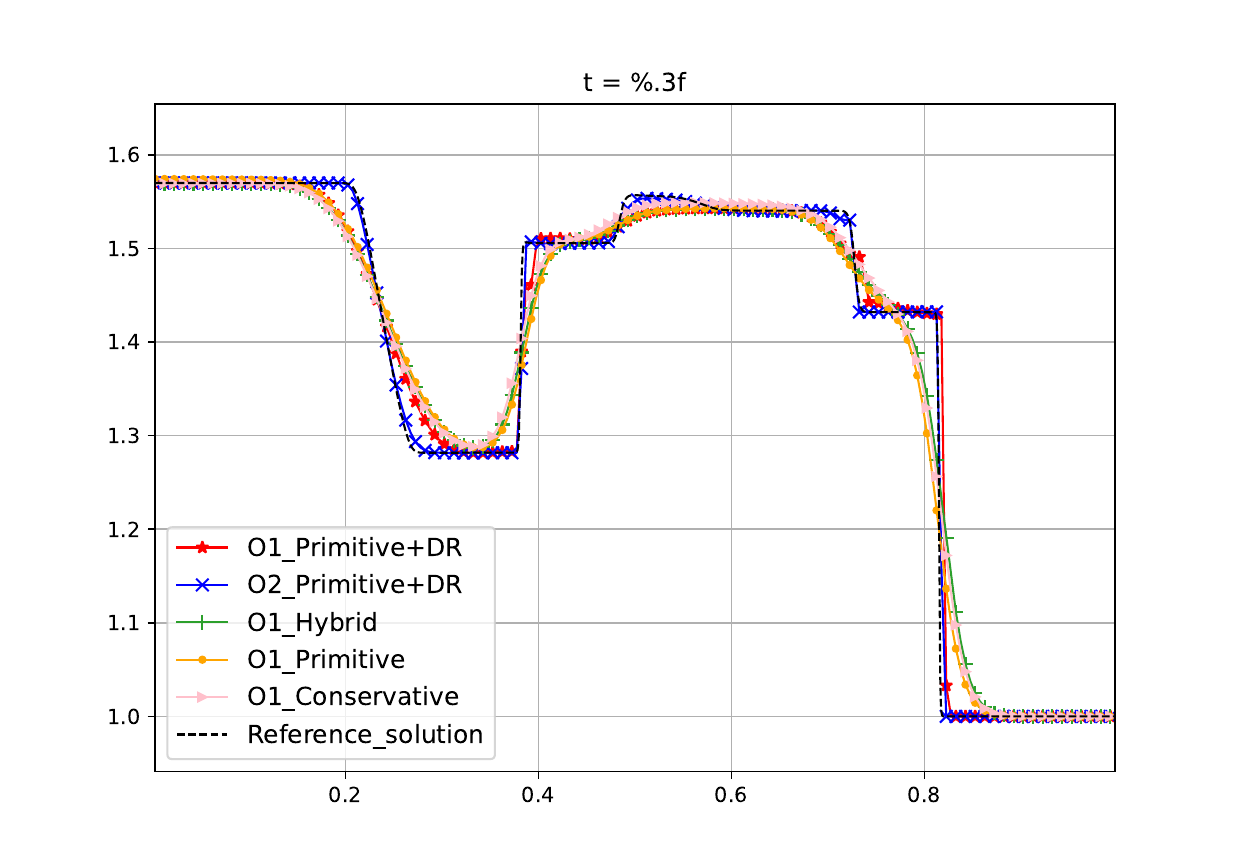}
    }
    \subfloat[Variable $p$. Zoom]{
        \includegraphics[width=0.23\textwidth]{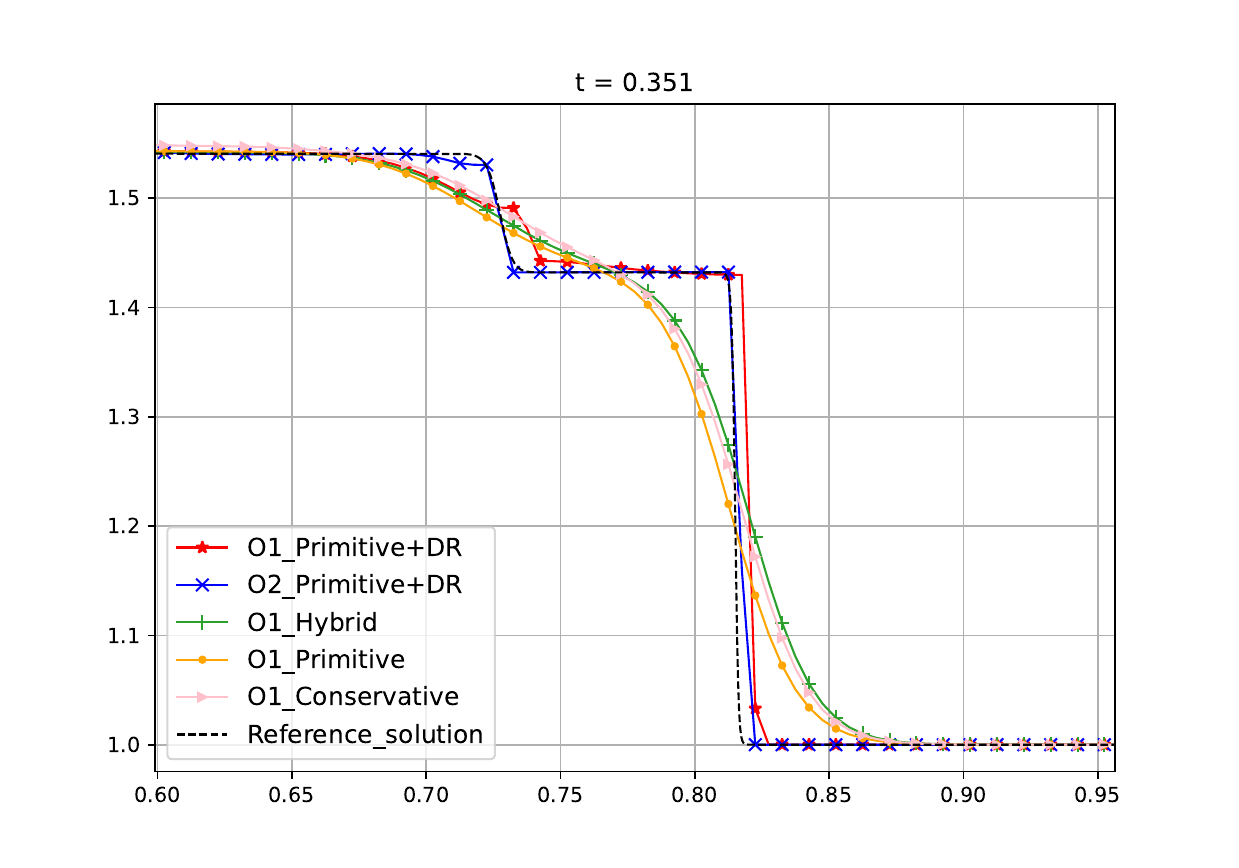}
    }
    \\
    \subfloat[Variable $v$]{
        \includegraphics[width=0.30\textwidth]{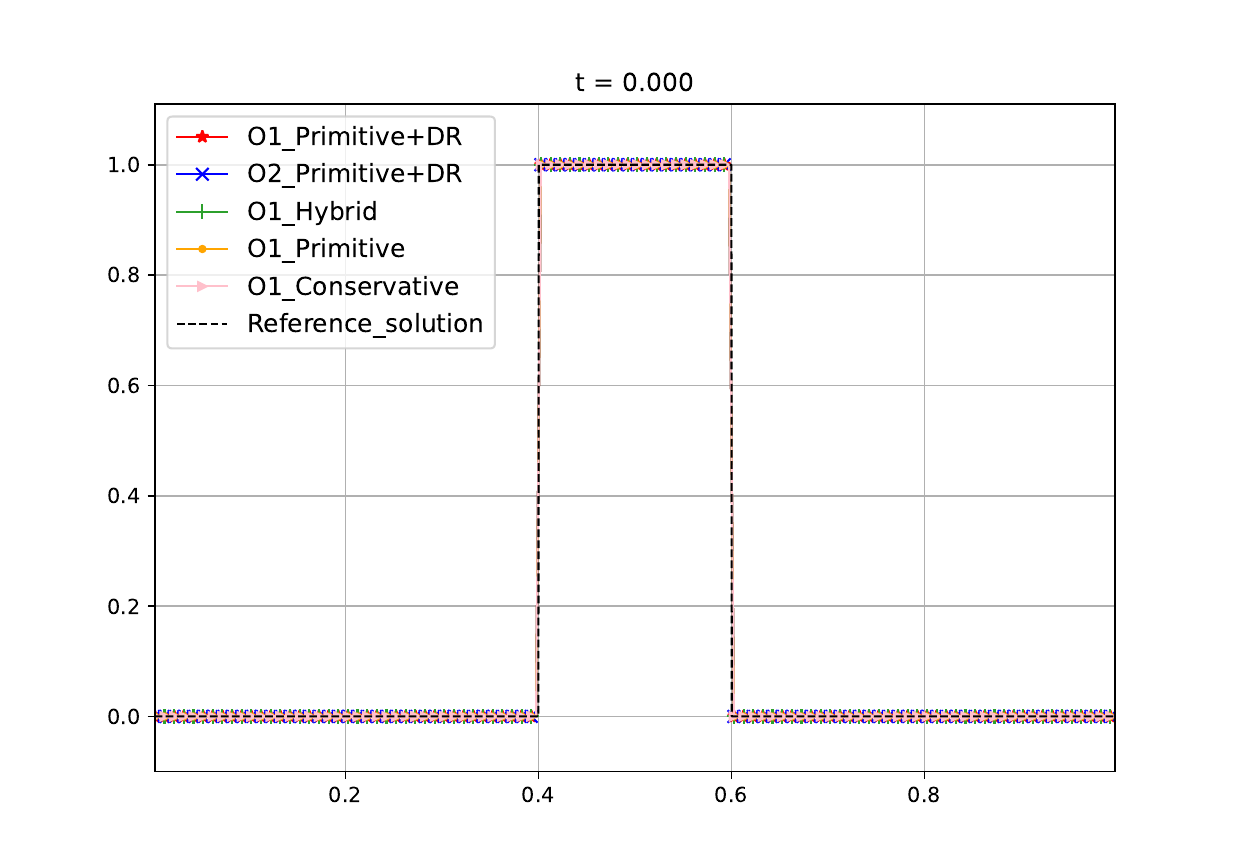}
    }
    \subfloat[Variable $v$]{
        \includegraphics[width=0.30\textwidth]{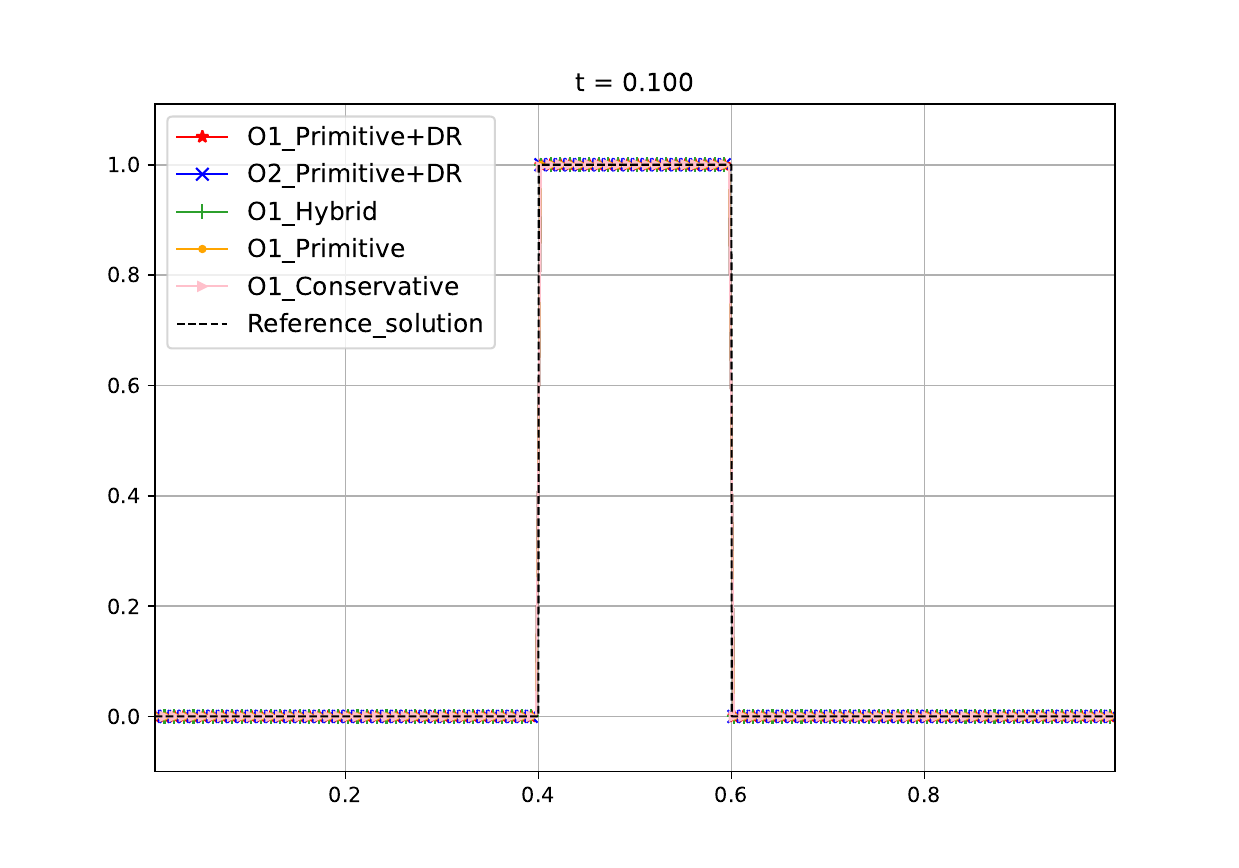}
    }
    \subfloat[Variable $v$]{
        \includegraphics[width=0.30\textwidth]{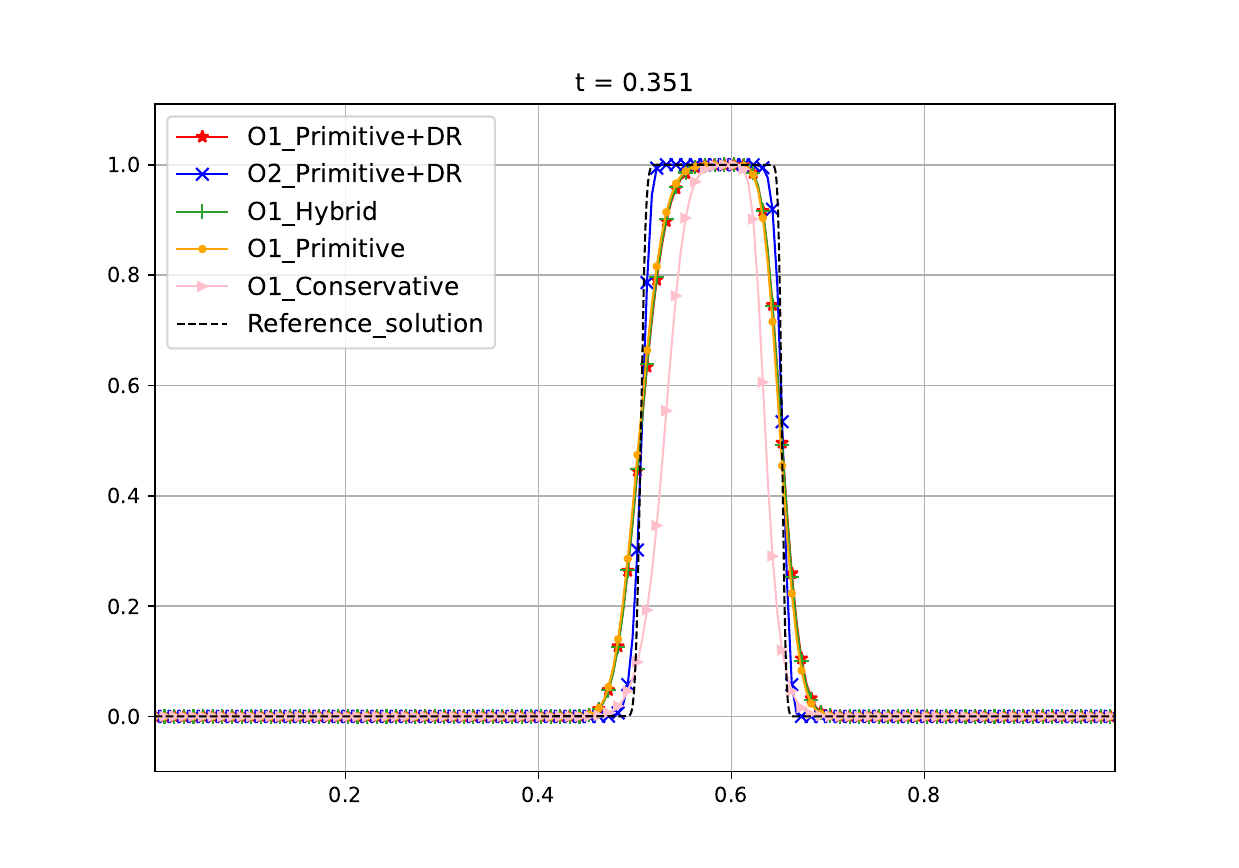}
    }

    \caption{Numerical results for Test: modified Quirk and Karni small shock-shear interaction.}
    \label{Test_Multispeed_Karni}
\end{figure}

This test is a modification of the one taken from Quirk and Karni in \cite{Quirk1996129}. The initial condition is given by
 $$
 \begin{aligned}
 \bw_{0}(x) & = (\rho_0(x), u_0(x), v_0(x), p_0(x), c_0(x))^{T} \\ & = \begin{cases}
     (1.3765, 0.3948, 0, 1.57,1)^{T} & \text{if $x<0.2$,}\\
     (1, 0, 0, 1,1)^{T} & \text{if $0.2\leq x<0.4$,}\\
     (0.138, 0, 1, 1,1)^{T} & \text{if $0.4\leq x\leq 0.6$,}\\
     (1, 0, 0, 1,1)^{T} & \text{if $x>0.6$.}
 \end{cases}
 \end{aligned}
 $$
This test corresponds to the interaction between a shock and a small shear wave. Figure \ref{Test_Multispeed_Karni} shows the numerical results. In this case, the only one numerical method that is not able to capture the correct solution is the standard path-conservative scheme without in-cell discontinuous reconstruction. It seems that the contact discontinuity related to the shear interface is not big enough to produce spurious oscillations in the fully conservative method. Again, the schemes with the in-cell discontinuous reconstruction operator perform better close to  shock waves. The second-order scheme captures better the rarefactions as expected.

\paragraph{Modified Quirk and Karni strong  shock-shear interaction}

\begin{figure}[h]
    \centering

    \subfloat[Variable $\rho$]{
        \includegraphics[width=0.23\textwidth]{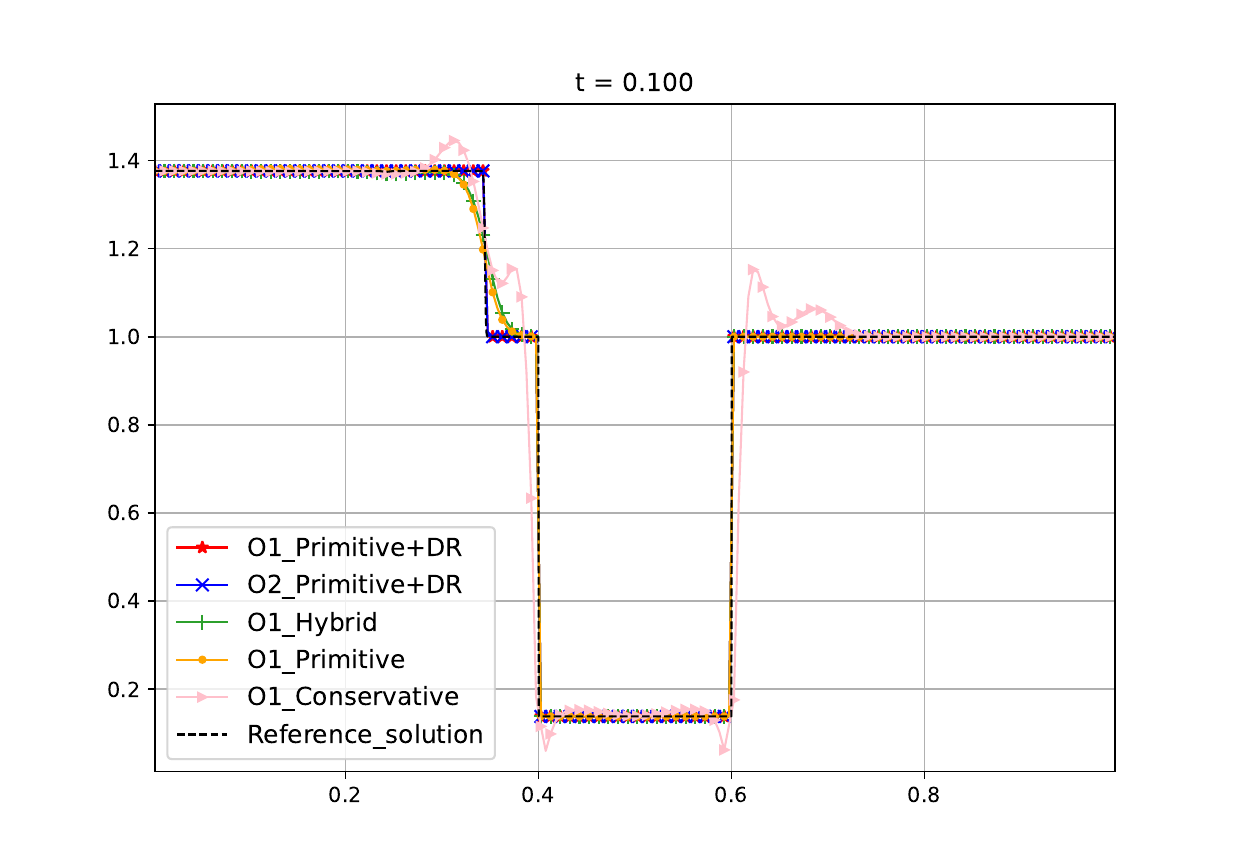}
    }
    \subfloat[Variable $\rho$. Zoom]{
        \includegraphics[width=0.23\textwidth]{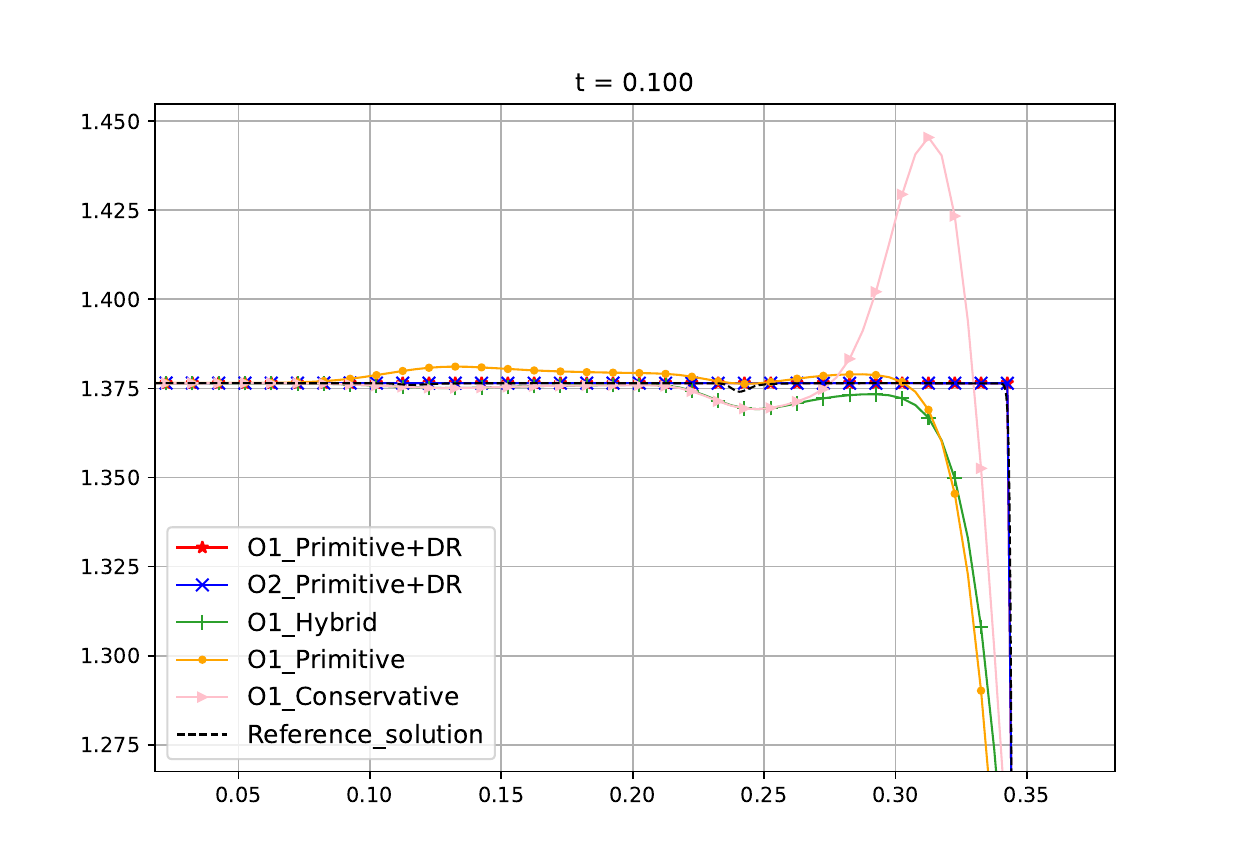}
    }
    \subfloat[Variable $\rho$]{
        \includegraphics[width=0.23\textwidth]{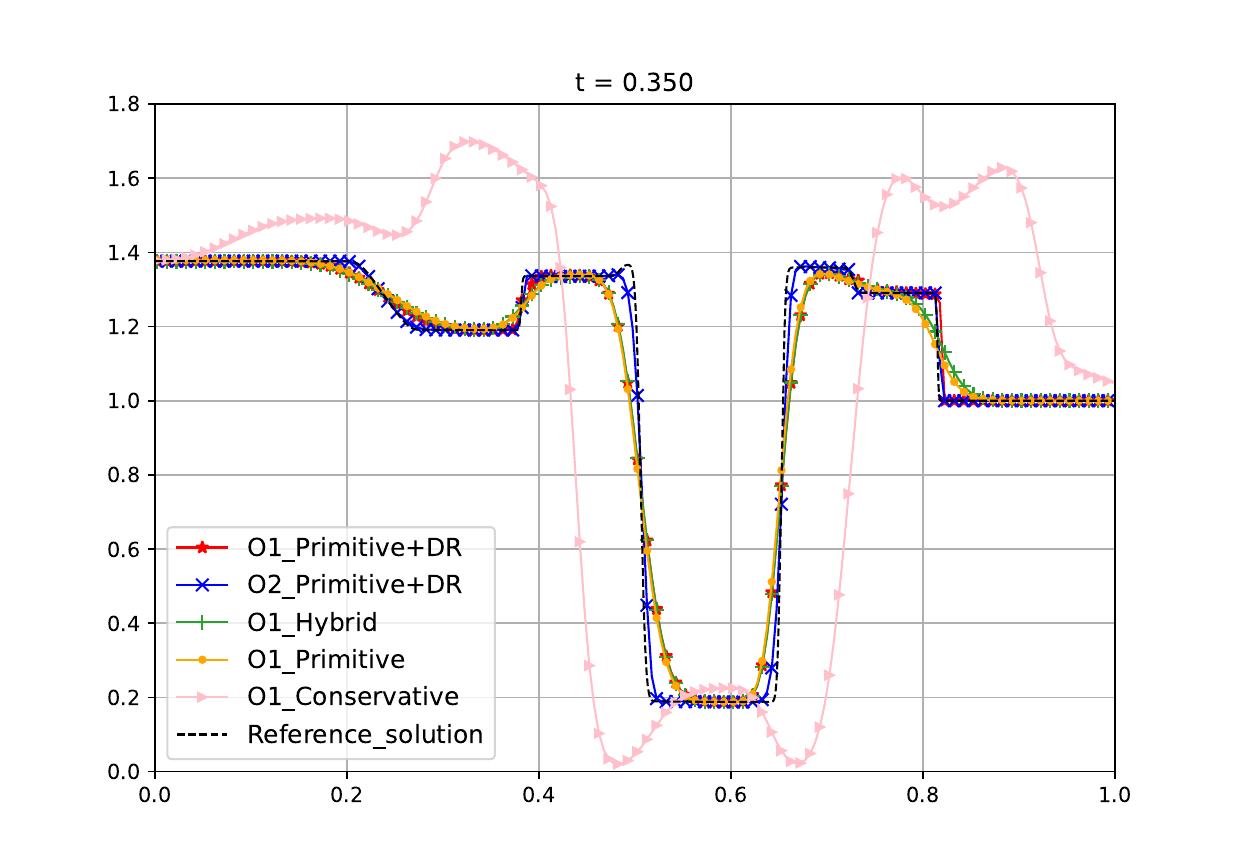}
    }
    \subfloat[Variable $\rho$. Zoom]{
        \includegraphics[width=0.23\textwidth]{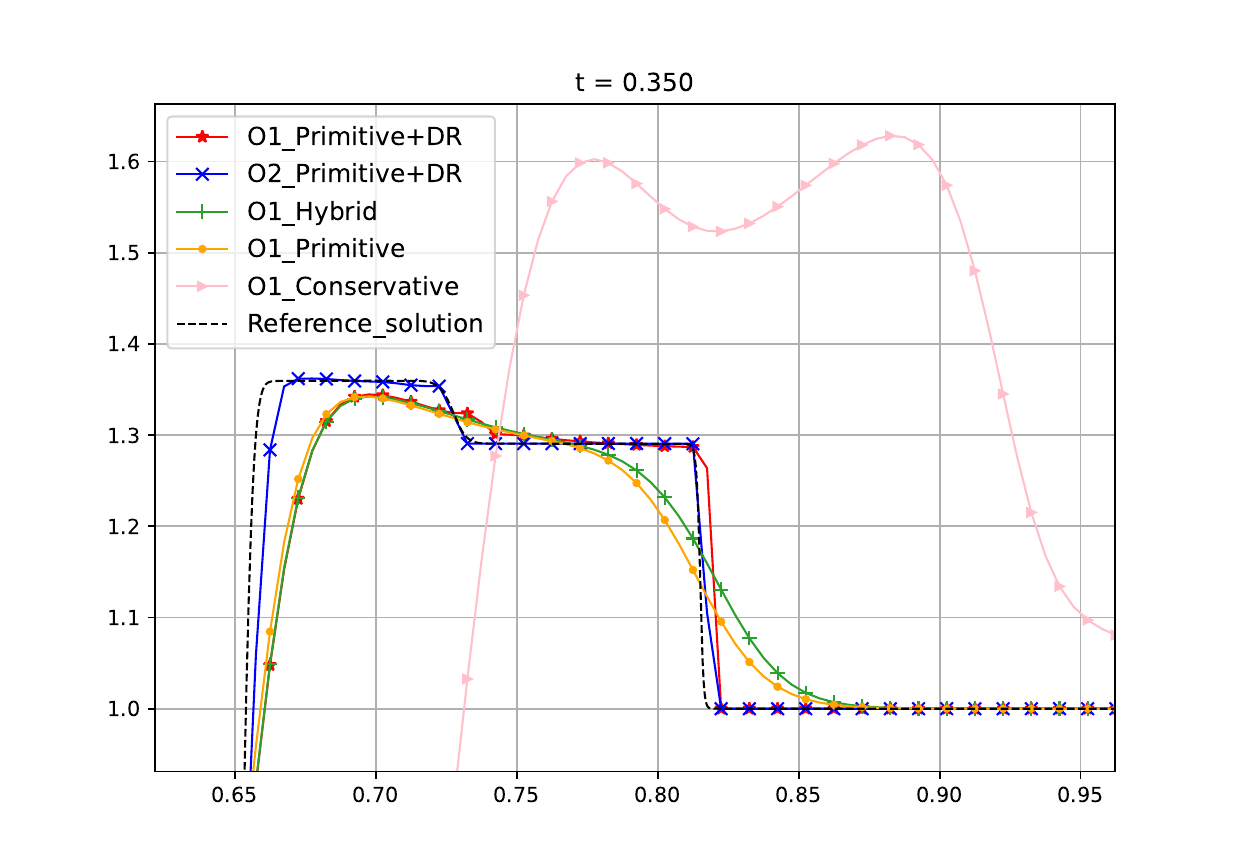}
    }

    \subfloat[Variable $u$]{
        \includegraphics[width=0.23\textwidth]{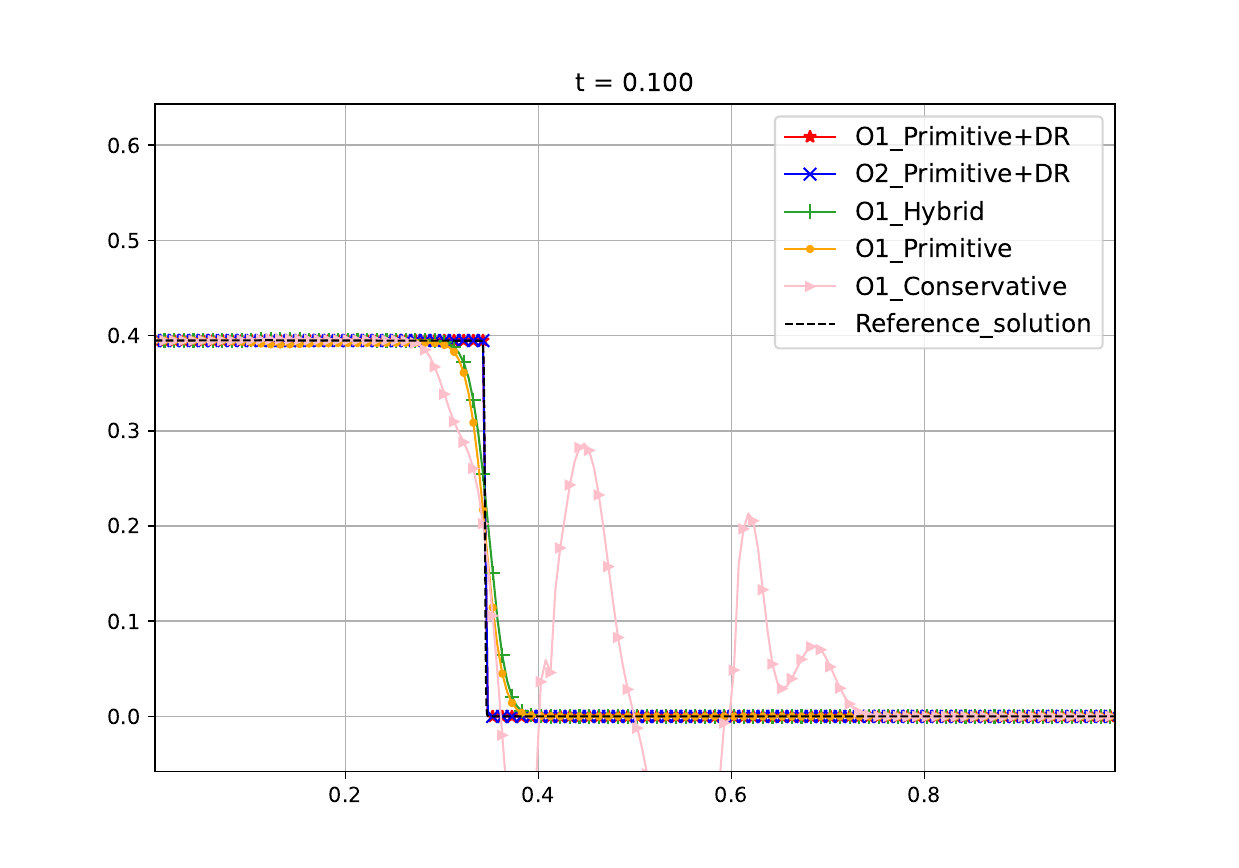}
    }
    \subfloat[Variable $u$. Zoom]{
        \includegraphics[width=0.23\textwidth]{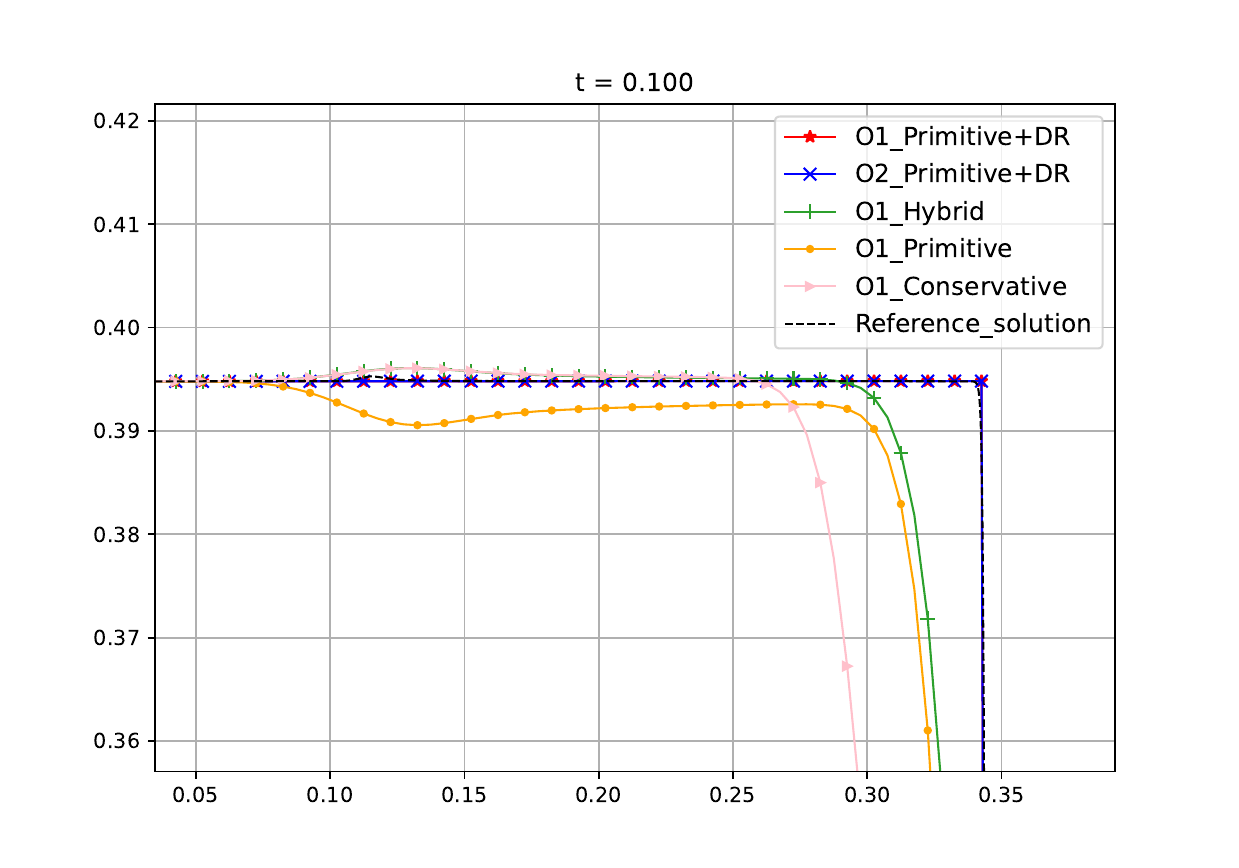}
    }
    \subfloat[Variable $u$]{
        \includegraphics[width=0.23\textwidth]{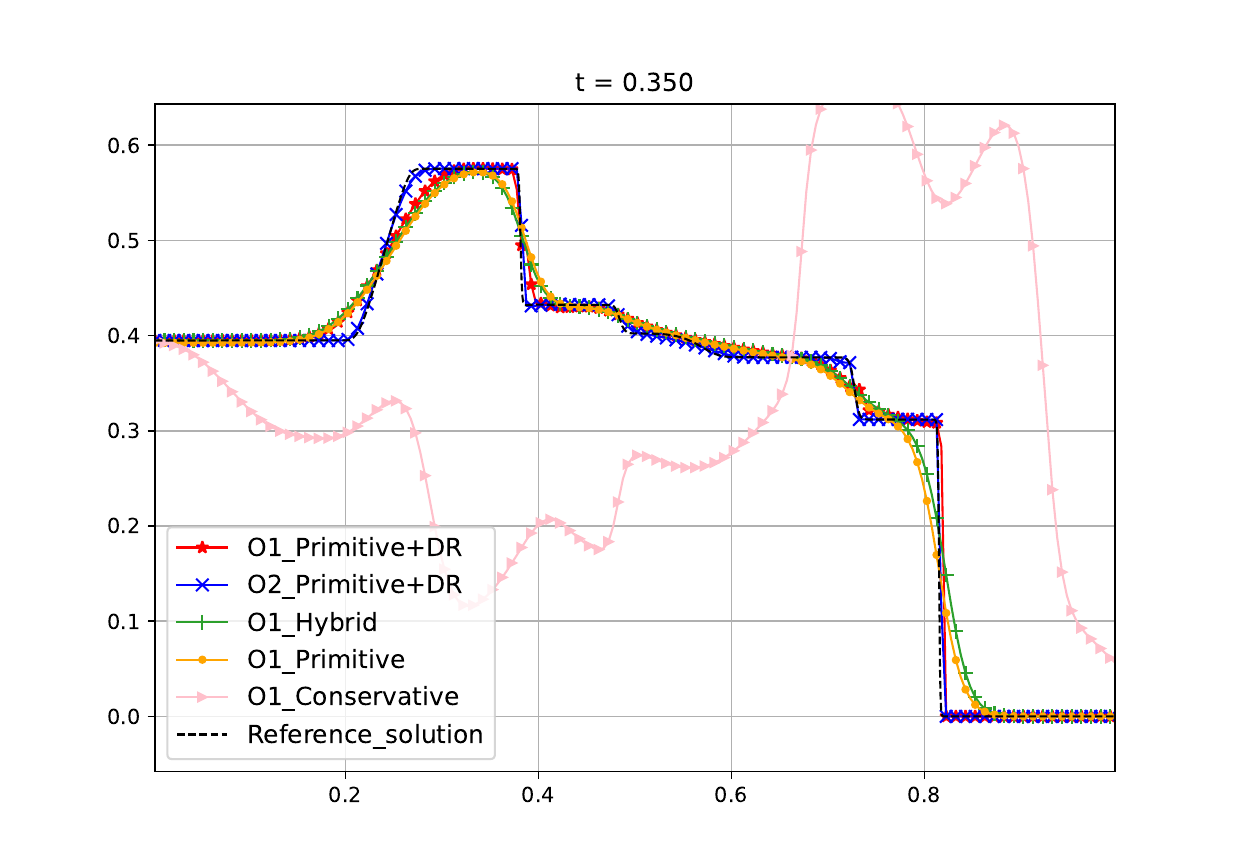}
    }
    \subfloat[Variable $u$. Zoom]{
        \includegraphics[width=0.23\textwidth]{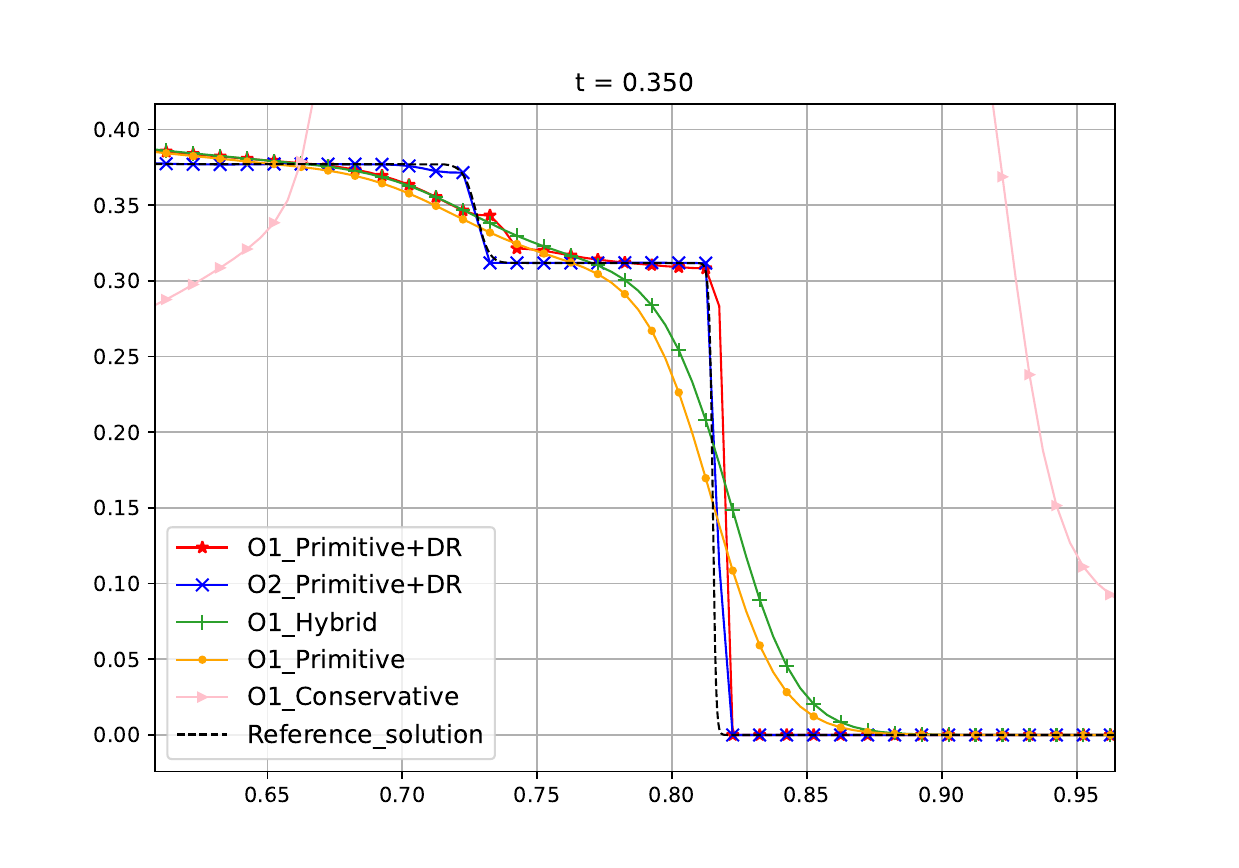}
    }

    \subfloat[Variable $p$]{
        \includegraphics[width=0.23\textwidth]{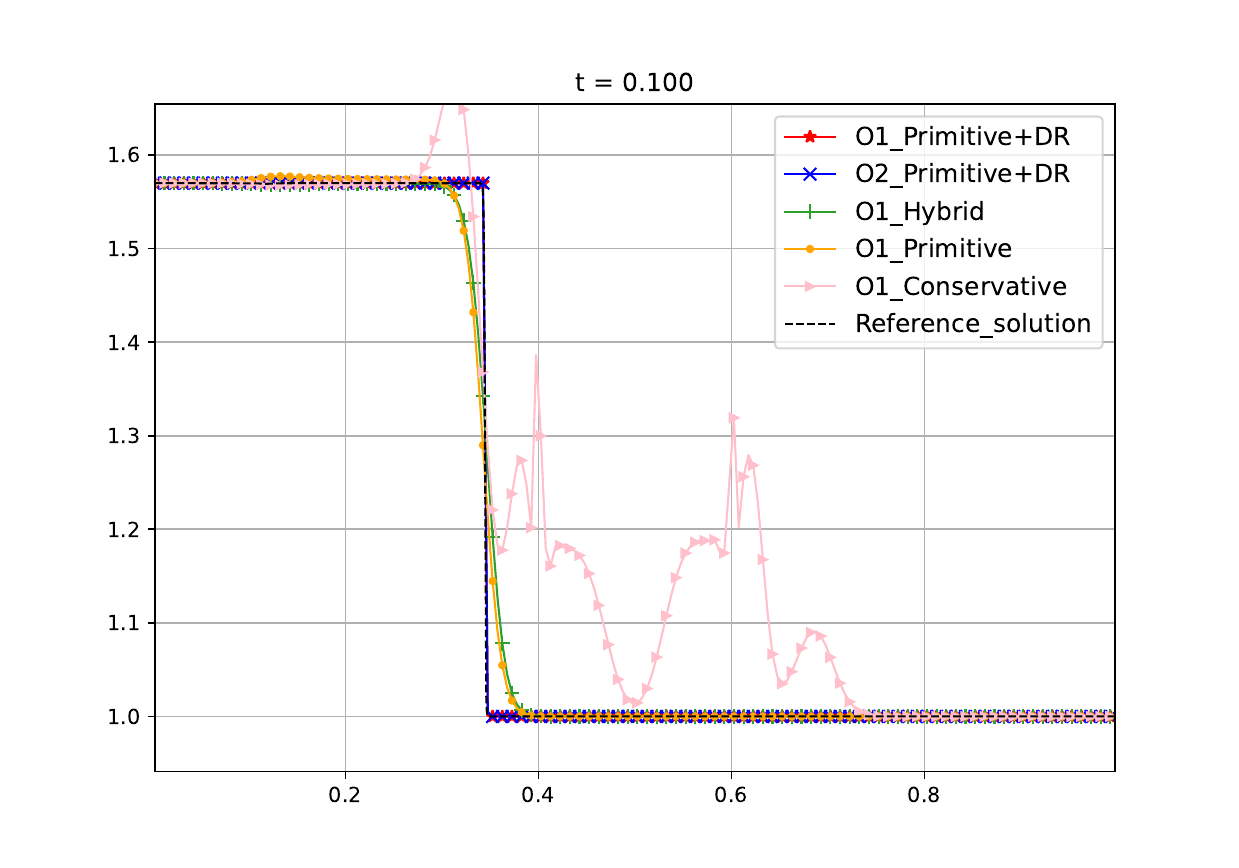}
    }
    \subfloat[Variable $p$. Zoom]{
        \includegraphics[width=0.23\textwidth]{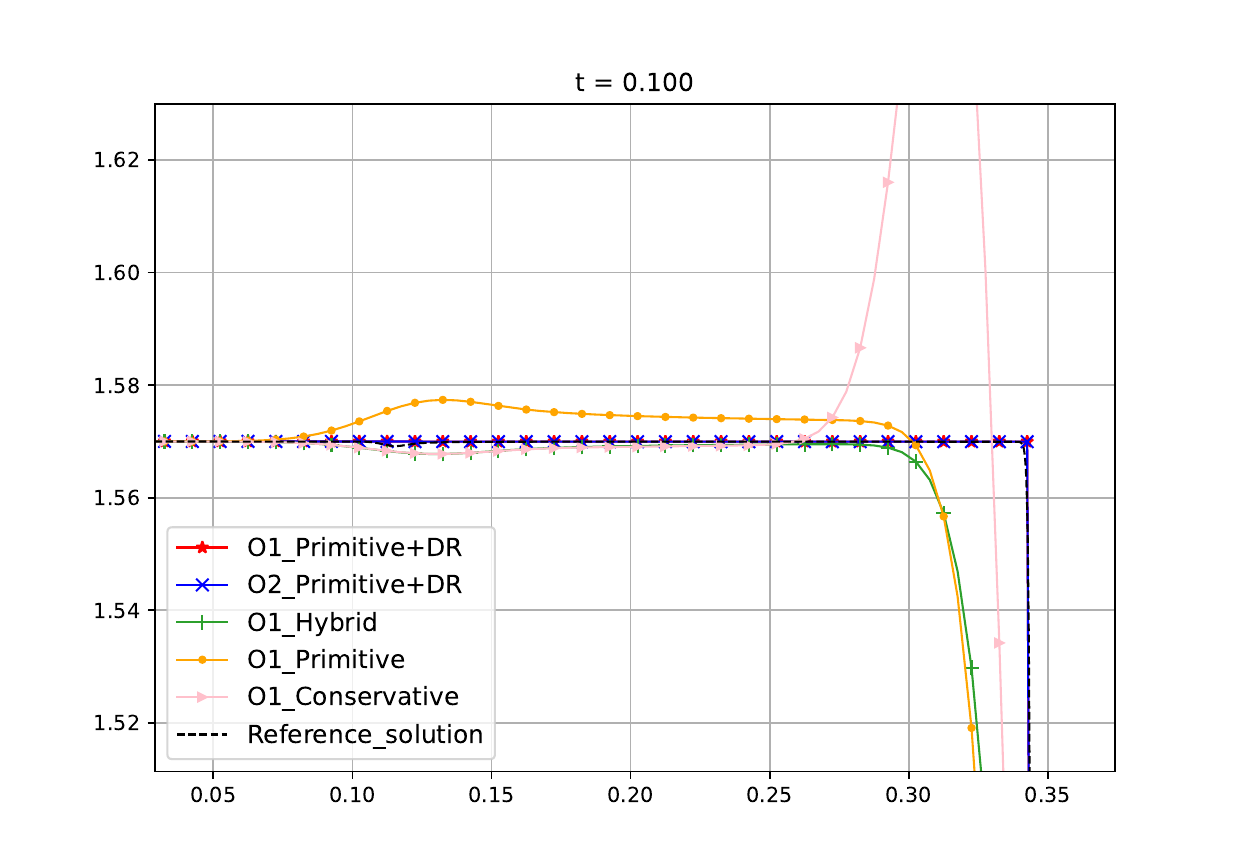}
    }
    \subfloat[Variable $p$]{
        \includegraphics[width=0.23\textwidth]{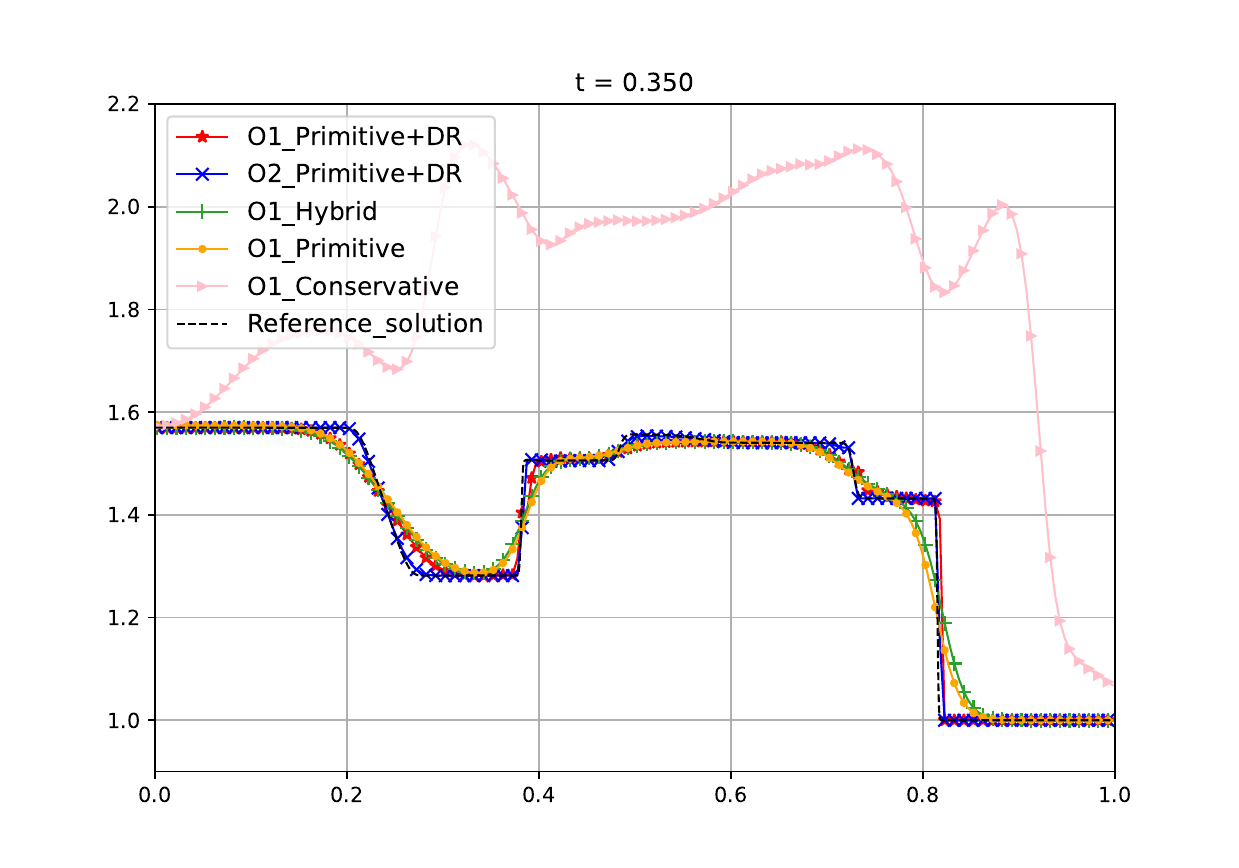}
    }
    \subfloat[Variable $p$. Zoom]{
        \includegraphics[width=0.23\textwidth]{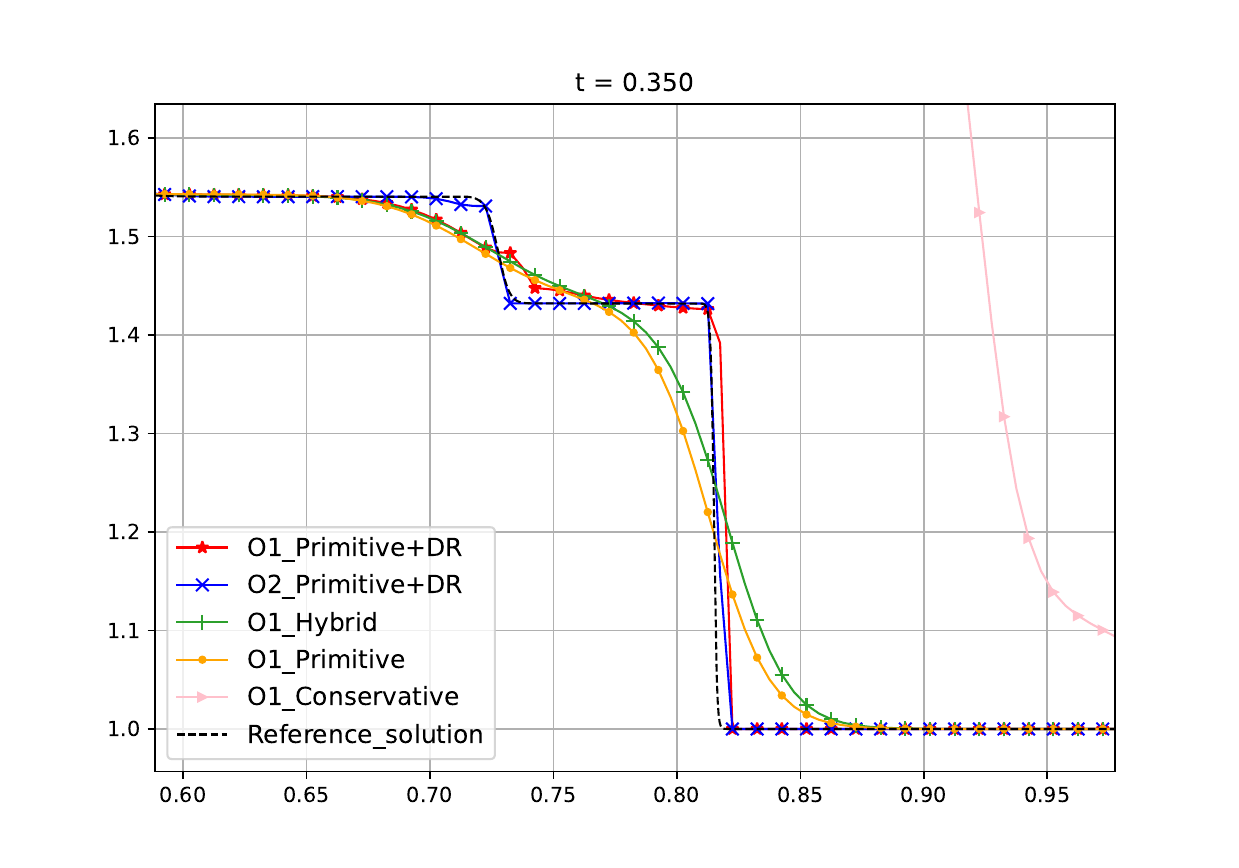}
    }
    \\
    \subfloat[Variable $v$]{
        \includegraphics[width=0.30\textwidth]{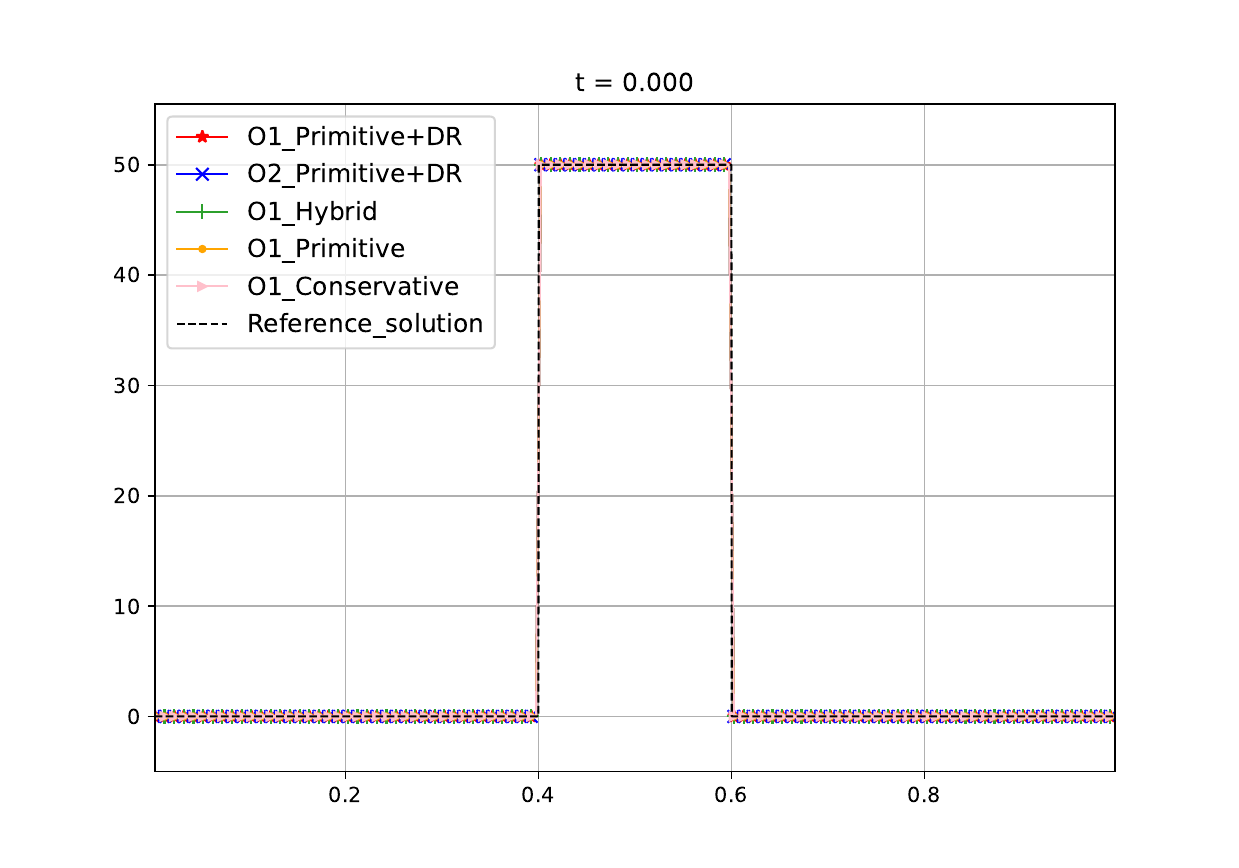}
    }
    \subfloat[Variable $v$]{
        \includegraphics[width=0.30\textwidth]{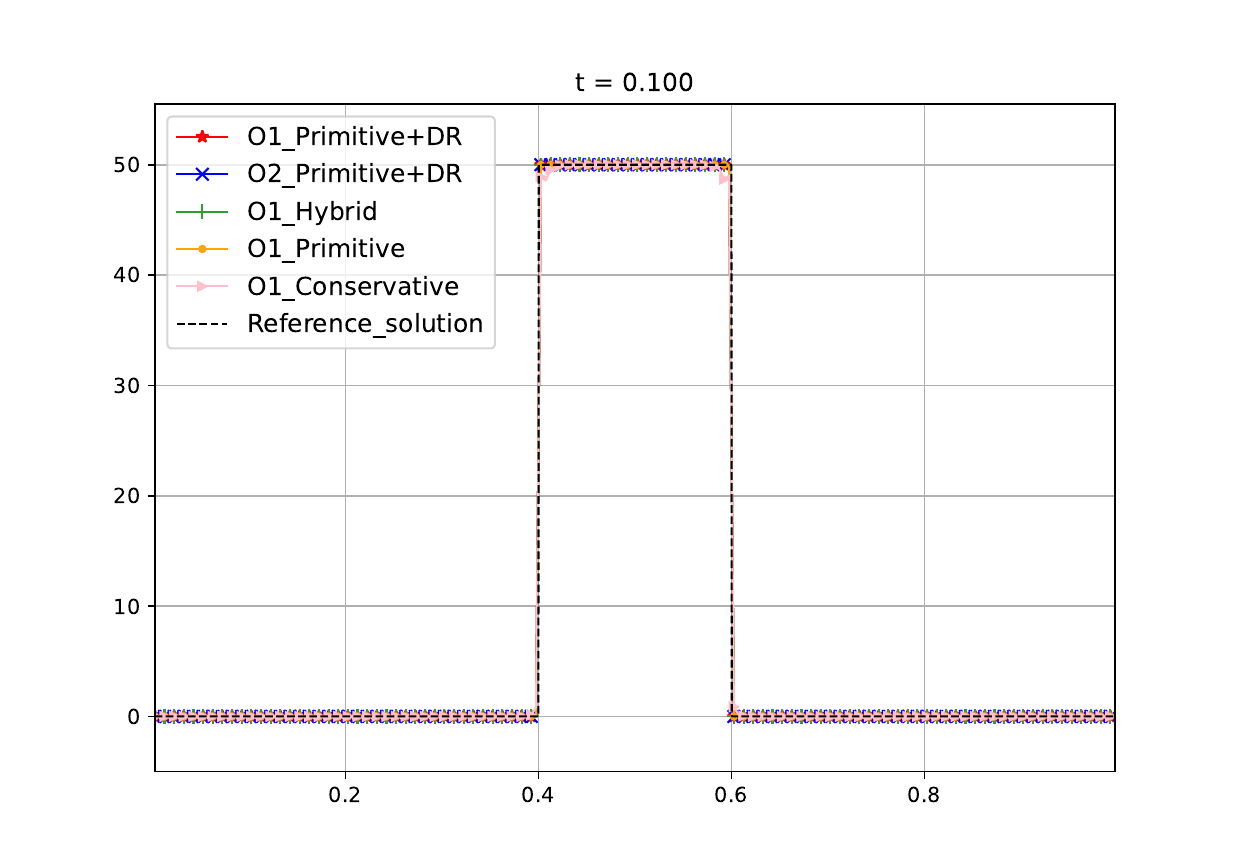}
    }
    \subfloat[Variable $v$]{
        \includegraphics[width=0.30\textwidth]{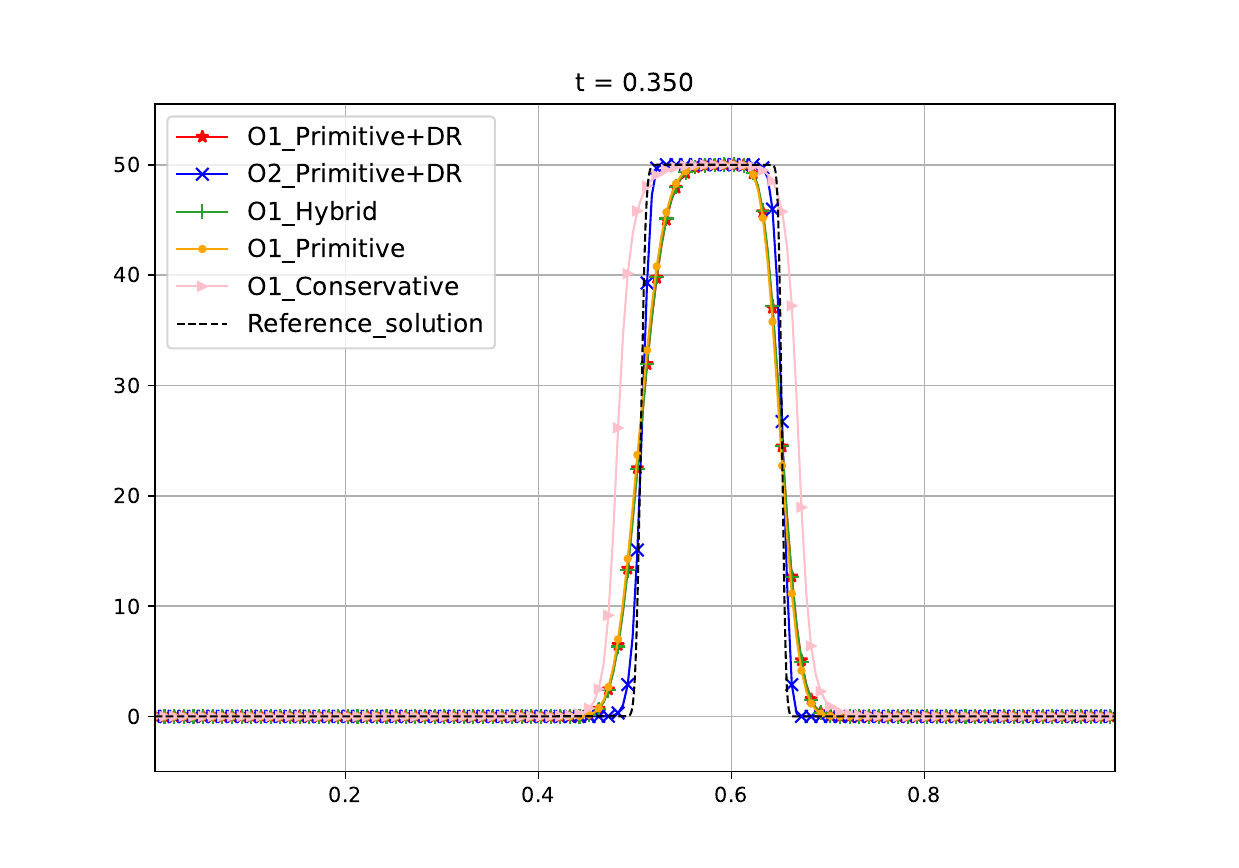}
    }

    \caption{Numerical results for Test: modified Quirk and Karni strong shock-shear interaction.}
    \label{Test_Big_Multispeed_Karni}
\end{figure}

This test is a modification of the previous one but now the shear wave in greater. The initial condition is given by
 $$
 \begin{aligned}
 \bw_{0}(x) & = (\rho_0(x), u_0(x), v_0(x), p_0(x), c_0(x))^{T} \\ & = \begin{cases}
     (1.3765, 0.3948, 0, 1.57,1)^{T} & \text{if $x<0.2$,}\\
     (1, 0, 0, 1,1)^{T} & \text{if $0.2\leq x<0.4$,}\\
     (0.138, 0, 50, 1,1)^{T} & \text{if $0.4\leq x\leq 0.6$,}\\
     (1, 0, 0, 1,1)^{T} & \text{if $x>0.6$.}
 \end{cases}
 \end{aligned}
 $$
Figure \ref{Test_Big_Multispeed_Karni} shows the numerical results. In this case, the strong contact discontinuities associated with the shear interface generate spurious oscillations in the fully conservative scheme, which completely ruins the numerical solution. For the remaining methods, similar conclusions can be drawn: standard path-conservative schemes fail to capture the correct shock waves, whereas the hybrid and in-cell discontinuous methods do. Among these, the schemes incorporating in-cell discontinuous reconstruction provide a more accurate resolution of shock waves. Once again, the second-order scheme exhibits superior performance in rarefaction regions.

\section{Conclusions}\label{sec:5}

In this work, we propose a methodology based on \cite{chalons2019path, pimentel2022cell, pimentel2024high}, for the numerical treatment of conservative systems formulated exclusively in primitive variables. The convergence difficulties associated with nonconservative products are addressed by identifying an appropriate family of paths and by incorporating in-cell discontinuous reconstructions within a path-conservative framework. In smooth regions or in the presence of contact discontinuities, a standard path-conservative method is used, while in regions where shocks appear, the in-cell discontinuous reconstruction is activated. The proposed approach is further extended to achieve high-order accuracy. Several numerical tests for the Euler multimaterial and shear interface demonstrate the improvement in the numerical solution: when using primitive variables exclusively together with this methodology, the scheme converges to the correct solution even in the presence of shocks, and spurious oscillations are effectively suppressed when dealing with contact discontinuities. No hybrid approaches are necessary to deal with these problems.

\appendix

\section{Paths and viscous profiles} \label{A1}
Although in general the choice of a family of paths for \eqref{nonconservative_system} that ensures that the jump condition \eqref{gen_R-H} is consistent with the physics of the problem may be difficult, when \eqref{nonconservative_system} comes from a conservative system \eqref{conservative_system} through a change of variables, there is a clear criterion for deciding which is the correct family of paths: it has to be such that a triple  $(\bw^- = T^{-1}(\bu^-), \bw^+ = T^{-1}(\bu^+), \sigma)$ satisfies \eqref{gen_R-H} if and only if  $(\bu^-, \bu^+, \sigma)$ satisfies \eqref{con_R-H}.
The interesting question is whether there is always a family of paths with this property. Let us assume that a triple $(\bu^-, \bu^+, \sigma)$ is an admissible discontinuity for \eqref{conservative_system} if $\bu^-$ and $\bu^+$ can be linked by a traveling wave of the viscous problem
\begin{equation}\label{viscous_problem}
    \bu_{t} + F(\bu)_{x} = \varepsilon \bu_{xx},
\end{equation}
i.e., if there exists a smooth function $V(\xi)$ satisfying
$$
\lim_{\xi \to \pm \infty} V(\xi) = \bu^{\pm}, \quad \lim_{\xi \to \pm \infty} V'(\xi) = 0,$$
such that
$$
\bu^\varepsilon (x,t) = V\left( \frac{x -\sigma t}{\varepsilon}\right)$$
is a solution of \eqref{viscous_problem} for every $\varepsilon$. The function $V$, which is called the viscous profile, has to solve the following boundary value problem:
\begin{equation}
\left \{
\begin{array}{l}
-\sigma V' + F(V)' = V'', \\
\displaystyle \lim_{\xi \to \pm \infty} V(\xi) = \bu^{\pm}, \quad \lim_{\xi \to \pm \infty} V'(\xi) = 0.
\end{array}
\right.
\end{equation}
In this case, the   path given by
$$
\Phi(\xi) = T^{-1}(V(\xi)), \quad \xi \in \mathbb{R}$$
is a compatible path, in the sense that it links
 $\bw^-$ and $\bw^+$ and satisfies \eqref{gen_R-H}. Indeed, one has:
\begin{eqnarray*}
0 & = & \frac{\partial \bw}{\partial \bu}\left( -\sigma V' + F(V)' - V'' \right) \\
& = & - \sigma \Phi' + \frac{\partial \bw}{\partial \bu}\frac{\partial F}{\partial \bu}\frac{\partial \bu}{\partial \bw} \Phi' - \frac{\partial \bw}{\partial \bu} V''\\
& = & - \sigma \Phi' + \mathcal{A}(\Phi) \Phi' - \frac{d\ }{d\xi} T^{-1}\left( V'\right).
\end{eqnarray*}
Integrating this equality in $\mathbb{R}$ we obtain:
\begin{eqnarray*}
\int_{-\infty}^{\infty} \mathcal{A}(\Phi(\xi)) \Phi'(\xi) \, d\xi & = & \sigma(\bw^+ -\bw^-) + T^{-1}(0) - T^{-1}(0) \\
 & = & \sigma(\bw^+ - \bw^-) ,
\end{eqnarray*}
as we wanted to prove. In order to obtain a family of paths with the sought property, first the paths found for states that can be linked by an admissible discontinuity have to be re-parametrized in the interval $[0,1]$, and then the family of paths have to be extended to arbitrary pair of states: this can be done, at least locally, on the basis of the solutions of Riemann problems as in \cite{MunozPares2007}. 

\section{Choice of the family of paths}\label{appendix_proof}
\paragraph{Proof of Theorem \ref{th:path_Euler}}
Let us see which choice of family of paths we should take to obtain  equivalence of \eqref{RK_conEuler} and \eqref{RK_nonconEuler}. We will focus on the second and third equations since the first one is the same in both cases.
Let us suppose that  \eqref{RK_nonconEuler} is equivalent to \eqref{RK_conEuler} for some family of paths $\Phi$. From \eqref{RK_conEuler} we have:
\begin{eqnarray*}
    \sigma(\rho_r u_r - \rho_l u_l)  =  \rho_r u_r^2 + p_r - \rho_l u_l^2 - p_l & \iff & \\\sigma(\rho_r u_r -\rho_l u_r + \rho_l u_r- \rho_l u_l) = \rho_r u_r^2 + p_r - \rho_l u_l^2 - p_l & \iff & \\
    u_r \sigma (\rho_r -\rho_l) + \rho_l\sigma(u_r- u_l) = \rho_r u_r^2 + p_r - \rho_l u_l^2 - p_l & \underset{\text{Using \eqref{RK_nonconEuler}}}{\iff} & \\
    u_r (\rho_r u_r - \rho_l u_l) + \rho_l\left(\frac{u_r^2-u_l^2}{2}\right) + \rho_l \int_0^1 \frac{\Phi_p'(s)}{\Phi_\rho(s)} ds= \rho_r u_r^2 + p_r - \rho_l u_l^2 - p_l & \iff & \\
    \int_0^1 \frac{\Phi_p'(s)}{\Phi_\rho(s)} ds= \frac{(u_r-u_l)^2}{2} + \frac{p_r - p_l}{\rho_l}.
\end{eqnarray*}
Now dividing the second equation of \eqref{RK_conEuler} by the first one, after some straightforward calculations, we obtain the equality
$$(u_r - u_l)^2 = (p_r-p_l)\left(\frac{1}{\rho_l} - \frac{1}{\rho_r}\right).$$
And using this equality, we finally obtain that the family of paths has to satisfy

$$
\begin{aligned}
\int_0^1 \frac{\Phi_p'(s)}{\Phi_\rho(s)} ds & = \frac{p_r-p_l}{2}\left(\frac{1}{\rho_l} - \frac{1}{\rho_r}\right)  + \frac{p_r - p_l}{\rho_r} \\ \iff \int_0^1 \frac{\Phi_p'(s)}{\Phi_\rho(s)} ds & = \frac{1}{2}\left(\frac{1}{\rho_l} + \frac{1}{\rho_r}\right)(p_r-p_l).
\end{aligned}
$$
A family of paths that satisfies this equality is given by
\begin{equation}
    \Phi_\rho(s) = \frac{1}{\frac{1}{\rho_l}+ s\left(\frac{1}{\rho_r}-\frac{1}{\rho_l}\right)}, \quad \Phi_p(s) = p_l + s (p_r-p_l).
\end{equation}
If we define
$$\Phi_u(s) =  u_l + s (u_r-u_l),$$
the equivalence of the third equations of \eqref{RK_conEuler} and \eqref{RK_nonconEuler} is proved using the same strategy of adding and subtracting $\rho_l\frac{u_r^2}{2}$.

\section*{Acknowledgments}
Ernesto Pimentel-García and Carlos Parés have been supported by projects PID2022-137637NB-C21 and PCI2024-155061-2 funded by \\MCIN/AEI/10.13039/501100011033/ and ERDF A way of making Europe. Both have also been supported by the Grant PPRO-FQM216-2023 funded by Consejería de Universidad, Investigación e Innovación and by ERDF Andalusia Program 2021-2027. MR is memebr of the CARDAMOM team at the INRIA research center at University of Bordeaux.

\bibliographystyle{siamplain}
\bibliography{references}
\end{document}